\documentclass[11pt, oneside]{article}

\RequirePackage[utf8]{inputenc}
\RequirePackage[T1]{fontenc}
\RequirePackage{xcolor,graphicx}
\RequirePackage{hyperref}
\hypersetup{%
  colorlinks=true,
  breaklinks=true,
  linkcolor=red,
  anchorcolor=black,
  citecolor=brown,
  filecolor=blue,
  menucolor=red,
  urlcolor=red,
}
\RequirePackage[authoryear]{natbib}
\RequirePackage{xr}
\RequirePackage{mathtools, amsfonts, amssymb,amsthm}
\mathtoolsset{showonlyrefs}
\RequirePackage{enumerate}
\RequirePackage{comment}
\RequirePackage{ifthen}
\RequirePackage{enumitem}
\RequirePackage{comment}
\RequirePackage{subcaption}
\RequirePackage{siunitx}
\RequirePackage[ruled,vlined]{algorithm2e}

\newtheorem{theorem}{Theorem}

\newtheorem{lemma}{Lemma}

\theoremstyle{remark}
\newtheorem{definition}{Definition}
\newtheorem{assumption}{Assumption}

\newcommand{\EE}{\mathbb{E}}
\newcommand{\EEi}{\mathbb{E}_i}

\newcommand{\EEMT}{\mathbb{E}_{M,T}}

\newcommand{\PP}{\mathbb{P}}

\newcommand{\NN}{\mathbb{N}}

\newcommand{\n}[1]{{n_{#1}}}

\makeatletter
\newcommand{\pushright}[1]{\ifmeasuring@#1\else\omit\hfill$\displaystyle#1$\fi\ignorespaces}
\newcommand{\pushleft}[1]{\ifmeasuring@#1\else\omit$\displaystyle#1$\hfill\fi\ignorespaces}
\makeatother

\newcommand{\degree}{\mathbf{d}}
\newcommand{\hdegree}{\hat {\mathbf{d}}}
\newcommand{\T}{{t}}

\newcommand{\Mmu}{\mathfrak{m}}

\newcommand{\Ostar}{\mathcal{O}_*}

\newcommand{\WN}{\mathcal W_N}

\newcommand{\Xtemp}[1]{%
  \ifthenelse{\equal{#1}{0}}{X^{(\n0)}}{
  \ifthenelse{\equal{#1}{1}}{X^{[\n1]}}{X^{(#1)}}
  }
}
\newcommand{\hXtemp}[1]{%
  \ifthenelse{\equal{#1}{0}}{\widehat X^{(\n0)}}{
  \ifthenelse{\equal{#1}{1}}{\widehat X^{[\n1]}}{\widehat X^{(#1)}}
  }
}
\newcommand{\Xp}[1]{X^{(#1)}}

\newcommand{\X}[1]{{\Xtemp{#1}}}

\newcommand{\Xn}{{\Xtemp{i}}}

\newcommand{\hXni}{{\hXtemp{i}}}

\newcommand{\Ynm}{{Y^{(i)}_m}}
\newcommand{\Tnm}{{T^{(i)}_m}}
\newcommand{\enm}{{\varepsilon^{(i)}_m}}
\newcommand{\unm}{{e^{(i)}_m}}

\newcounter{experiment}[section]

\title{Adaptive estimation of irregular mean and covariance functions}
\author{%
Steven Golovkine\footnote{MACSI, Department of Mathematics and Statistics, University of Limerick, Ireland; steven.golovkine@ul.ie}
\and 
Nicolas Klutchnikoff\footnote{Univ Rennes, CNRS, IRMAR - UMR 6625, F-35000 Rennes, France; nicolas.klutchnikoff@univ-rennes2.fr}
\and
Valentin Patilea\footnote{Univ Rennes, Ensai, CNRS, CREST - UMR 9194, F-35000 Rennes, France; valentin.patilea@ensai.fr}
}
\date{\today}

\begin{document}

\maketitle

\begin{abstract}
Nonparametric estimators for the mean and the covariance functions of functional data are proposed. The setup covers a wide range of practical situations. The random trajectories are, not necessarily differentiable, have unknown regularity, and are measured with error at discrete design points. The measurement error could be heteroscedastic. The design points could be either randomly drawn or common for all curves. The  estimators depend on the local regularity of the stochastic process generating the functional data. We consider a simple estimator of this local regularity which exploits  the replication and regularization features of functional data. Next, we use the ``smoothing first, then estimate'' approach for the mean and the covariance functions. 
They can be applied with both sparsely or densely sampled curves, are easy to calculate and to update, and perform well in simulations. Simulations built upon an example of real data set, illustrate the effectiveness of the new approach.


\end{abstract}

%%%%%%%%%%%%%%%%%%%%%%%%%%%%%%%%%%%%%%%%%%%%%%
%%%% Main text entry area:
% !TeX root=bj-sample.tex

\section{Introduction}

Motivated by a large number of applications, there is a great interest in models for observation entities in the form of a sequence of measurements recorded intermittently at several discrete points in time. Functional data analysis (FDA) considers such data as being values on the trajectories of a stochastic process, recorded with some error, at discrete random times. The mean and the covariance functions play a critical role in FDA. 

To formalize the framework, let $\mathcal T$ be a compact interval, typically $[0,1]$. Data consist of random realizations of sample paths from a second-order stochastic  process $X = (X_t : t\in \mathcal T)$ with continuous trajectories. The mean and covariance functions are $\mu(t) = \EE[X_t]$ and
$$
\Gamma (s,t) =  \EE\left[ \{X_s - \mu(s)\}  \{X_t-\mu(t)\}\right] = \EE\left[X_s X_t\right]-\mu(s)\mu(t) ,\quad s,t\in \mathcal T,
$$
respectively. If the independent realizations $\Xp{1},\dotsc,\Xp{i},\dotsc,\X{N}$ of $X$ were observed, the ideal estimators would be
$$
\widetilde \mu_N (t) = \frac{1}{N} \sum_{i=1}^N \Xp{i}_t \quad \text{ and } \quad 
\widetilde \Gamma _N(s,t) = \frac{1}{N-1} \sum_{i=1}^N \left\{\Xp{i}_s -\widetilde \mu_N(s) \right\}\left\{\Xp{i}_t -\widetilde \mu_N(t)\right\} ,\quad s,t\in \mathcal T.
$$

In applications, the curves are rarely observed without error and never at each value $t\in\mathcal T$. This is why we consider the following common and more realistic setup. For each $1\leq i \leq N$, and given a positive integer $M_i$, let $\Tnm\in\mathcal T$, $1\leq m \leq  M_i$, be the %possibly random 
observation times for the curve $\Xp{i}$. The observations associated with a curve, or trajectory, $\Xp{i}$ consist of  the pairs  $(\Ynm , \Tnm ) \in\mathbb R \times \mathcal T $ where %$\Ynm $ is defined as
\begin{equation}\label{model_eq}
  \Ynm = \Xn(\Tnm) + \enm,
  \quad 1\leq m \leq M_i ,\quad 1\leq i \leq N,
\end{equation}
and  $\enm$  is an independent (centered) error variable. Here, and in the following, we use the  notation $\Xp{i}_t$ for the value at a generic point $t\in\mathcal T$ of the realization  $\Xp{i}$ of $X$,  while $\Xp{i}(\Tnm)$ denotes the measurement  at $\Tnm$ of this realization.

A commonly used idea is to build feasible versions of $\widetilde \mu_N (\cdot) $ and $\widetilde \Gamma _N(\cdot,\cdot)$ using nonparametric estimates of $\Xp{i}_t$ and $\Xp{i}_s\Xp{i}_t$, such as those obtained by smoothing splines or local polynomials. This approach, usually called ``smoothing first, then estimate''  or ``two-stage procedure'',  is considered, among others, by \cite{hall2006} and \cite{zhang_statistical_2007}. In general, the sample trajectories are required to admit at least second-order derivatives over $\mathcal T$. \cite{li2010}, \cite{zhang_sparse_2016} and \cite{zhang_wang2018} propose an alternative local linear smoothing approach where the  estimators are determined by suitable weighting schemes which involve the whole sample of curves. This idea exploits the so-called replication and regularization features of functional data, see \citet[Ch.~22]{ramsay_functional_2005}. In this alternative approach, the regularity assumptions are imposed on the mean and covariance functions, which are required to admit second, or higher, order derivatives over the  domain. Since, in general, the mean and covariance functions are more regular than the sample trajectories, the approach based on  weighting schemes using all the sample curves might be preferable. However, in some cases, for instance in energy, chemistry, physics, astronomy or medical applications, the mean and covariance  could be quite irregular, of unknown irregularity. 

\cite{cai2011} and \cite{cai2010} derive the optimal rates of convergence, in the minimax sense, for the mean and covariance functions, respectively, and proposed optimal estimators. The  estimator of the mean function proposed by  \cite{cai2011} is a smoothing spline estimator which could be built only if the regularity  of the sample paths  is known.  \cite{cai2010} use the representation of the covariance function in a tensor product reproducing kernel Hilbert  space. Under some assumptions, they derive estimators for $\Gamma (s,t)$ using a low-dimensional version of this representation obtained by a regularization procedure, provided the values $M_i$ are not very different. This procedure involves  numerical optimization. 
See also \cite{wong2019}. The  optimal rates for the mean and covariance functions are defined by the sum of two types of terms. One corresponds to the rate of convergence of the ideal estimators $\widetilde \mu_N (\cdot) $ and $\widetilde \Gamma _N(\cdot,\cdot)$, which is the standard rate of convergence for empirical means and covariances. The other contribution to the optimal rates is given by the differences between $\widetilde \mu_N (\cdot) $ and $\widetilde \Gamma _N(\cdot,\cdot)$ and their feasible versions. The optimal rates of the differences depend on the regularity of sample trajectories, because the minimax lower bounds should also take into account the case where the functions to be estimated have the same regularity as the trajectories. 

The estimation of the mean and covariance functions presents another specific feature. The  optimal  rates of convergence depend on the nature of the measurement times  $\Tnm$. Up to now, two situations have been  investigated in the literature. On the one hand, the so-called \emph{independent design} case where, given the $M_i$'s, the $\Tnm$ are obtained as a random sample of size $M_1+\cdots+M_N$ from the same continuous distribution. On the other hand,  the so-called \emph{common design} case where the $M_i$ are all equal to some integer value $\Mmu$, and the $\Tnm$, $1\leq m\leq  \Mmu$, are the same across the curves $\Xp{i}$. In both cases, the best rates for the nonparametric estimators depend on the regularity of the sample trajectories. These rates also depend on the number of different observation times $\Tnm$, that is equal to $M_1+\cdots+M_N$ with independent design, and equal to $\Mmu$ with common design. In other words, the replication feature of functional data is less impactful with common design. See \cite{cai2011} for the case of the mean function, and \cite{cai2010} and \cite{cai2016minimax} for the case of the covariance function. 

In this paper, we propose   data-driven   ``smoothing first, then estimate'' type methods, based on  $1-$dimensional smoothing. 
%, which  achieve minimax optimal rates of convergence. 
The  process is  allowed to have a varying, unknown regularity. 
Our method does not require complex numerical optimization. It applies in the same way to common and independent design situations, and allows for general  heteroscedastic measurement errors $\enm$. Moreover, our approach is suitable with both sparsely or densely sampled curves. The definition of sparse and dense regimes is recalled in Section \ref{sec:methodo}. 

Let $ \hXni$ be a suitable  nonparametric estimator of   $\Xp{i}$ applied to the $M_i$ pairs $(\Ynm , \Tnm )$, for instance a  kernel  estimator.  This estimator will be suitable because it will take into account the regularity of the  process $X$ and the final estimation purpose, i.e. the mean or the covariance function. 
These features can be achieved in an easy, data-driven way, as explained below.
With at hand, the $ \hXni$'s tuned for the mean function estimation, we define
\begin{equation}\label{est_mean0}
\widehat \mu_N (t) = \frac{1}{N} \sum_{i=1}^N \hXni_t,\quad t\in \mathcal T. 
\end{equation}
For the covariance function, we distinguish the diagonal from the non-diagonal points. With at hand, the $ \hXni$'s tuned for the covariance function estimation, and for some diagonal set $\mathcal D \subset \mathcal T^2\coloneqq   \mathcal T\times \mathcal T $ that we shall determine using the data, let us define 
\begin{equation}\label{est_cov1}
\widehat \Gamma _N(s,t) = \frac{1}{N} \sum_{i=1}^N \hXni_s \hXni_t - \widehat \mu_N(s)\widehat \mu_N(t) ,\quad (s,t)\in \mathcal T^2 \setminus \mathcal D.
\end{equation}
It is well known that the variance function $\Gamma(s,s)$ induces a singularity when estimating the covariance function $\Gamma(\cdot,\cdot)$. See, for instance, \citet[Rem.~4]{zhang_sparse_2016}. 
We propose a simple way to build the diagonal set $\mathcal D$, which asymptotically reduces to the diagonal segment according to a data-driven rule that we provide in the following. Given  $\mathcal D$, the estimates of $\Gamma(\cdot,\cdot)$ on $\mathcal D$ are directly obtained from the estimates $\widehat \Gamma _N(s,t)$ for the closest $(s,t)$ on the boundary of $\mathcal D$. 

Although  the methodology we propose is general and can be used with different types of smoothers, we focus on the case where the $\hXni_t$ are obtained by kernel smoothing. In this case, tuning the $ \hXni$'s means  suitably determining the rate of decrease and the constant defining the bandwidth. Our approach is completely data-driven and relies on the minimization of a new and suitable risk function.

To the best of our knowledge, there is no contribution  
which considers estimators of the curves $\Xp{i}$ adapted to their regularity and to the purpose of estimating mean or covariance functions.  Trajectory-by-trajectory adaptive optimal smoothing, for instance using 
the \cite{GL2011} method, generally yields sub-optimal rates of convergence for $\widehat \mu_N (t) $ and $\widehat \Gamma _N(s,t) $. The reason is that trajectory-by-trajectory smoothing ignores the information contained in the other $N-1$ curves in the sample, generated according to the same stochastic process $X$. See  \cite{cai2011} for a discussion on the differences with the usual nonparametric rates.
One can also use cross-validation for choosing the bandwidth with the suitable weighting schemes, such as proposed by \cite{li2010} or \cite{zhang_sparse_2016}. This would however require significant computational effort, and, to the best of our knowledge, the idea has not yet received  a theoretical justification.  Using the replication and regularization features of functional data, we consider an effective estimator for  the local regularity of the process $X$, a probabilistic concept which determines the analytic regularity of the trajectories of $X$. 
The local regularity estimator,  a version of the one introduced by \cite{golo2020},  combines information both across and within curves.
Moreover, it allows for general heteroscedastic measurement errors, does not involve any optimization and is obtained after a fast, possibly parallel, computation. With at hand the local regularity estimator, we derive the  
 suitable estimators $\hXni_t$, 
and finally our optimal mean and covariance functions estimators. 
The smoothing parameter used to build  the $\hXni_t$ depends on  $M_i$ and $N$,  but can be easily computed given the  estimate of the local regularity of $X$. We assert that the replication feature of the functional data makes the concept of local regularity of the process a more meaningful parameter than the usual curve regularity,  which is an analytic concept designed for a single function.

In Section \ref{sec:methodo}, we provide insight on why the local regularity of the process $X$ is a natural feature to be considered. Moreover, we explain why the  ``smoothing first, then estimate'' approach could achieve optimal rates  when the  regularity of $X$ is known. In Section \ref{subsec:golo_new}, we formally define the local regularity of the process $X$.

Moreover,  we introduce the  estimator for this regularity and present exponential bounds for the concentration under mild conditions. In particular, both independent and common designs are allowed, and the process regularity  is allowed to vary with $t$. Section \ref{subsec:golo_new} ends with a  discussion on the relationship between the process regularity  and the analytical regularity of the trajectories.  In Section \ref{sec:LP_est}, we use the regularity estimate to build sharp bounds of the pointwise quadratic risk function between our estimators and the unfeasible estimators  $\widetilde\mu_N$ and $\widetilde \Gamma_N$, respectively. The bounds depend on quantities which could be estimated by sample averages. Minimizing the risk bounds with respect to the bandwidth, we derive the optimal bandwidth for the  kernel estimates of the trajectories. These estimates are  further used to  estimate the mean and covariance functions. 
Our mean and covariance estimators, and the local regularity estimator, are computed on the same sample of curves. In other words, no data splitting is necessary with our approach. The finite sample performance of the new estimators is illustrated in Section \ref{sec:emp} using simulated samples generated according to the setup of a real data set on the power consumption of households. The simulation method which we introduce in Section \ref{sec:emp} is a simple device allowing functional data to be generated with regularity features similar to those observed in real applications. 
Some conclusions and discussions are given in  Section \ref{sec:dics}. A few proofs are presented in the Appendix. The Supplement \cite{supplement} contains more  technical arguments and simulation results.

% !TeX root=bj-sample.tex

\section{From unfeasible to feasible optimal  estimators}\label{sec:methodo}

The novelty of our approach is based on the local regularity of $X$, a mild condition on the second-order moments of the local increments of the process $X$. Before giving formal definitions, we first explain why the 
local regularity of the process generating the curves is a meaningful concept,  and why our approach can achieve good performance.
For this purpose, we analyze the difference $\widehat \mu_N (t) -\widetilde \mu_N (t)$, $s,t\in \mathcal T,$ but similar ideas apply to the covariance function estimation. 

The data $(\Ynm , \Tnm ) \in\mathbb R \times \mathcal T $ are generated according to \eqref{model_eq} with
\begin{equation}\label{def_err}
\enm = \sigma(\Tnm,\Xn(\Tnm)) \unm   , \quad 1\leq m \leq  M_i,\quad 1\leq i \leq N,
\end{equation}
where  the $\Xn$ are independent trajectories of $X$,   $\unm$ are independent copies of a centered variable $e$ with unit variance, and $\sigma(t,x)$ is some unknown    bounded  function which takes into account possible heteroscedastic measurement errors. 
The integers $M_1, \dotsc,M_N$ represent an independent sample of an integer-valued random variable $M$ with expectation $\Mmu$ which increases with $N$. Thus, $M_1, \dotsc,M_N$ is the $N$th line of a triangular array of integers. 
In the independent design case, for each $1\leq i \leq N$, the  observation  times $ \Tnm $ are  random copies of a  variable $T\in \mathcal T$. The realizations of $X$, $e$, $M$ and $T$ are assumed to be mutually independent.   Let $\mathcal T_{obs}^{(i)} $ be the set of observation times $ \Tnm $, $1\leq m \leq M_i$, over the trajectory $\Xp{i}$.  
With a common design, $M\equiv \mathfrak m$, and the $\mathcal T_{obs}^{(i)} $  are the same for all $i$.
If not stated differently, the issues discussed in this section apply to both independent and common design cases. 

Let  
$$\EEi [\cdot] = \EE [\cdot \mid M_i,  \mathcal T_{obs} ^{(i)}  , \Xp{i}]  \quad \text{ and } \quad \EEMT[\cdot] = \EE [\cdot \mid M_i, \mathcal T_{obs}^{(i)} ,1\leq i \leq N ].$$
For any $t\in\mathcal T$, we consider a generic, linear, nonparametric  estimator:
\begin{equation}\label{LP_est_v1}
\hXni_t  = \sum_{m=1}^{M_i} \Ynm W_{m}^{(i)}(t), 
\quad 1\leq i \leq N.
\end{equation}
The weights $W_{m}^{(i)}(t)$ are defined as functions of the elements in $\mathcal T_{obs}^{(i)} $. 
The example we keep in mind, and investigate in detail in Section \ref{sec:LP_est}, is that of kernel smoothing  with a compactly supported kernel. 
Let
\begin{equation}\label{dec_s}
\hXni_t - \Xp{i}_t = B_t^{(i)} + V_t^{(i)} , \quad t \in \mathcal{T},
\end{equation}
where 
$$
B_t^{(i)} \coloneqq \EEi[\hXni_t ]  -  \Xp{i}_t  
\quad 
\text{ and } \quad
 V_t^{(i)} \coloneqq \hXni_t - \EEi [\hXni_t ] = \sum_{m=1}^{M_i} \enm  W_{m}^{(i)}(t).
$$
The pairs of random variables $(B_t^{(i)}, V_t^{(i)})$, $1\leq i\leq N$, are independent and we could reasonably assume that they are squared integrable for all $t$. For the mean, we can then write
\begin{equation}\label{dec_mean}
\widehat \mu_N(t) - \widetilde \mu _N(t) = \frac{1}{N} \sum_{i=1}^N  B_t^{(i)} + \frac{1}{N} \sum_{i=1}^N    V_t^{(i)}.
\end{equation}
All the variables  $ \enm  $ are centered and conditionally independent, with bounded conditional variance, given all $M_i$, $\mathcal T_{obs}^{(i)}$ and $\Xp{i}$. Thus, 
\begin{equation}\label{eq:var_8}
\EEMT\left[\left\{ N^{-1} \sum_{i=1}^N   V_t^{(i)} \right\}^2 \right]  
\leq N^{-1} \sup_{x} \sigma^2(t,x) \times N^{-1}  \sum_{i=1}^N \left\{ \max_{m} \left\lvert W_{m}^{(i)}(t)\right\rvert \times  \sum_{m=1}^{M_i} \left\lvert W_{m}^{(i)}(t)\right\rvert\right\}.  
\end{equation}
For the Nadaraya-Watson (NW) estimator with bandwidth $h$, under some mild conditions,  the rate of decrease of the right-hand side in the last display  is $O_{\mathbb P} ( (N\mathfrak m h)^{-1})$.

For illustration, we suppose here that the trajectories are not differentiable.  The  case of smooth  paths is discussed in Section \ref{sec:dics}. On the bias part, by Cauchy-Schwarz inequality, we then have   
\begin{equation}\label{eq:reta_bias2}
\EEMT\left[\left\{ N^{-1} \sum_{i=1}^N   B_t^{(i)} \right\}^2 \right] 
  \leq   N^{-1} \sum_{i=1}^N \left\{ \sum_{m=1}^{M_i} \left\lvert W_{m}^{(i)}(t)\right\rvert  \sum_{m=1}^{M_i}\EEMT \left[ \left\{ \Xn(\Tnm)- \Xp{i}_t \right\}^2  
  %\mid \mathcal  T_{obs}^{(i)}
   \right] \left\vert W_{m}^{(i)}(t)\right\rvert \right\} .
\end{equation}
It now becomes clear that the rate of the square of the bias term in $\widehat \mu_N(t) - \widetilde \mu _N(t)$  is determined  by the second-order moment of the increments $   \Xn(\Tnm)- \Xp{i}_t  $. If,  for $u,v \in \mathcal T$ close to $t$, 
\begin{equation}\label{eq:reg0}
\EE [ \left\{  X_{u}- X_{v} \right\}^2  ] \approx L_t^2 \lvert u-v \rvert^{2H_t},
\end{equation}
with  some  $0< H_t\leq 1$ and $L_t>0$, then  the rate of the right-hand side  in \eqref{eq:reta_bias2}  is bounded by
\begin{equation}\label{eq:reg01}
 N^{-1}  \sum_{i=1}^N  \left\{\sum_{m=1}^{M_i} \left\lvert W_{m}^{(i)}(t)\right\rvert \times  \sum_{m=1}^{M_i}L_t^2 \left\lvert \Tnm -t \right\rvert^{2H_t}  \left\lvert W_{m}^{(i)}(t)\right\rvert\right\} .
\end{equation}
For the NW estimator with bandwidth $h$, this has the rate
$ O_{\mathbb P} ( h^{2H_t} )$.

Gathering facts, we  deduce that, in the case of non-differentiable trajectories, with the NW estimator and 
\begin{equation}\label{eq:hinsightNW}
  h \sim (N\Mmu)^{-1/(1+2 H_t)},
\end{equation}
one can expect
\begin{equation}\label{r_diff_o}
\EEMT [ \{\widehat \mu _N(t) - \widetilde \mu _N(t) \}^2 ] = O_{\mathbb P} \left(  (N\Mmu)^{- \; \frac{2H_t}{1+2H_t}} \right) .
\end{equation}
Given the local regularity $ H_t$, the estimator $\widehat \mu_N (t)$ can thus achieve the minimax optimal rate for the estimation of the mean function $\mu(t)$. See \cite{cai2011}. 

In some cases, in particular with kernel smoothing, the estimator defined in \eqref{LP_est_v1} could be degenerate, i.e. the weights $W_{m}^{(i)}(t)$ are not well defined because $h$ is too small. The trajectories for which this happens could change with $t$. The estimator $\widehat \mu_N (t) $ is then defined as an average over the trajectories for which  the estimator \eqref{LP_est_v1} is not degenerate. This is more likely to happen in the so-called \emph{sparse} regime, where $\mathfrak m^{2H_t}  \ll  N$. A similar phenomenon occurs with estimators determined by suitable weighting schemes, see for instance, \citet[Eq.~(2.1)]{li2010}, or  \citet[Eq.~(2.3)]{zhang_sparse_2016}. However, in the independent case, one could benefit from the replication feature of functional data, because only a fraction of trajectories will yield non-degenerate estimators $\hXni_t $. The size of this fraction plays a central role in the sparse regime. This  aspect is taken into account in  Sections \ref{sec:mean_mn} and \ref{sec:LP_est_b}, where we choose the bandwidths while penalizing the number of trajectories which yield degenerate estimators.  

The case of common design requires some special attention. For simplicity, let us assume the common design points are equidistant and consider that kernel smoothing uses a  kernel  supported on  $[-1,1]$.  In this case, the bandwidth cannot have a rate  smaller than $\mathfrak m^{-1}$, otherwise the weights $W_{m}^{(i)}(t)$ could all  be equal to zero. This means that with a common design, the optimal bandwidth is given by the minimization of $ h^{2  H_t} + (N\Mmu h)^{-1}   $ under the constraint that $\Mmu h $ stays away from zero. Without loss of generality, we could set $h=k/\Mmu$ with  $k$ a positive integer and search $k$ which minimizes $ h^{2H_t} + (N\Mmu h )^{-1} $. Balancing the two terms, one expects the optimal $k/\Mmu$ to have the rate  $(N\Mmu)^{-1/\{1+2  H_t\}}$. If $\mathfrak m^{2 H_t} $ is larger than $N$, i.e., in the so-called \emph{dense} regime,  the optimal $k$ is well defined and $k\sim (\mathfrak m^{2 H_t}  / N)^{1/  \{1+2H_t\}}$ and, with this optimal choice, $\EEMT [ \{\widehat \mu _N(t) - \widetilde \mu _N(t) \}^2 ] = o_{\mathbb P} ( N^{-1}).$ 
If $\mathfrak m^{2  H_t}  \ll  N$, then the constraint that $k\geq 1$ becomes binding, and it is no longer possible to balance the squared bias term and the variance term. The rate of $ h^{2  H_t}$ dominates the rate $(N\Mmu h)^{-1} $. 
The minimal rate for $\EEMT [ \{\widehat \mu _N(t) - \widetilde \mu _N(t) \}^2] $ then  corresponds to $k=1$, and is $O_{\mathbb P} (  \Mmu^{- 2  H_t} ) .$ 
Gathering facts, we recover the optimal rate for mean estimation with common design, that  is $O_{\mathbb P} (  \Mmu^{-   H_t}+N^{-1/2} ) ,$  see \cite{cai2011}.  Finally, let us recall the somehow surprising message from \citet[p.~2332]{cai2011}: the interpolation is rate optimal when $\mathfrak m^{2   H_t} \gg N$ in the case of common design; smoothing does not improve convergence rates. Our contribution to this aspect is a data-driven rule for the practitioner which supplements this theoretical fact. The adaptive bandwidth rule proposed in Section \ref{sec:LP_est} automatically chooses between smoothing and interpolation.
% with finite sample sizes.  

We learn from the above that the ``smoothing first, then estimate'' approach can lead  to optimal rates of convergence for estimating the mean function with independent and common design, as derived by \cite{cai2011}, provided  the  local regularity   parameter $H_t$ in \eqref{eq:reg0} is known. In the next section, 
we introduce a simple estimator of this parameter. Under mild conditions, the estimator concentrates around $H_t$ faster than a suitable negative power of $ \log (\Mmu)$. 
This suffices to guarantee that our mean and covariance functions estimators achieve the same rates as when the local regularity is known. 

Let us end this section with a discussion of the differences with the weighting schemes approach, as for instance considered by \cite{li2010} and \cite{zhang_sparse_2016}. If the regularity of $\mu(\cdot)$ is known, one  could define $B_t^{(i)}$ and 
$V_t^{(i)}$ in \eqref{dec_s} centering by the mean function instead of the trajectory $\Xp{i}_t$,  derive the bound of 
$\EE [ \{\widehat \mu _N(t) -  \mu (t) \}^2 ] $, and find the bandwidth which minimizes this bound.
These steps can be found in  \cite{li2010} and \cite{zhang_sparse_2016}, where $\mu(\cdot)$ is assumed to be twice differentiable. 
However, the estimation of the regularity of  $\mu(\cdot)$ remains an open problem.

% !TeX root=bj-sample.tex

\section{Local regularity estimator}\label{subsec:golo_new}

Our approach is based on the general regularity condition \eqref{eq:reg0}, which is a local property that we formally define below. Our definition is closely related to the notion of local intrinsic stationarity introduced by \citet[p.~2060]{hsing2016}. Similar concepts of local regularity are common in continuous-time  processes, see for example, \citet[Sec.~5]{Jirak17} and the references therein, but are still little investigated in the context of functional data.   In Section~\ref{ssec:est}, we propose an estimator of $H_t$ and in Section~\ref{ssec:conc_est}, we provide theoretical guarantees. Our estimator is related to the estimation of the Hurst function of a multifractional Brownian motion, and to other regularity parameters studied in the stochastic process theory. See for example, \cite{Hsing20}, \cite{Corcuera13}, and the references therein. The existing estimators are however usually built with noiseless observations from a single sample path. Finally, given the local regularity, the Kolmogorov Continuity  Theorem allows to determine the analytic regularity of the trajectories of $X$. Details are provided in Section~\ref{sec:reg_traj}.  Hereafter,  $t\in \mathcal T$ is an arbitrary fixed point.

\subsection{Local regularity in quadratic mean}\label{lemma_NicoH}

Let $H : u\mapsto H_u \in (0,1)$ and $L : u\mapsto L_u >0$ be Lipschitz functions defined on $\mathcal T$. Let $ \Delta_* >0$ and $\Ostar(t) =[t-\Delta_*/2, t + \Delta_*/2] \cap \mathcal T$.

\begin{definition}\label{def_loc_reg}
The class $\mathcal{X}(H, L ; t, \Delta_*  )$ is the set of stochastic processes $X$
satisfying the following   conditions.
 
\begin{enumerate}[label=(H\arabic*), ref=(H\arabic*)]
\item\label{H:moments}
Constants $\mathfrak{a}>0$ and $\mathfrak{A}>0$ (not depending on $t$) exist such that, for any  $p\geq 1$,
   \begin{equation*}
                 \sup_{u\in\Ostar(t) }  \EE[ \lvert X_{u}-X_{t}\rvert^{2p}] \leq
        \frac{p!}{2} \mathfrak{a} \mathfrak{A}^{p-2}. 
    \end{equation*}    
\item\label{H:equivalent} Positive constants $S$ and  $\beta$ (not depending on $t$)  exist such that
    \begin{equation*}
        \left\lvert
        \EE[ ( X_{u} -  X_{v})^{2}]
        -  L_t^2 \lvert u-v \rvert^{2H_t}
        \right\rvert
        \leq
        S^2 \lvert u-v\rvert^{2H_t}\Delta_{*}^{2\beta},
        \quad
         u,v\in \Ostar(t),\quad u\leq t \leq v.
    \end{equation*}
\end{enumerate}
The quantity $ H_t$ is the local regularity of the process over $\Ostar(t) $,   while $L_t$ is the Hölder constant.
% of the trajectories.  
\end{definition}

The condition \ref{H:moments} serves to derive the exponential bound for the concentration of the local regularity estimator.  \cite{blanke2014} and \cite{golo2020} provide several examples and references on processes satisfying the mild condition in \ref{H:equivalent}. Examples include, but are not limited to stationary or stationary increment processes. It is worth noting that for  some common processes with the  ordered eigenvalues  of the covariance operator such that, for some $1 < \nu < 3 $, $\lambda_j \sim j^{-\nu}$, $j \geq 1$, one has $H \equiv (\nu-1)/2$. \cite{golo2020} also considers the case of differentiable trajectories, in which case the local regularity $H_t$ refers to the highest order derivative of the sample path in the neighborhood of $t$. 

Let us now consider a general class of processes satisfying  \ref{H:equivalent}, which is also used for our simulation study. We start from the  multifractional Brownian motion (MfBm) process. See, for instance, \cite{balanca2015} and the references therein for the formal definitions and the properties of this large class of Gaussian processes. 
An MfBm, say $(W_t)_{t\geq 0}$, with Hurst index function, say $t\mapsto H_t \in(0,1)$, is a centered Gaussian process with covariance function 
\begin{equation}\label{eq:cov_mfbm_b}
	C(s,t) = \EE[W(s)W(t)] =  D(H_s,H_t )\left[ s^{H_s+H_t} +  t^{H_s+H_t} - \lvert t-s \rvert^{H_s+H_t}\right] , \quad s, t\geq 0,
\end{equation}
where
\begin{equation}\label{eq_D}
	D(x,y )=\frac{\sqrt{\Gamma (2x+1)\Gamma (2y+1)\sin(\pi x)\sin(\pi y)}} {2\Gamma (x+y+1)\sin(\pi(x+y)/2)} , \quad D(x,x) = 1/2, \quad x,y >0.
\end{equation}
To make the MfBm class even more general, we also consider a deterministic time deformation. The time deformation is defined here by $t\mapsto A(t)\geq 0$, a strictly increasing, continuously differentiable function defined on $[0,\infty)$. Moreover, the derivative $A^\prime (t)$ is strictly positive on any compact interval. Let $A^{-1}(\cdot)$ denote the inverse of $A(\cdot)$, and let 
$H_{A,t} = H_{A^{-1}(t)}$. 
We consider the  MfBm $(W_{A,t})_{t\geq 0}$ with Hurst index function $H_{A,t}$. Given the Hurst index function $H$ and the time deformation function $A$, the process we consider is  
\begin{equation}\label{eq:X_md}
	X(t) = W_A (A(t)), \quad t\geq 0,
\end{equation}
with covariance function
\begin{equation}\label{eq:cov_mfbm_c}
	C_A(s,t) = \EE[X(s)X(t)] =  D(H_s,H_t )\left[ A(s)^{H_s+H_t} +  A(t)^{H_s+H_t} - \lvert A(t)-A(s)\rvert^{H_s+H_t}\right] .
\end{equation}

\begin{lemma}\label{lemma_mfbm}
	Assume $ t\mapsto H_t \in (0,1)$ is twice continuously differentiable, and $ t\mapsto L_t >0$ is continuous, $t\geq 0$. Then, $X$ defined in \eqref{eq:X_md} satisfies condition \ref{H:equivalent} with local regularity $H_t$ and Hölder constant $L_t$, provided that, for some $A(0)\geq 0$, the time deformation is 
	$$
	A(t) = A(0)+ \int_{0}^t L_s^{1/H_{s}} ds, \quad t\geq 0.
	$$
\end{lemma}

\subsection{The local regularity estimation method}\label{ssec:est}
 
Assume that $X$ belongs to $\mathcal{X}(H,L;t, \Delta_*  )$. Our first goal is to construct an estimator of $H_t$. 
For simplicity, for $u,v\in\Ostar(t)$, $u\leq t \leq v$, let us denote
\begin{equation*}
    \theta(u,v) = \EE[ (X_u -  X_v)^{2} ]
    \approx L_t^2 \lvert u-v\rvert^{2H_t} \quad \text{if }  \Delta_* \text{ is small}.
\end{equation*}
Now, let $t_1$ and $t_3$ be such that $[t_1,t_3]\subset\Ostar(t)$ and $t_3-t_1=\Delta_*/2$. Let $t_2=(t_1+t_3)/2$ and define
\begin{equation}\label{eq:tilde-Hd}
    \widetilde H_{t} =  \frac{\log(\theta(t_1,t_3)) - \log(\theta(t_1,t_2))}{2\log(2)} \quad \text{if }  \Delta_* \text{ is small}.
\end{equation}
When $t$ is offset from the left and right endpoints of $\mathcal{T}$ by more than $\Delta_*/2$, we set $t_2 = t$. Otherwise, we set $t_1 = \min\mathcal T $ or $t_3 = \max\mathcal T$, respectively. Since $H$ is Lipschitz continuous and, by construction, $\lvert t_2-t \rvert\leq\Delta_*/2$, the quantity $\widetilde H_{t}$ is a proxy of $H_t$. A complementary discussion of the choice of $t_1,t_2$ and $t_3$ in \eqref{eq:tilde-Hd} is provided in \cite{supplement}. 

Given nonparametric estimators $\widetilde{X}^{(i)}_u$ of $X^{(i)}_u$, we define a natural estimator of $\widetilde H_t$, and thus of $H_t$, as
\begin{equation}\label{eq:hat-h}
    \widehat H_t = \frac{\log\big( \widehat\theta(t_1,t_3)\big) - \log \big(\widehat\theta(t_1,t_2)\big)}{2\log (2)}, \quad \text{where} \quad
    \widehat\theta(u,v) = \frac1{N} \sum_{i=1}^{N}
    \left(
    \widetilde{X}^{(i)}_u - \widetilde{X}^{(i)}_v
    \right)^2,
    \quad u,v\in\Ostar(t).
\end{equation}
The estimate of $L_t$ is readily obtained given  $\widehat H_t$, the details being provided in Section \ref{sec:imp_asp}.

\subsection{Concentration properties of the local regularity estimator}\label{ssec:conc_est}

The local regularity  estimator $\widehat H_t$ in~\eqref{eq:hat-h} is studied by \cite{golo2020} in the case of constant functions $H$ and $L$ in a neighborhood of $t$. 
The quality of $\widehat H_t$ depends on the generic nonparametric estimator $\widetilde{X}_u$ of $X_u$. To quantify their behavior, we consider the local $\mathbb{L}^p$-risk 
\begin{equation*}
    R_{\mathfrak m}(t;p) = \sup_{u\in\Ostar(t)} \EE[\lvert\widetilde{X}_u - {X}_u\rvert^p ],
    \quad p\geq 1.
\end{equation*}
The risks $R_{\mathfrak m}(t,p)$ depend on $\mathfrak m$, the average number of points  on each curve. Our methodology applies to  any type of nonparametric estimator $\widetilde{X}$ (local polynomials, splines, etc.) as soon as,  for any $p\in\NN$, its $\mathbb{L}^p$-risk  is  suitably bounded. The following mild condition is satisfied by common estimators, see for instance, \citet[Th.~1]{GAIFFAS2007} for the case of local polynomials. 

\begin{assumption}\label{LP:1}
There exist two positive constants $\mathfrak{c}$ and $\mathfrak{C}$ such that, for any $p\geq 1$,
    \begin{equation*}
        R_{\mathfrak m}(t;2p)
        \leq 
        \frac{p!}{2} \mathfrak{c} \mathfrak{C}^{p-2}, \quad \mathfrak m \geq 1.
    \end{equation*}
\end{assumption}
We can now state a non-asymptotic concentration result for the estimator $\widehat H_t$.

\begin{theorem}\label{thm:estimation-alpha}  
  Assume that  $X$ belongs to $\mathcal{X}(H,L;t, \Delta_*  )$, and that Assumption~\ref{LP:1} holds true.   
   Assume also that there exist $\tau>0$ and $B>0$ such that $R_{\mathfrak m}(t;2)  \leq B\Mmu^{-\tau}$. Let 
     $1< \varrho$ and $0< \gamma$,    and consider 
    \begin{equation*}
       \varphi(\Mmu) = \log^{-\varrho}(\Mmu) \quad  \text{and}  \quad \Delta_* / 2  = \exp(-\log^{\gamma}(\Mmu)) .  
    \end{equation*}
    Then, for any $\Mmu$ larger than some constant $\Mmu_0$ depending on $B$, $\tau$,  $\gamma$, $\rho$,   $H$, $\beta$ and for  some constant $\mathfrak{f} $, %we have
    \begin{equation*}
        \PP\left[ \left\lvert\widehat H_t - H_t\right\rvert > \varphi(\Mmu), \, \widehat H_t>0 \right]
        \leq
        \exp\left( -\mathfrak{f} N\varphi^2(\Mmu) \Delta_*^{4H_t} \right).
    \end{equation*}
\end{theorem}

The proof of Theorem \ref{thm:estimation-alpha} follows the lines of \citet[Th.~5]{golo2020} and is thus omitted. Let us however point out that the three quantities $R_{\mathfrak m}(t;2)$, $\Delta_* $ and $ \varphi(\Mmu)$ are required to decrease to zero, as $\Mmu$ tends to infinity with $N$, in such a way that 
\begin{equation}\label{hertm}
R_{\mathfrak m}(t;2)/\Delta_*^a + \Delta^{1/a}_*/\varphi(\Mmu)\rightarrow 0, \;\, \text{ for some } \;\; a>0.
\end{equation}
First, the choice of  $\varphi(\Mmu)$ is such that the effect of estimating $H_t$ does not deteriorate the pointwise rates for mean and covariance function estimation.
  Imposing the mild condition that $\log(N)/\log(\mathfrak m)$ is bounded, since $\Mmu^{1/\log(\Mmu)} = e$, the effect of using the bandwidth in~\eqref{eq:hinsightNW} with   $H_t$  replaced by $\widehat H_t$ is negligible as soon as $\varphi(\Mmu) \ll \log^{-1}(\Mmu)$. 
Second, since $\tau>0$ could be arbitrarily small, the rate imposed on the pilot estimator $\widetilde{X}$ of $X$, is a very mild consistency requirement. It is achieved by the common nonparametric estimators under general conditions on the smoothing parameters, with random or fixed design, and mild conditions on the distribution of the $M_i$. See, for instance, \cite{tsybakov2009} and \cite{BELLONI2015}. In particular, the required rate for the $\widetilde{X}$ can be obtained under general forms of heteroscedasticity. 
These facts explain the  choice of $\Delta_*$ which makes \eqref{hertm} to hold true. 
In conclusion, the only practical choice is that of $\gamma$, and we set   $\gamma = 1/3$ in our empirical study.

\subsection{From the regularity of the process to the regularity of the trajectories}\label{sec:reg_traj}

Let us now connect the probabilistic concept of local regularity with the regularity of the sample paths considered as functions.  When  $H$ is constant in a neighborhood of $t$,  by Assumptions~\ref{H:moments}, \ref{H:equivalent}, using the version of Kolmogorov's criterion stated in \citet[Th.~2.1]{MR1725357}, it can be proven that almost all sample paths of $X$  belong to any Hölder space of functions defined over the neighborhood of $t$, with the Hölder exponent less than $H$. As an example, the Brownian motion has a constant local regularity equal to 1/2. Moreover, almost surely, the sample paths of the Brownian motion belong to any  Hölder space of local regularity less than $1/2$, but cannot be Hölder continuous with exponent greater than or equal to $1/2$. The general form of the Hölder regularity of the MfBm sample paths, which depends on $H_t$, can be found for example in \cite{balanca2015}.

Hence, the probability theory indicates that imposing assumptions on the regularity of the sample paths could be a delicate issue. Indeed, even for some widely used examples, this regularity is not well defined in the sense required by the nonparametric statistics  theory. Since the sample paths have a regularity which can be arbitrarily close to the local regularity of the process $X$ as defined above,  the probabilistic concept  of local regularity  seems more appropriate for establishing the rates of convergence for the mean and covariance estimators.

% !TeX root=bj-sample.tex

\section{Adaptive mean and  covariance function estimators}\label{sec:LP_est}

We now explain how to select data-driven bandwidths for kernel smoothing of the trajectories and build adaptive mean and covariance function estimates.
Hereafter,  $\widehat  H_t$ will be the estimator  of $H_t$ defined in \eqref{eq:hat-h}, considered on the event $\{\widehat H_t>0\}$. Let $\widehat\Mmu = N^{-1} \sum_{i=1}^N M_i$. Let us consider a class of linear smoothers of the sample paths. For each $1\leq i \leq N$, using the measurements $(\Ynm , \Tnm )$, $1\leq m \leq M_i$, of the trajectory $\Xp{i}$, we define $\hXni_t$ as in \eqref{LP_est_v1}, 
where $W_{m}^{(i)}(t)$ are  weights depending on the $\Tnm$'s only, and on some smoothing parameter. In the following, we focus on the case of Nadaraya-Watson (NW) estimators,   but also indicate how to adapt the construction for local linear smoothing.  Given the bandwidth $h$, with the convention $0/0=0$, the weights of the NW estimator of  $X^{(i)}$ are
\begin{equation}\label{LP_weights}
W_{m}^{(i)}(t)  =W_{m}^{(i)}(t;h)  =  K\!\left( (\Tnm-\T)/h \right)\left[\sum_{m'=1}^{M_i} K\left( (T^{(i)}_{m'}-\T)/h \right)\right]^{-1},
\quad 1\leq m \leq M_i.
\end{equation}
Herein, $K$ is a bounded density with the support $[-1,1]$, and $\inf K>0$ on a sub-interval of the support.

\subsection{Adaptive optimal mean  estimation}\label{sec:mean_mn}

With finite samples it may happen that $\hXni_t$ is degenerate, i.e. $W_{m}^{(i)}(t) =0$ for all $1\leq m \leq M_i$. In such a case, the $i$th curve will be dropped for the mean and covariance estimations. With kernel smoothing, in the common design case, the number of degenerate estimates $\hXni_t$ is either equal to $N$ or to zero. In the independent design case, this number could be any integer between 0 and $N$. A suitable bandwidth rule should be based on a risk which includes a penalty for the number of curves which are not considered for the estimation. In the following, we propose a  natural way to  penalize which adapts to the sparse and dense regimes. Moreover, the two types of design are  handled automatically. For this purpose, let $\mathbf{1}\{\cdot\}$ denote the indicator function, define
\begin{equation}\label{def_wi}
w_i (t;h) = 1 \quad \text{if } \; \sum_{m=1}^{M_i} \mathbf{1}\left\{\lvert\Tnm-t\rvert\leq h\right\} \geq 1 \quad \text{and} \quad w_i (t;h) = 0  \quad\text{otherwise},
\end{equation}
 and let
$
\WN(t;h) =\sum_{i=1}^N w_i (t;h).
$
By construction, $w_i (t;h) = 0$ if and only if $W_{m}^{(i)}(t;h) =0$ for all $1\leq m \leq M_i$.

Our adaptive mean function estimator is 
\begin{equation}\label{mu1_d}
\widehat \mu_N^* (t) = \widehat \mu_N (t;h_{\mu} ^*) \quad \text{with} \quad \widehat \mu_N (t;h ) = \frac{1}{\WN(t;h ) }\sum_{i=1}^N  w_i (t;h )  \hXni_t,
\end{equation}
where $h_{\mu} ^* $  is  a suitable bandwidth defined below. The mean estimator $\widehat \mu_N (t;h) $ is a version of that defined in \eqref{est_mean0} which takes into account that some trajectories have no observation times between  $t- h$  and  $t+ h$.  The normalization of the mean estimator by $\WN(t;h)$ is also implicitly used in the definition of the estimators proposed by \cite{li2010} and  \cite{zhang_sparse_2016}. 

To introduce our bandwidth rule, for any  $h>0$, $\alpha >0$, let  
\begin{equation}\label{eq:ci}
c_i (t;h)= \sum_{m=1}^{M_i} \left\lvert W_{m}^{(i)}(t;h)\right\rvert,  
\quad c_i (t;h,\alpha )  = \sum_{m=1}^{M_i} \left\lvert (\Tnm - t)/h\right\rvert^{\alpha} \left\lvert W_{m}^{(i)}(t;h)\right\rvert,
\end{equation}
and 
\begin{equation}\label{eq:ci2}
 \overline{C} (t;h,\alpha)= 
\frac{1}{\WN(t;h)  } \sum_{i=1}^{N}w_i (t;h)   c_i (t;h)c_i (t;h,\alpha ). 
\end{equation}
In \eqref{eq:ci},  $W_{m}^{(i)}(t;h)$ can be the weights corresponding to local polynomial smoothing. 
With the NW estimator,  $w_i (t;h)   c_i (t;h) =w_i (t;h) $.  

Moreover, 
\begin{equation}\label{eq:ci2b}
\overline{C}(t;h,\alpha)\approx \int|u|^{ \alpha}K(u)du.
\end{equation}
The details for \eqref{eq:ci2b} are provided in \cite{supplement}. Using the equivalent kernels idea, see \citet[Sec. 3.2.2]{fan_local_1996}, the same approximation could be used in the case of local linear estimators.  
The accuracy of the approximation \eqref{eq:ci2b} could be high since it involves the $\Tnm$ to be close to $t$ for all the curves with $w_i (t;h)=1$. 
Next, using the rule $0/0=0$, let 
\begin{equation}\label{eq:Ncal1}
\mathcal N_{i} (t;h) =  \frac{w_i (t;h)}{\max_{1\leq m\leq M_i} \lvert W_{m}^{(i)}(t;h) \rvert } \quad \text{and} \quad
\mathcal N_{\mu} (t;h) = \left[\frac{1}{\WN^2(t;h)}\sum_{i=1}^N  w_i (t;h)\frac{c_i (t;h) }{\mathcal N_{i} (t;h)}\right]^{-1} .
\end{equation}
With the NW estimator, $\mathcal N_{\mu} (t;h) $ is equal to  $\WN(t;h) $ times the harmonic mean of $\mathcal N_{i} (t;h) $, over the curves with $w_i (t;h) = 1 $. 

Let $\mathcal H_N$ be a bandwidths range.  We define the  bandwidth for computing $\widehat \mu_N^* (t) $ such that it minimizes the mean squared difference between $\widehat \mu_N (t;h) $ and $\widetilde \mu_N (t)$.  This leads us to define the optimal  bandwidth 
\begin{equation}\label{eq:opt_h_mean0b}
h_{\mu} ^* =h_{\mu} ^* (t) =   \arg\min_{h\in\mathcal H_N} \mathcal{R}_\mu (t;h) ,
\end{equation}
with 
\begin{equation}\label{eq:opt_h_mean0c}
\mathcal{R}_\mu(t;h)= 
\!q_1^2h^{2\widehat H_t} + \frac{ 
q_2^2}{\mathcal N_\mu (t;h)}+
 q_3^2\left[ \frac{1}{\WN(t;h)} - \frac{1}{N}\right],
\end{equation}
and 
\begin{equation}\label{eq:les_q}
  q_1^2 =  \overline{C} (t;h,2\widehat H_t)    
  \widehat L_{t}^2\;,
  \quad  
q_2^2 = 
\sigma^2_{\max}\;
, 
\quad
q_3 ^2= 
\operatorname{Var}[X_t],
\end{equation}
where $\sigma_{\max}$ is a bound for the function $\sigma(t,x)$  in \eqref{def_err} and  $ \widehat L_{t}$ is an estimate of the  Hölder constant $L_t$ from \ref{H:equivalent}.  
In Section \ref{sec:emp}, we propose  a simple procedure to build $\widehat L_{t}$ based on $\widehat H_t$. 
We show in the Appendix that
$\mathcal{R}_\mu (t;h)/2$ is a sharp bound for $\EEMT[ \{   \widehat \mu_N (t;h) - \widetilde \mu_N (t)  \}^2 ]$.
The minimization of $\mathcal{R}_\mu (t;h)$ can be easily  performed on a grid of  $h$ values  in the range $\mathcal H_N$. 
 
The bandwidth rule \eqref{eq:opt_h_mean0b} could be used with both independent and common design. With common design, the $\Tnm \equiv T_m$ and $W_{m}^{(i)}(t;h)\equiv W_{m}(t;h)$ no longer depend on $i$ and the solution $h_{\mu} ^*$ will always be a value in the set of $h$ such that $\WN(t;h) =N$. Moreover, for the NW estimator, whenever $\WN(t;h) =N$, we have
\begin{equation}\label{eq:qdf1}
\overline{C} (t;h,2\widehat H_t) = \sum_{m=1}^{\mathfrak m} \left\lvert(T_m - t)/h\right\rvert^{2\widehat H_t} W_{m}(t;h)  
\quad \text{and} \quad
 \mathcal N_{\mu}^{-1} (t;h) = N^{-1} \max\limits_{1\leq m\leq \mathfrak m} W_{m}(t;h).
\end{equation}
In a  data-driven way, the bandwidth $h_{\mu} ^*$ automatically chooses between interpolation and smoothing. 

The following result states that our estimator $\widehat \mu_N^* (t)$ 
achieves the best rates one can expect.   We assume
\begin{equation}\label{eq:cdg_H}
N\mathfrak m \min \mathcal H_N /\log(N\mathfrak m) \rightarrow \infty\quad \text{and}\quad \max \mathcal H_N \rightarrow 0,
\end{equation}
a minimal condition for the bandwidth range.
For simplicity,  we also assume that
\begin{equation}\label{mfka}
\limsup_{N,\mathfrak m\rightarrow \infty}  \{\log(N)/\log(\mathfrak m)\} <\infty,
\end{equation}
a technical condition which is realistic in applications. Moreover, we impose the following mild technical condition in the independent design case:
\begin{equation}\label{mfkb}
\exists c_L, C_U>0 \quad \text{such that} \quad c_L \leq M_i\mathfrak m^{-1} \leq C_U, \quad \text{for } 1\leq i \leq N \quad \text{and all } N.
\end{equation}
With a common design where $M_i\equiv \mathfrak m$ and  the $T_1^{(i)},\ldots,T_{\mathfrak m} ^{(i)}$ are not changing with $i$, we suppose that:
\begin{equation}\label{mfkb_2}
\exists C_U \geq 1 \quad \text{such that} \quad  \max_{1\leq m\leq \mathfrak m -1}\{T_{m+1}^{(i)} -T_{m}^{(i)}  \}\leq  C_U\min_{1\leq m\leq \mathfrak m - 1}\{T_{m+1}^{(i)} -T_{m}^{(i)}  \}   .
\end{equation}
Below,  $\asymp_\PP$ means left side is bounded above and below by positive constants times the right side, with probability tending to 1, provided $N$ (and thus $\mathfrak m$) tends to infinity.

\begin{theorem}\label{thm:optim_mu}
 Assume the conditions of Theorem \ref{thm:estimation-alpha}, and assume \eqref{eq:cdg_H}, \eqref{mfka}  hold true. Assume that $\Tnm$ are either independently drawn, with a Hölder  continuous density which is bounded away from zero and \eqref{mfkb} holds true, or $\Tnm$ are the points of a common design satisfying \eqref{mfkb_2}. Then, in the independent design case,
\begin{equation}\label{eq:hmu_opt}
h_{\mu} ^*  \asymp_\PP  (N\Mmu)^{- \; \frac{1}{1+2H_t }} ,
\end{equation}
and the estimator 
$\widehat \mu_N^* (t) = \widehat \mu_N (t;h_{\mu} ^*)$ defined by \eqref{mu1_d} and \eqref{eq:opt_h_mean0b} satisfies 
$$
\widehat \mu_N^* (t) - \widetilde \mu_N (t) = O_{\mathbb P} \left(  (N\Mmu)^{- \; \frac{H_t }{1+2H_t }} \right) \quad  \text{and} \quad  \widehat \mu_N^* (t) - \mu (t) = O_{\mathbb P} \left((N\Mmu)^{- \; \frac{H_t }{1+2H_t }}+  N^{- 1/2} \right).
$$
Meanwhile,  with the common design,
$$
\widehat \mu_N^* (t) -   \mu (t) = O_{\mathbb P} \left( \max\left\{  (N\Mmu)^{- \; \frac{H_t }{1+2H_t }} , \Mmu^{- H_t }\right\}+ N^{-1/2}\right) =   O_{\mathbb P} \left(  \Mmu^{- H_t }+  N^{- 1/2}\right)  . 
$$

\end{theorem}

The rates of $\widehat \mu_N^* (t) $ are the best one could expect in view of the results of \cite{cai2011}. The difference between the common and independent designs comes from the fact that, in order to avoid a degenerate mean estimator, the bandwidth cannot decrease faster than~$\Mmu^{-1}$.

\subsection{Adaptive covariance function estimates}\label{sec:LP_est_b}

For any $s, \T \in \mathcal T $, $s\neq t$, define
\begin{equation}\label{def_wi_b}
w_i (s,t;h) = w_i (s;h) w_i (t;h) \quad \text{and} \quad \WN(s,t;h) =\sum_{i=1}^N w_i (s,t;h),
\end{equation}
with $w_i (s;h)$ and $ w_i (t;h)$  as in \eqref{def_wi}. 
Our adaptive covariance function estimator is 
\begin{equation}\label{cov1_d}
\widehat \Gamma ^*_N(s,t)=\widehat \Gamma _N(s,t;h_{\Gamma} ^*)\quad \text{with}\quad 
\widehat \Gamma _N(s,t;h)=  \widehat \gamma _N(s,t;h) - 
\widehat \mu_N (s;h) \widehat \mu_N (t;h), 
\end{equation}
where  $\widehat \mu_N (s;h) $, $\widehat \mu_N (t;h)$ 
are defined according to \eqref{mu1_d},  and 
\begin{equation}\label{mu1_d2}
 \widehat \gamma_N(s,t;h)  = \frac{1}{\WN(s,t;h) }\sum_{i=1}^N w_i (s,t;h)\hXni_s \hXni_t .
\end{equation}
Here,  $\hXni_s$ and $\hXni_t$ are the NW estimators built with some suitable bandwidth  $h_{\Gamma} ^* $ which is  defined below. This covariance function estimator is a practical version of the one defined in \eqref{est_cov1}. The normalization of the covariance estimator by $\WN(s,t;h)$ is also implicitly used in the definition of the estimators proposed by \cite{li2010} and \cite{zhang_sparse_2016}. 

We  define the  bandwidth for computing $ \widehat \gamma_N(s,t;h) $, and eventually $\widehat \Gamma_N^* (s,t) $, such that it minimizes the mean squared difference between $\widehat \gamma_N (s,t;h) $ and the unfeasible estimator 
$$
\widetilde \gamma_N (s,t) = N^{-1} \sum_{i=1}^N \Xp{i}_s \Xp{i}_t ,
$$ 
of  $\EE[X_s X_t]$.  To this aim, we define modified versions of $\mathcal N_{i} (t;h)$ and $\mathcal N_{\mu} (t;h)$, see \eqref{eq:Ncal1}, only taking into account  the curves with $w_i(s,t;h)=1$:
\begin{equation}%\label{eq:Ncal2def1}\label{eq:Ncal2def2}
    \mathcal N_{i} (t \mid s;h) =  \frac{w_i (s,t;h)}{\max_{1\leq m\leq M_i} \lvert W_{m}^{(i)}(t,h)\rvert},
\quad
\text{and}
\quad
    \mathcal N_{\Gamma} (t \mid s;h) =  \left[\frac{1}{\WN^2(s,t;h)}\sum_{i=1}^N \frac{ w_i (s,t;h) }{\mathcal N_{i} (t|s;h)}\right]^{-1}.
\end{equation}
This idea leads us to define the optimal bandwidth,   in some range $\mathcal H_N$,   as 
\begin{equation}\label{eq:opt_h_cov0b}
h_{\Gamma} ^* = h_{\Gamma} ^* (s,t) = h_{\Gamma} ^* (t,s) =   \arg\min_{  h\in \mathcal H_N  }\{ \mathcal{R}_\Gamma (s \mid t;h) +\mathcal{R}_\Gamma (t \mid s;h) \},
\end{equation}
with
\begin{equation}\label{eq:opt_h_cov0c}
 \mathcal{R}_\Gamma (t \mid s;h)=  \mathfrak{q}_1^2 (t \mid s) h^{2\widehat H_t}
+  \frac{\mathfrak q_2^2(t \mid s)} {\mathcal N_\Gamma(t \mid s;h)} 
+ \mathfrak{q}_3^2 
\left[ \frac{1}{\WN(s,t;h)} - \frac{1}{N}\right]. 
\end{equation}
The $\mathfrak{q_\ell}$, $1\leq l \leq 3$, are defined by:
\begin{equation}\label{eq:les_q_frak}
\mathfrak{q}_1^2 (t \mid s) =   2  \EE[X_s^2]  \overline{\mathfrak{C}} (t \mid s;h,2\widehat H_t)    \widehat L_{t}^2 ,
\quad
\mathfrak{q}_2^2  (t|s)= \sigma^2_{\max}\EE[X_s^2],\quad
\mathfrak{q}_3^2 = \frac{\operatorname{Var}[X_sX_t]}{2}, 
\end{equation}
where 
\begin{equation}\label{eq:cbar_g}
   \overline{\mathfrak{C}} (t \mid s;h,\alpha)=    \frac{ \sum_{i=1}^{N}w_i (s,t;h) 
c_i(t;h,\alpha)  }{\WN(s,t;h)  } \approx\int|u|^{ \alpha}K(u)du,
\end{equation}
and the approximation \eqref{eq:cbar_g} is valid with NW or local linear estimators. The details for \eqref{eq:opt_h_cov0c} are provided in \cite{supplement}.

We show in \cite{supplement} that  the function of $h$ minimized in \eqref{eq:opt_h_cov0b}  is a sharp bound for $\EEMT[ \{   \widehat \gamma_N (s,t;h) - \widetilde \gamma_N (s,t)  \}^2 ]/2$,  which is the leading term of $\EEMT[ \{   \widehat \Gamma_N (s,t;h) - \widetilde \Gamma_N (s,t)  \}^2 ]/2$ with $\widetilde \Gamma_N (s,t) = \widetilde \gamma_N (s,t) - \widetilde \mu_N (s)\widetilde \mu_N (t)$. 
The sum of the first two terms in the expressions of  $\mathcal{R}_\Gamma (s \mid t;h)$ and $\mathcal{R}_\Gamma (t \mid s;h)$  represents the  quadratic risk of our estimator of $\EE[X_sX_t]$ compared to the unfeasible one based on the true values $\Xp{i}_s\Xp{i}_t$ from the curves yielding non-degenerate estimates $\hXni_s\hXni_t$. 
The third term in  \eqref{eq:opt_h_cov0c} penalizes for the number of curves which are dropped when calculating our estimator. The minimization in \eqref{eq:opt_h_cov0b} can be  done on a grid of values $h$.  

Like for the mean function, 
the definition \eqref{eq:opt_h_cov0b} can be used with both independent and common design. Indeed, with common design,  $h_{\Gamma} ^*$ will always be a value in $\mathcal H_N$ such that $\WN(s,t;h) =N$. In a completely data-driven way, $h_{\Gamma} ^*$ will  choose between interpolation and smoothing. 

\begin{theorem}\label{thm:optim_gamma}
Let 
$s\neq t$. Assume    $N\{\mathfrak m \min \mathcal H_N \}^2/\log^2(N\mathfrak m) \rightarrow \infty$,   $\sup_{t\in\mathcal T }\EE(X_t^4)<\infty $, and the conditions of Theorem \ref{thm:optim_mu} hold true. Let $H(s,t) = \min \{H_s,H_t \}$.  Then, in the independent design case,
\begin{equation}\label{eq:hgamma_opt}
  h_{\Gamma} ^* \asymp_\PP   \max\left\{  (N\Mmu^2)^{- \; \frac{ 1}{2 \{H(s,t)+1\}} }, \;   (N\Mmu)^{- \; \frac{ 1}{2 H(s,t)+1} } \right\},  
\end{equation}
and the estimator 
$\widehat \Gamma_N^* (s,t) = \widehat \Gamma_N^* (s,t;h_{\Gamma} ^*)$ defined by \eqref{cov1_d} and \eqref{eq:opt_h_cov0b} satisfies 
\begin{equation}\label{eq:rate_G*}
   \widehat \Gamma_N^* (s,t) -   \Gamma (s,t)  =O_{\mathbb P} \left(  (N\Mmu^2)^{- \; \frac{ H(s,t)}{2 \{H(s,t)+1\}} } +   (N\Mmu)^{- \; \frac{ H(s,t)}{2 H(s,t)+1} } +  N^{-1/2}\right).
\end{equation}
Meanwhile  with the common design,
\begin{equation}\label{eq:rate_CG*}
 \widehat \Gamma_N^* (s,t) - \Gamma (s,t)  
= O_{\mathbb P} \left(  \Mmu^{-  H (s,t)}\!+ N^{-1/2}\right)  .
\end{equation}
\end{theorem}

In view of the results of \cite{cai2010}, the rate achieved by $\widehat \Gamma_N^* (s,t) $ is the best one could expect in the case of common design. However, with independent design, 
\begin{equation}\label{zdwr}
\Mmu^{2H(s,t)} \ll N \quad \text{if and only if} \quad 
N^{-1/2}\ll (N\Mmu)^{- \; \frac{ H(s,t)}{2 H(s,t)+1} }
\ll  (N\Mmu^2)^{- \; \frac{ H(s,t)}{2 \{H(s,t)+1\}} },
\end{equation}
and thus the rate of $\widehat \Gamma_N^* (s,t) $ is slower than one may expect.  Even if this rate seems sub-optimal in the sparse case, with respect to the minimax rate for the  $\mathbb{L}^2$-risk obtained by \cite{cai2010}, we conjecture that $\widehat \Gamma_N^* (s,t) $ achieves the optimal pointwise rate.
We leave the clarification of this subtle  aspect for future theoretical work. 

\subsection{The estimator on the diagonal of the covariance function}

As mentioned in \eqref{est_cov1}, we propose to use the estimator of $\Gamma(s,t)$, defined in  \eqref{mu1_d2}, only 
outside a diagonal set $\mathcal D$. It remains to give a data-driven rule for choosing $\mathcal D$ shrinking to the diagonal segment $\{(s,s):s\in\mathcal T\}$, and to propose an estimator for the covariance function on the diagonal set.  Let us fix $t\in\mathcal T$, and consider $\mathfrak d_t \leq \Delta_*/2 $, with $\Delta_*$ like in Theorem \ref{thm:estimation-alpha}.
Under mild moment assumptions, 
we show in the Appendix that, a constant $C$ exists such that 
\begin{equation}\label{Dtilde_d}
\mathbb E\left[   \left\{  \widetilde \Gamma_N ( t- \mathfrak  u_1 , t+ \mathfrak {u}_2) -  \widetilde \Gamma_N (t,t) \right\}^2  \right]  \leq C  \mathfrak d_t^{2H_t}, \quad  0 \leq \mathfrak  u_1,  \mathfrak  u_2 \leq   \mathfrak d_t  .
\end{equation}
On the other hand, for proving Theorem \ref{thm:optim_gamma}, we need a bandwidth smaller than $|s-t|/2$ for a kernel with support in $[-1,1]$. Taking into account these aspects, our estimator of $\Gamma_N ( t- \mathfrak  u_1 , t+ \mathfrak {u}_2)$ and $\Gamma_N ( t+ \mathfrak  u_2 , t- \mathfrak {u}_1)$, for $0\leq \mathfrak u_1,\mathfrak u_2 \leq  \mathfrak d_t$, is defined as
\begin{equation}\label{eq:Gam_di}
\widehat \Gamma_N ( t- \mathfrak  u_1 , t+ \mathfrak {u}_2) = \widehat \Gamma_N ( t+\mathfrak  u_2 , t- \mathfrak {u}_1) =\widehat \Gamma_N(t -\mathfrak d _t , t +\mathfrak d_t ).
\end{equation}
The quantity $\mathfrak d_t$ can be the smallest value $d$ which is larger than the   bandwidth $h_{\Gamma} ^*(t -  d  , t + d )$ defined in \eqref{eq:opt_h_cov0b}. In practice, one can simply consider the points $(t-d, t+d)$ on a grid, for decreasing values of $d$. The value $\mathfrak d_t$ is then the smallest  $d$ for which $d\geq h_{\Gamma} ^* (t-d, t+d)$.

\subsection{Implementation aspects}\label{sec:imp_asp}

The risks $\mathcal R _\mu$ and $\mathcal R_\Gamma$ defined in  \eqref{eq:opt_h_mean0c} and \eqref{eq:opt_h_cov0c}, respectively, depend on the second order moments of $X_t$ and $X_sX_t$,  on $ L^2_t$, and on the conditional variance bound $\sigma^2_{\max}$. For the second order moments, we simply use empirical moments with $X$ replaced by $\widetilde X$. 
To obtain the presmoothed curves $\widetilde X$ introduced in Section \ref{ssec:est}, we use  the NW estimator using the bandwidth defined in \cite{BERTIN2004225} and the triangular kernel 
$K(t) = \left(1 - \lvert t\rvert \right)\mathbf{1}_{[-1, 1]}(t)$.

In view of \ref{H:equivalent}, with  $[t_1,t_3]\subset\Ostar(t)$ and $t_3-t_1=\Delta_*/2$, if $t_2$ is the midpoint of $[t_1, t_3]$, 
\begin{equation}\label{eq:tilde-Ld}
    L_t^2 \approx  \frac{\theta (t_2, t_3)}{\lvert t_3 - t_2\rvert^{2H_t}} \approx \frac{\theta (t_1, t_2)}{\lvert t_2 -t_1 \rvert^{2H_t}}.
\end{equation}
Given the estimate $\widehat H_t$ and estimates $\widehat\theta (t_2, t_3)$ and $\widehat\theta (t_1, t_2 )$ as in \eqref{eq:hat-h}, 
we then define the estimate % of $L^2_t$ as
\begin{equation}\label{eq:hat-l}
    \widehat{L}_t^2 = \frac{1}{2}\left(
\frac{\widehat \theta (t_2, t_3)}{\lvert t_3 - t_2\rvert^{2\widehat H_t}} + \frac{\widehat \theta (t_1, t_2)}{\lvert t_2 -t_1 \rvert^{2\widehat H_t}}
 \right) .
\end{equation}

To estimate the conditional variance bound, let us first consider the case where $\sigma^2(t,x)$ does not depend on $x$. In this case, one can compute  
\begin{equation}\label{eq:sigma2_estim}
\widehat{\sigma}^2(t)= \frac{1}{N}\sum_{i = 1}^{N} \frac{1}{2\lvert\mathcal{S}_i(t) \rvert}\sum_{m \in  \mathcal S_i (t)} \left\{Y_{m}^{(i)} - Y_{m-1}^{(i)}\right\}^2,
\end{equation}
where $\mathcal S_i (t)$ is a subset of indices $m$ for the $i$th trajectory, and $\lvert\mathcal S_i (t)\rvert$ denotes its cardinal. 
When the variance of the errors is considered constant, $\mathcal S_i (t)$ can be the set equal to $\{2,3,\ldots, M_i\}$ for all $t$.  
When the variance depends on $\T$, one could define $\mathcal S_i(t)$ as the set of indices corresponding to the $ K_0$ values $ T_{m}^{(i)}$ closest to $t$. The theory allows for a  choice such as 
$ K_0 = \lfloor \mathfrak{m}\exp(-\{\log\log  \mathfrak{m} \}^2)\rfloor$. Then $\sigma^2_{\max} $ could be $ \max_{t\in\mathcal T} \widehat{\sigma}^2(t)$, and this choice is used in our empirical investigation.

% !TeX root=bj-sample.tex

\section{Empirical study}\label{sec:emp}

To investigate the finite sample properties of our adaptive nonparametric estimators of the mean and covariance functions, we proceed to an extensive simulation study,  using the class of MfBm with time deformation introduced in  Section \ref{lemma_NicoH}.  We next use the functions $u\mapsto H_u$ and $u\mapsto L_u$ estimated from a real data set and Lemma \ref{lemma_mfbm} to choose a  process in this class. Finally, we add the estimated mean function from the real data set, and thus define the simulated data generating process.

% !TeX root=bj-sample.tex

\subsection{Simulation design} % (fold)
\label{sub:simulation_design}

Our simulation study is based on the Household Active Power Consumption data set which was sourced from the UC Irvine Machine Learning Repository (\url{https://archive.ics.uci.edu/ml/datasets/Individual+household+electric+power+consumption}). This data set contains diverse energy related features gathered in a house located near Paris, every minute between December 2006 and November 2010. In total, this represents around $2$ million data points. We focus here on the daily voltage and we only consider the days without missing values in the measurements. The extracted data set contains $708$ voltage curves with an uniform common design with  $1440$ points, normalized such that $\mathcal{T} = [0, 1]$. We aim to simulate data sets using the data generating process defined in Section \ref{lemma_NicoH}, with a Hurst index function $H_t$ and a time deformation function $A_t$ estimated on the Power Consumption data set, to which we add a mean curve also fitted to the real data set.
For the fitted mean curve, we consider the  model 
\begin{equation}\label{setup_simu_mu}
\mu(t) = \beta_0 t + \sqrt{2}\sum_{1\leq k \leq 50}%^{50}
\left\{\beta_{1, k}\cos(2k\pi t) + \beta_{2, k}\sin(2k\pi t)\right\}, \quad t\in[0,1].
\end{equation}
The coefficients $\beta$ are obtained by LASSO regression with the \texttt{\textbf{R}} package \texttt{\textbf{glmnet}}. The outcomes are given by the 1440 values of the empirical mean of the $708$ curves, and the $t$ are sampled on the regular grid of 1440 points.  The regularity of the mean function is controlled using the penalty parameter $s$ of the \texttt{predict.glmnet} function.

\begin{figure}
		\vspace{0.2cm}
    \centering
    \begin{subfigure}[b]{0.3\linewidth}
        \centering
        \includegraphics[scale=0.3]{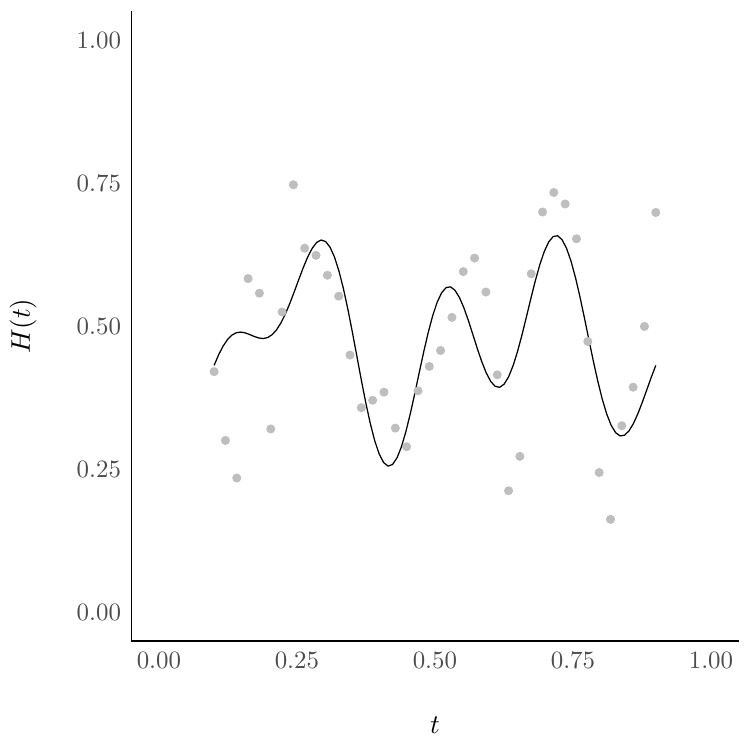}
        \caption{Estimation of $H(t)$}
        \label{fig:estim_H}
    \end{subfigure}
    \hfill
    \begin{subfigure}[b]{0.3\linewidth}
        \centering
        \includegraphics[scale=0.3]{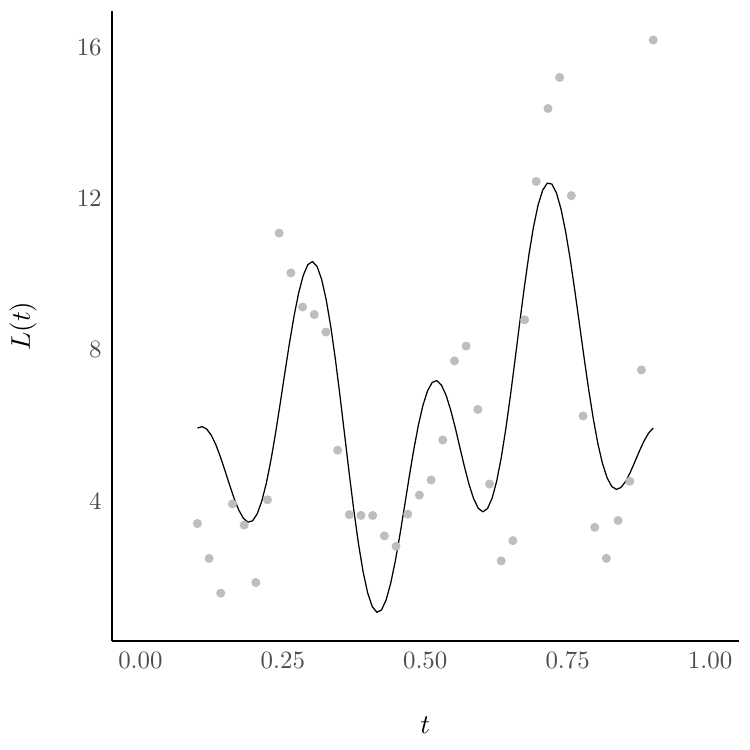}
        \caption{Estimation of $L(t)$}
        \label{fig:estim_L}
    \end{subfigure}
    \hfill
    \begin{subfigure}[b]{0.3\linewidth}
        \centering
        \includegraphics[scale=0.3]{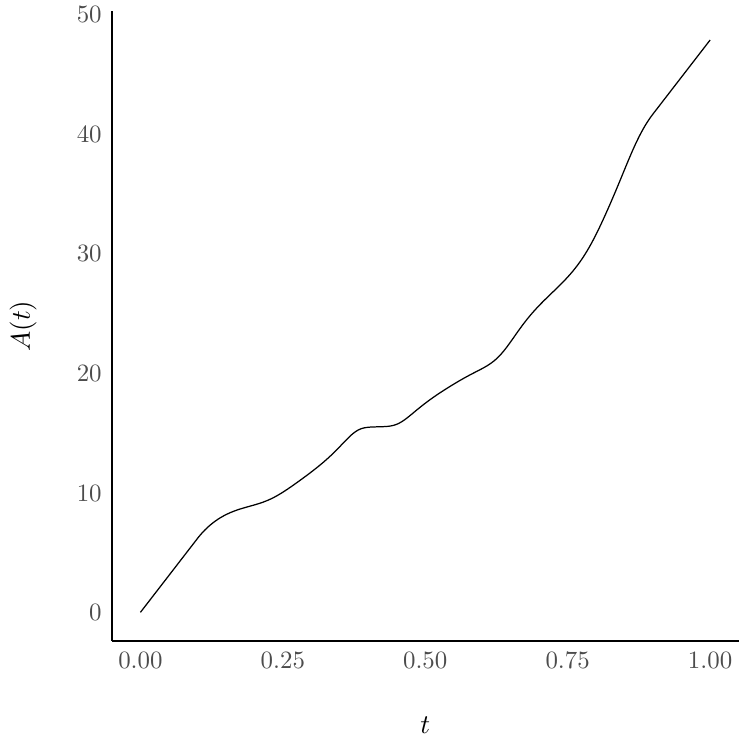}
        \caption{Estimation of $A(t)$}
        \label{fig:estim_A}
    \end{subfigure}
    \caption{Estimation of the different quantities for the data generating process.}
    \label{fig:estim_quantites}
\end{figure}

For the estimation of the Hurst index function $H_t$ and of the Hölder constant function $L_t$ using the Power Consumption data set, we apply  \eqref{eq:hat-h} and \eqref{eq:hat-l}, respectively.   The estimated values of $H_t$ and $L_t$ are smoothed using few functions from the Fourier basis. The resulting smoothed functions $H$ and $L$ are plotted in Figures \ref{fig:estim_H} and \ref{fig:estim_L}. The time deformation function $A(t)$ is then estimated using Lemma \ref{lemma_mfbm}, and the result is in Figure \ref{fig:estim_A}. Using these quantities, we estimate the covariance $C_A(\cdot, \cdot)$ of a MfBm process, as  defined in \eqref{eq:cov_mfbm_c}. Finally, to prevent each curve from starting from the same point, we add a random shift $X(0) \sim \mathcal{N}(0, \varpi^2)$. The final covariance of the process is thus given by $\Gamma(s, t) = \varpi^2 + C_A(s, t)$ for all $s, t \in \mathcal{T}$. We next generate samples of independent paths $X^{(i)}$ from the Gaussian process characterized by $\mu$ and $\Gamma$. Finally, to obtain the simulated functional data, we add  Gaussian noise of variance $\sigma^2$ at any observation time $\Tnm$.

\begin{figure}
    \centering
    \begin{subfigure}[b]{0.3\linewidth}
        \centering
        \includegraphics[scale=0.3]{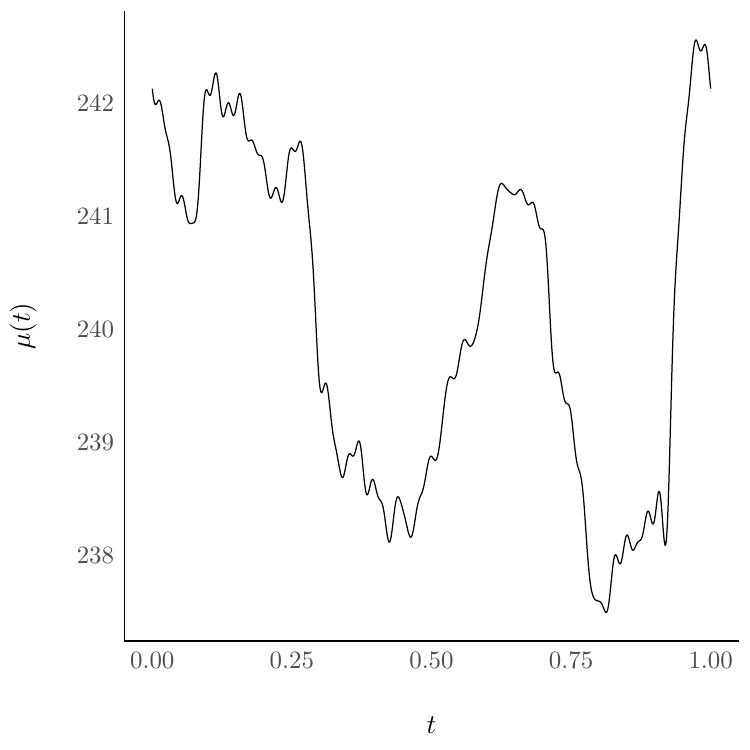}
        \caption{Mean curve $\mu(\cdot)$}
        \label{fig:mean_curve}
    \end{subfigure}
    \hfill
    \begin{subfigure}[b]{0.3\linewidth}
        \centering
        \includegraphics[scale=0.3]{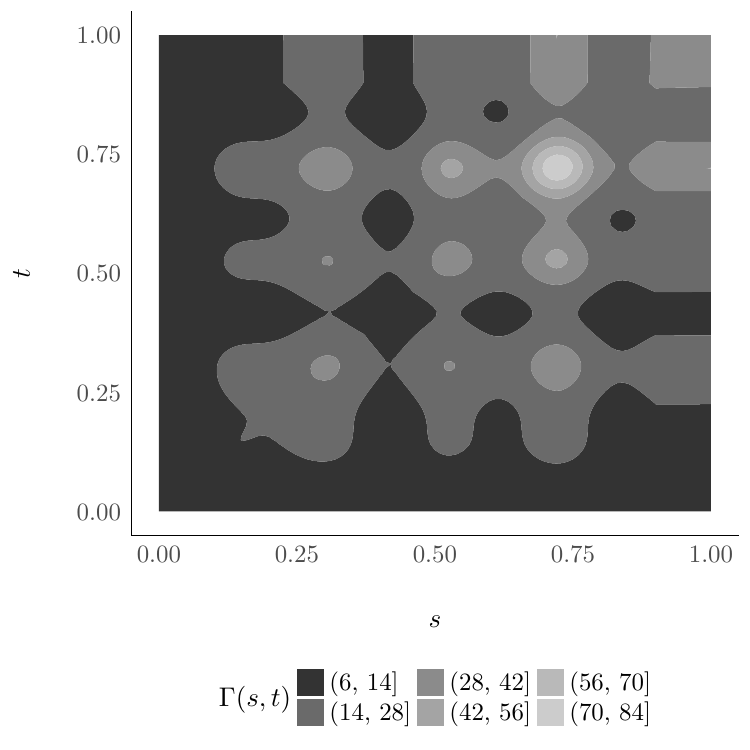}
        \caption{Covariance surface $\Gamma(\cdot, \cdot)$}
        \label{fig:cov_curve}
    \end{subfigure}
    \hfill
    \begin{subfigure}[b]{0.3\linewidth}
        \centering
        \includegraphics[scale=0.3]{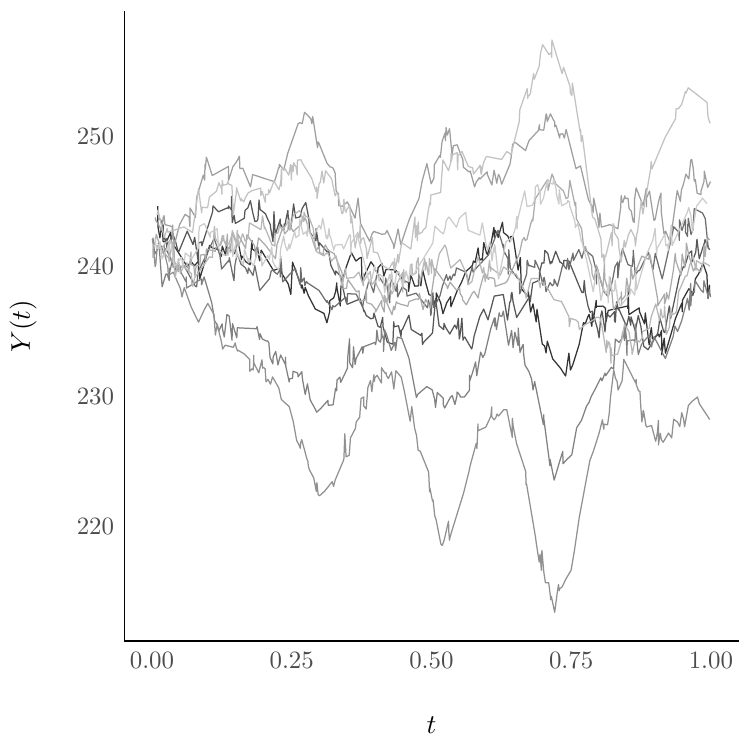}
        \caption{Noisy curves $Y^{(i)}$}
        \label{fig:sample_noisy_curve}
    \end{subfigure}
    \caption{Description of the simulated data set.}
    \label{fig:dataset_description}
\end{figure}

We consider eight experiments, each of them replicated 500 times. For each experiment, except specifically specified, we consider $N \in \{50, 100, 200\}$, $\mathfrak{m} \in \{20, 30, 40, 50\}$ and that the number of points per curve $M_i$ has a Poisson distribution with mean $\mathfrak{m}$. In \emph{Experiment 1}, we assume that the sampling points are uniformly distributed in $\mathcal{T}$, the standard deviation of the noise is $\sigma = 0.5$, the regularity of the mean function is $s = \exp(-6)$, the number of Fourier basis functions for the estimation of $H_t$ and $L_t$ is $9$, and $\varpi = 2.5$. All the other experiments are designed starting  from \emph{Experiment 1} and modifying one parameter at a time.  The mean and covariance functions corresponding to \emph{Experiment 1} are plotted in  Figure \ref{fig:mean_curve} and Figure \ref{fig:cov_curve}, respectively.
The noisy versions $Y^{(i)}$ of a random sample of ten curves $X^{(i)}$ generated according to \emph{Experiment 1} are plotted in Figure \ref{fig:sample_noisy_curve}.

In \emph{Experiment 2} and \emph{Experiment 3}, we consider $\sigma = 0.25$ and $\sigma = 1$, respectively. We set $s = \exp(-3)$ for \emph{Experiment 4} resulting in a smoother mean function $\mu$. We use only $7$ functions in the Fourier basis in \emph{Experiment 5}, that is a smoother estimation of $H_t$ and $L_t$. For \emph{Experiment 6}, the distribution of the sampling points is a mixture of beta distributions $0.5  \mathcal{B}(1, 2) + 0.5 \mathcal{B}(2, 1)$. For \emph{Experiment 7}, we set $\varpi = 1$. Finally, in \emph{Experiment 8}, we apply our approach to the case of differentiable trajectories that we obtain by integrating the sample paths generated as in \emph{Experiment 1}.
The results from \emph{Experiment 1} are presented below, those of the other seven experiments, and some additional implementation details, are provided in \cite{supplement}.  
An implementation of the methods used in all experiments is available as an \texttt{\textbf{R}} package on Github at the URL adress: \url{https://github.com/StevenGolovkine/funestim}.

% subsection simulation_design (end)

\subsection{Mean estimation} % (fold)
\label{sub:mean_estimation}

For the adaptive estimation of the mean curve, we first compute  $\widehat{H}_t$, according to  \eqref{eq:hat-h}, on a uniform grid of 20 points $t_2$ between $0.2$ and $0.8$, with $t_3 - t_1 = \Delta_* / 2   = \min(\exp(-\log(\mathfrak m)^{1/3}), 0.2) $. The local regularity being a local property, we constrain $\Delta_*$ to sufficiently small values. For each value of the $20$ estimates $\widehat{H}_t$, we compute the optimal bandwidths $h^*_\mu$ by minimization with respect to $h$ over a geometric grid $\mathcal H_N$ of $151$ points. We then estimate the mean function on $101$ regularly spaced points in $[0,1]$. The $101$ bandwidth values used for our estimator are then obtained from the $20$ optimal bandwidths $h^*_\mu$ by linear interpolation.

Our mean estimator, denoted $\widehat{\mu}_{GKP}$, is compared to that of \cite{cai2011}, denoted $\widehat{\mu}_{CY}$, and \cite{zhang_sparse_2016}, denoted $\widehat{\mu}_{ZW}$. To compute $\widehat{\mu}_{CY}$, we use the \texttt{smooth.splines} function in the \texttt{\textbf{R}} package \texttt{\textbf{stats}}, with the $M_1 + \ldots + M_N$ data points $(\Ynm , \Tnm)$. To obtain $\widehat{\mu}_{ZW}$, we use the \texttt{\textbf{R}} package \texttt{\textbf{fdapace}}, see  \cite{fdapace_b}. To compare the estimators, we use the integrated squared error (ISE) risk. For any $\varepsilon\in[0,1/2)$, if $f$ and $g$ are real-valued functions defined on $[0,1]$, let
\begin{equation}
\text{ISE}_{\varepsilon}(f, g) = \int_{ [\varepsilon, 1 - \varepsilon]} \{f(t) - g(t)\}^2 \mathrm{d}t.
\end{equation}
The integral is approximated by the trapezoidal rule with an equidistant  grid. 
For each configuration $(N, \mathfrak m)$, and each of the 500 samples, we compute   the ratios 
$$\frac{\text{ISE}_{\varepsilon}(\widehat{\mu}_{GKP}, \mu)}{\text{ISE}_{\varepsilon}(\widehat{\mu}_{CY}, \mu)} \quad\text{and}\quad \frac{\text{ISE}_{\varepsilon}(\widehat{\mu}_{GKP}, \mu)}{\text{ISE}_{\varepsilon}(\widehat{\mu}_{ZW}, \mu)},$$
and compare them to $1$. 

\begin{figure}
		\vspace{0.2cm}
    \centering
    \includegraphics[scale=0.38]{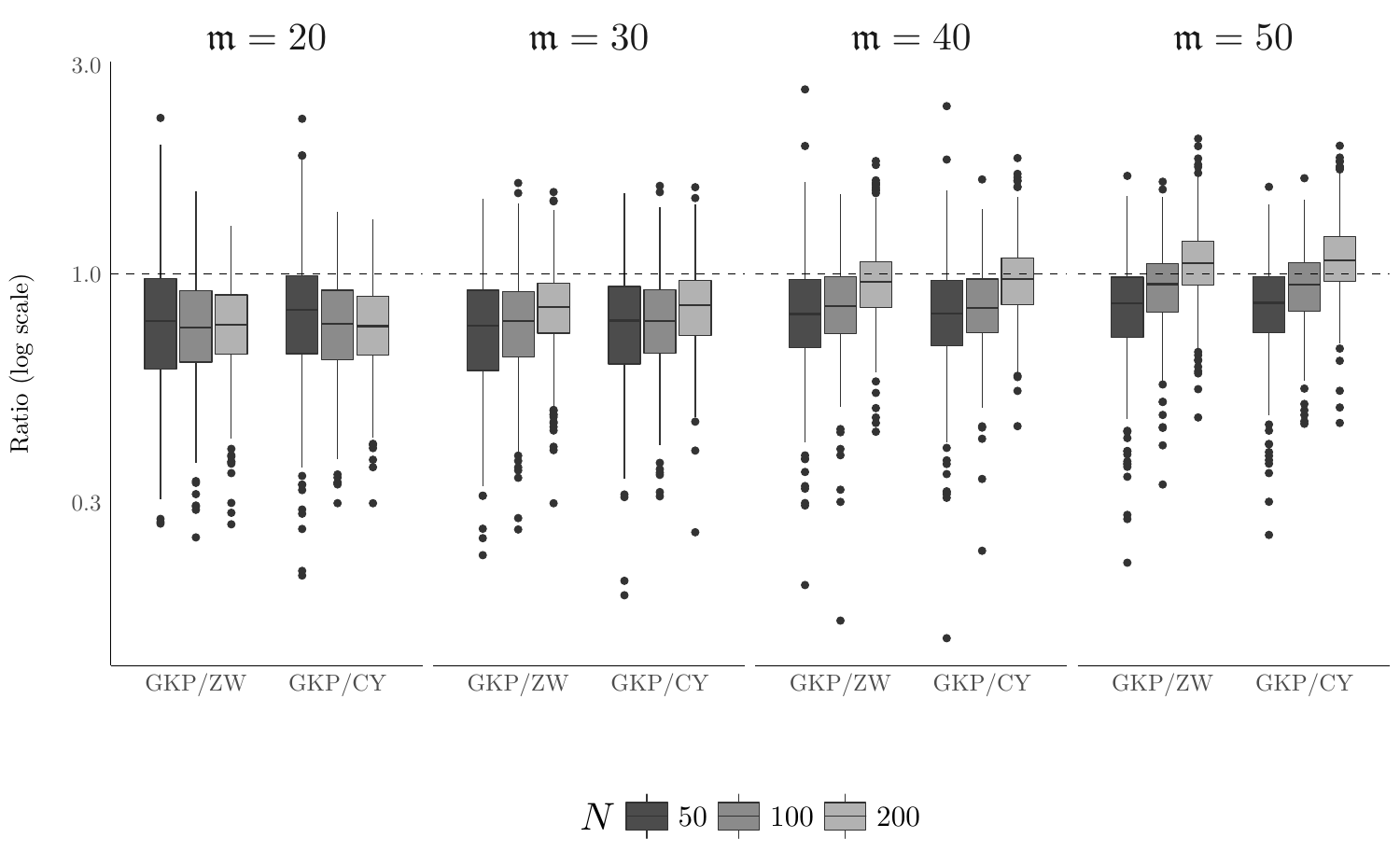}
    \caption{Results for the estimation of $\mu$ in \emph{Experiment 1}. The ratios are computed using $\text{ISE}_{0}$.}
    \label{fig:mean_true_exp1}
\end{figure}

The results for the $\text{ISE}_0$ ratios obtained in \emph{Experiment 1} are plotted in  Figure \ref{fig:mean_true_exp1}, on a logarithmic scale. To account for a possible boundary effect, we also computed the $\text{ISE}_{0.05}$ ratios, for which the results are similar, and reported in \cite{supplement}.
Our mean function estimator reveals good performance. Except for some cases where $N\mathfrak m$ is large,  our estimator outperforms the competitors. In those cases, the three estimators have similar performance. The fact that the advantage of our estimator vanishes when $N\mathfrak m$ is large could be explained by the fact that the approaches of \cite{cai2011} and \cite{zhang_sparse_2016} smooth over the pooled observations $(\Ynm , \Tnm)$. Similar conclusions are drawn from \emph{Experiments 2} to \emph{7}. In the setup with a more regular mean function (\emph{Experiment 4}), the advantage of our estimator diminishes.

% subsection mean_estimation (end)

\subsection{Covariance estimation} % (fold)
\label{sub:covariance_estimation}

For the adaptive estimation of the covariance function, we use the estimates $\widehat{H}_t$  computed for the mean function on the grid of 20 points between 0.2 and 0.8. For each of the 190 pairs $(s,t)$, $s< t$, on the grid, we compute the optimal bandwidths $h_\Gamma^*(s,t)$ by minimization over a logarithmic grid of $41$ points.
 %between $0.01$ and $0.1$.  
 We then estimate the covariance on a $101 \times 101$ regular grid. The $101 \times 101$ bandwidth values used for our estimator are obtained from the 190 optimal bandwidths $h^*_\Gamma(s,t)$ by symmetry and linear interpolation.

\begin{figure}
	\vspace{0.2cm}
    \centering
    \includegraphics[scale=0.38]{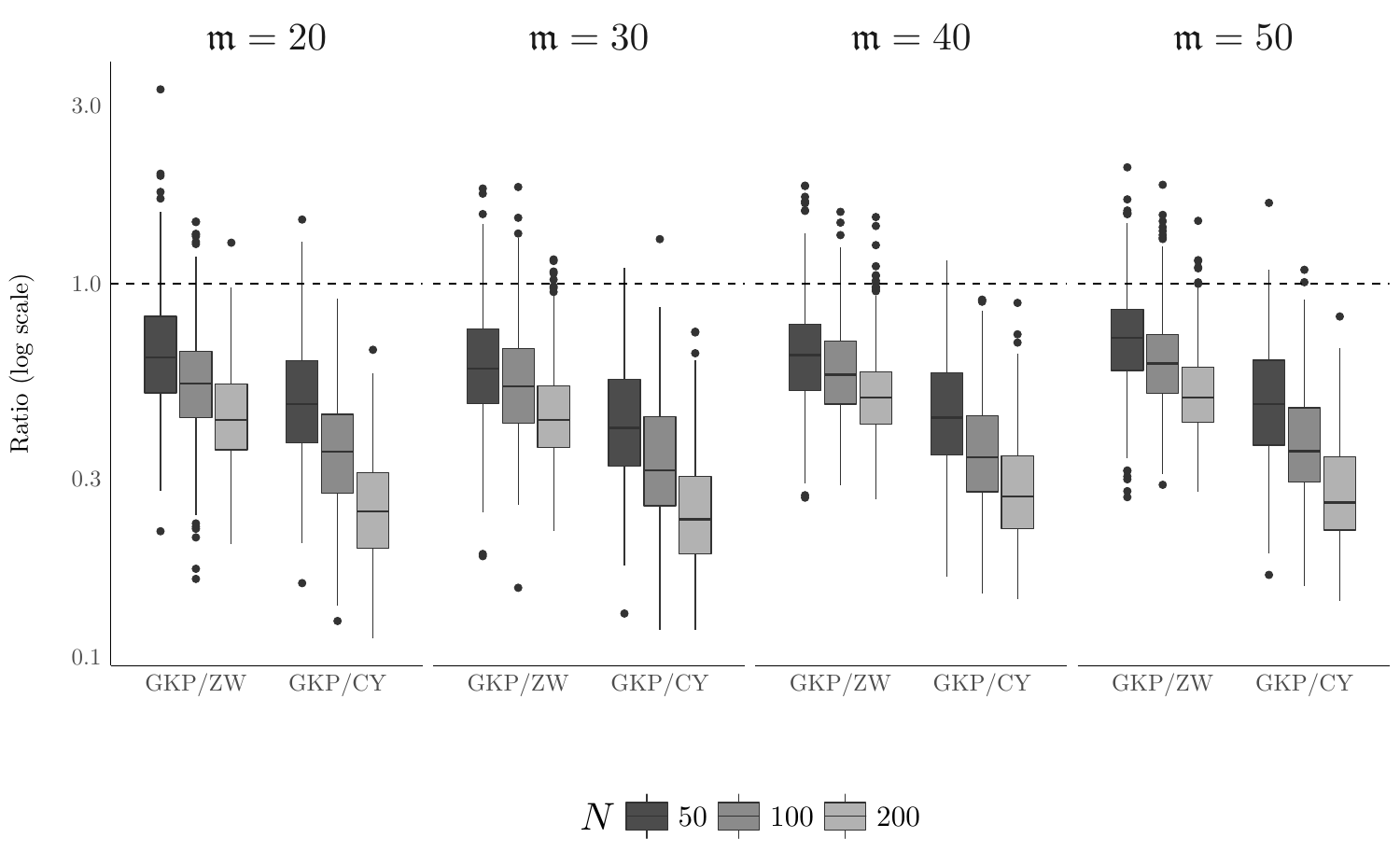}
    \caption{Estimation of $\Gamma$ in \emph{Experiment 1}. The ratios are computed using $\text{2-ISE}_{0}$.}
    \label{fig:cov_true_exp1}
\end{figure}

Our covariance estimator, denoted  $\widehat{\Gamma}_{GKP}$, is compared to that of \cite{cai2010}, denoted $\widehat{\Gamma}_{CY}$, and from \cite{zhang_sparse_2016}, denoted $\widehat{\Gamma}_{ZW}$. We compute $\widehat{\Gamma}_{CY}$ using the \texttt{\textbf{R}} package \texttt{\textbf{ssfcov}}, see  \cite{cai2010}. For $\widehat{\Gamma}_{ZW}$, we use the \texttt{\textbf{R}} package \texttt{\textbf{fdapace}}, see  \cite{fdapace_b}. To compare the accuracy of the estimators, we use the $2$-dimensional ISE risk. For any $ \varepsilon\in[0, 1/2)$, if $f$ and $g$ are real-valued functions defined on $[0, 1] \times [0, 1]$, let
\begin{equation}
\text{2-ISE}_{\varepsilon}(f, g) = \int_{[ \varepsilon, 1 - \varepsilon]}\int_{[ \varepsilon, 1 - \varepsilon]} \{f(s, t) - g(s, t)\}^2 \mathrm{d}s\mathrm{d}t,
\end{equation}
and the integral is approximated by the trapezoidal rule. For each configuration $(N,\mathfrak m)$, and each replication, we compute the $\text{2-ISE}_{\varepsilon}$'s with respect to the true covariance function $\Gamma$. We then compute the ratios 
$$\frac{\text{2-ISE}_{\varepsilon}(\widehat{\Gamma}_{GKP}, \Gamma)}{\text{2-ISE}_{\varepsilon}(\widehat{\Gamma}_{CY}, \Gamma)} \quad\text{and}\quad \frac{\text{2-ISE}_{\varepsilon}(\widehat{\Gamma}_{GKP}, \Gamma)}{\text{2-ISE}_{\varepsilon}(\widehat{\Gamma}_{ZW}, \Gamma)}.$$
The results for the ratios obtained with $\text{2-ISE}_0$ in \emph{Experiment 1} are plotted in  Figure \ref{fig:cov_true_exp1}, on a logarithmic scale. 
Those obtained with $\text{2-ISE}_{0.05}$, presented in \cite{supplement}, are similar. Our estimator shows better accuracy for estimating $\Gamma$ than $\widehat{\Gamma}_{ZW}$ and $\widehat{\Gamma}_{CY}$ in all cases considered. The advantage of our approach increases with $N$. 

% !TeX root=bj-sample.tex

%\newpage

\section{Discussion and conclusions}\label{sec:dics}

We propose new nonparametric estimators for the mean and covariance functions. They are built using a novel ``smoothing
first, then estimate'' strategy based on univariate kernel smoothing. The main novelty comes from the fact that the optimal bandwidths are selected by minimization of suitable penalized quadratic risks. The penalized risks for the mean and the covariance functions are quite similar and could easily be  built from data, and optimized on a grid of bandwidths. What distinguishes them from the usual sum between the squared bias and the variance, is a penalty for the fact that not all the curves have enough observation points to be included in the final estimator. Removing curves from the nonparametric estimators of the mean and covariance functions is an aspect which characterizes practically all smoothing-based approaches. Indeed, to entirely benefit from the replication feature of functional data, one has to determine the amount of smoothing for the mean and covariance functions estimation using all the curves. In this case, some curves could present too few observation points and thus will be dropped. This is more likely to happen in the so-called sparse regime. To the best of our knowledge, our bandwidth choice  is the first attempt to explicitly account for this aspect. We thus build estimators which achieve optimal rates of convergence in a completely adaptive, data-driven way. The theoretical results are derived under very mild conditions. In particular, the curves could be observed with heteroscedastic errors at discrete observations points. These points could be common to all curves or they could change randomly from one curve to another. In the case of the common observation points, our procedure automatically chooses between smoothing and interpolation, the latter being known to be rate optimal, but is not necessarily the best solution with finite samples. 

Our nonparametric estimation approach relies on a probabilistic concept of local regularity for the sample paths of the process generating the curves.  In some common examples, this local regularity is related to the polynomial decrease rate of the eigenvalues of the covariance operator, a characteristic of the data generating process widely used in the literature and usually supposed to be known.  The local regularity  also determines the regularity of the trajectories, the usual concept used in nonparametric regression.
It is well-known that the optimal rates, in the minimax sense, for estimating the mean and covariance functions, depend on the regularity of the paths. Moreover, the so-called sparse and dense regimes  in functional data analysis, are defined using the regularity of the trajectories, which is usually supposed to be known. We therefore consider a simple estimator of the local regularity of the process and use it to build our penalized quadratic risks. Applied to real data sets, the local regularity estimator reveals that the regularity of the trajectories could be quite far from what is usually assumed in the literature. However, in some applications, assuming smooth trajectories seems reasonable. The mean and covariance functions estimation approach based on local regularity extends to such situations.  In the case where the sample paths of $X$ admit derivatives up to the order, say $\alpha$, condition \eqref{eq:reg0} has to be stated for the increments of the $\alpha$th derivative of the sample path. \citet[App.~D]{golo2020} investigates this extension and proposes an estimator of $\alpha + H_t$, for which they derive a concentration bound. The mean and covariance functions can next be estimated using the estimates $\hXni_t$  built with local polynomial weights $W_{m}^{(i)}(t)$. The risk bounds derived in Section \ref{sec:LP_est} above can be extended to this case using standard arguments. See, for instance, \cite{tsybakov2009}. A discussion and an illustration of the performance of our adaptive estimation of the mean function with smooth sample paths is provided in \citet[Sec.~\ref{sub:case_of_differentiable_curves}]{supplement}.

Our regularity $H_t$ estimator is related to existing estimators in the stochastic process literature, which are generally based on noiseless observations from a single path. In the context of functional data, the estimators designed for a single path would require (very) large sets of design points on the curves. Investigating how to modify these estimators and use data points from several sample paths to reduce the number of design points required on each path is an interesting topic for future work.

The methodology we propose also allows for other types of adaptive results for the mean and covariance functions, such as convergence in distribution. A Central Limit Theorem (CLT) for the pointwise mean  estimator is provided in \cite{supplement}. In the context of functional data, the mean function may be more regular than the sample paths, but the limit case where the mean and the sample paths have the same regularity has to be considered from a minimax point of view. This means that from a minimax perspective, the rate of convergence of the distribution for estimating the mean function must take into account the regularity of the sample paths, otherwise the convergence of the distribution may fail due to the large bias.

Our method performs  well in simulations and outperforms the main competitors when  the mean and covariance functions have a regularity close to that of the trajectories. The approach is still satisfactory when these 
%mean or the covariance 
functions are more regular than the trajectories. The reason is that, in some sense, our nonparametric estimators are  close  to the empirical mean and covariance, respectively, which are the ideal estimators if the trajectories were observed at any point without error. In the case where the mean and covariance functions are smoother than the trajectories, our penalized quadratic risks should be built using the regularity of the mean or covariance functions, instead the regularity of the trajectories. However,  the estimation of the regularity of the mean or covariance function remains an open problem.

%%%%%%%%%%%%%%%%%%%%%%%%%%%%%%%%%%%%%%%%%%%%%%
%% Single Appendix:                         %%
%%%%%%%%%%%%%%%%%%%%%%%%%%%%%%%%%%%%%%%%%%%%%%
\setcounter{equation}{0}
\renewcommand{\theequation}{A.\arabic{equation}}
\begin{appendix}

% !TeX root=bj-sample.tex

\section{Technical details and proofs}

\begin{proof}[Details on (\ref{eq:opt_h_mean0c})]

Let
$\widetilde \mu_W (t;h) = \WN(t;h)^{-1} \sum_{i=1}^N  w_i (t;h)\Xp{i}_t $
be the unfeasible estimator of $\mu(\cdot)$ using only the curves for which $\hXni_t$ is well-defined. In the following,  we write $w_i$ and $\WN$ instead of $w_i (t;h) $ and $\WN(t;h) $, respectively.   By~\eqref{dec_s} and \eqref{mu1_d}, 
\begin{multline*}\label{eq:dec_emp_r}
\EEMT  [ \{  \widetilde \mu_N (t) - \widehat \mu_N (t;h) \}^2 ]
=    \EEMT\left[ \left\{ \widetilde \mu_N (t)  - \widetilde \mu_W  -   \frac{1}{\WN} \sum_{i=1}^N  w_i  \left(B^{(i)}_t +V^{(i)}_t \right)\right\}^2  \right] \\
\leq 2 
 \EEMT[ \{ \widetilde \mu_N (t)  - \widetilde \mu_W (t;h)\}^2  ] + 2\EEMT\left[ \left\{  \frac{1}{\WN} \sum_{i=1}^N  w_i \left(B^{(i)}_t +V^{(i)}_t \right)\right\}^2  \right]\eqqcolon 2E_1+2E_2.
\end{multline*}
Since 
$$
 \widetilde \mu_W (t;h)  - \widetilde \mu_N (t)  = \frac{1}{\WN  } \sum_{i=1}^N \left\{ X_t^{(i)} - \mu(t) \right\} \left\{w_i  - \WN  /N \right\},
$$
the trajectories of $X$ are drawn independently, and  independently of the $M_i$ and the $\Tnm$,   we have
\begin{equation}\label{eq:Ncal2b1}
E_1
 =\frac{\operatorname{Var}[X_t]}{\WN^2  } \sum_{i=1}^N  \left\{w_i    - \frac{\WN  }{N}  \right\}^2= q^2_3   \left\{  \frac{1}{\WN  } - \frac{1}{N } \right\} .
\end{equation}

For $E_2$,  let us first look at the bias part.   
By Theorem \ref{thm:estimation-alpha}, there   exists $\varrho >1$ such that the probability of the event $\{|\widehat H_t - H_t | > \log^{-\varrho}(\mathfrak m) \} $ is exponentially small.   Hence, by \eqref{mfka} and \eqref{mfkb}, we have $h^{2\widehat H_t } = h^{2H_t}\{1+o_{\mathbb P} (1) \}$, uniformly over the range $\mathcal H_N$  (i.e. the $o_{\mathbb P} (1)$ does not depend on $h$).  Next, 
in the case of non-differentiable sample paths,
by \ref{H:equivalent},
\begin{equation}\label{ez:equiv}
\EEMT\left[\left\{ \Xn(\Tnm)-   \Xp{i}_t \right\}^{2}   \right]=
\EE \left[ \left\{ \Xn(\Tnm)-   \Xp{i}_t \right\}^{2}  \Bigm| \mathcal  T_{obs}^{(i)} \right]
=\{1+o_{\mathbb P} (1) \} L_t^2 \left\lvert  (\Tnm-t)/h\right\rvert^{2H_t}.
\end{equation}
Similarly to \eqref{eq:reta_bias2} and \eqref{eq:reg01},  for the NW estimator $\hXni_t$ and $ \overline{C}$ defined in  \eqref{eq:qdf1},  we then have
\begin{multline}\label{eq:Ncal2z}
\EEMT\left[ \left\{  \frac{1}{\WN} \sum_{i=1}^N  w_i B^{(i)}_t \right\}^2  \right] 
   \leq   L_t^2h^{2H_t} %\{1+o_{\mathbb P} (1) \}   \times  
   \frac{1+o_{\mathbb P} (1)}{\WN}  \sum_{i=1}^N w_i \left\{ \sum_{m=1}^{M_i} W_{m}^{(i)}(t)  \times  \sum_{m=1}^{M_i}\left\lvert  \frac{\Tnm-t}{h}\right\rvert^{2H_t} W_{m}^{(i)}(t) \right\}  
\\  =L_t^2h^{2\widehat H_t} \times  \overline{C} (t;h,2\widehat H_t)  \times \{1+o_{\mathbb P}  (1) \} = L_t^2h^{2\widehat H_t}   \times \int \lvert u \rvert^{ 2\widehat H_t}K(u)du \times  \{1+o_{\mathbb P}  (1) \}   .
\end{multline}
Using the equivalent kernels, see \citet[Sec~3.2.2]{fan_local_1996}, the bound on the last line of the last display could be extended to the case of local linear estimators.

To complete the bound for $E_2$, note that by construction, $\EEMT[V_t^{(i)}  B_t^{(i)}] =0$,  
$$
\EEMT [ V_t^{(i)}  B_t^{(j)} ] =\EEMT [V_t^{(i)}  V_t^{(j)}  ] =0,\quad 1\leq i\neq j\leq N.
$$
The variance part in $E_2$ can be bounded as in  \eqref{eq:var_8}. Up to negligible terms, %we then have 
for the NW estimator,
\begin{equation}\label{eq:riskE3}
E_2 \leq  h^{2\widehat H_t}  L_t^2   \overline{C} (t;h, 2\widehat H_t)
+ \frac{\sigma^2_{\max}}{\WN^2}\sum_{i=1}^N  w_i\mathcal N^{-1}_{i} (t;h) =  q_1^2 h^{2\widehat H_t} + q_2^2\mathcal N^{-1}_{\mu} (t;h) .
\end{equation}
\end{proof}

%\section{Proofs}

\begin{proof}[Proof of Theorem \ref{thm:optim_mu}]
First,  note that if $\WN(t;h)=0$, then necessarily $\mathcal N_\mu (t;h) =0$. Moreover, it will be shown below that $\inf_{h\in\mathcal H_N}\EE[\WN(t;h)]$ stays away from zero, and $\WN(t;h)$ uniformly concentrates to $\EE[\WN(t;h)]$ with high probability. We therefore, in the following, work on the event $\{\inf_{h\in\mathcal H_N}\WN(t;h)\geq 1\}$. 
First, let us prove that a constant $C>0$ exists such that 
\begin{equation}\label{th_q1}
0\leq \WN(t;h)^{-1} - N^{-1} \leq C \max\left\{ h^{2 H_t} ,\mathcal N^{-1}_\mu (t;h) \right\} \{1+r_N(h)\}, 
\end{equation}
with $\sup_{h} \lvert r_N(h) \rvert = o_{\mathbb P}(1)$.
Property \eqref{th_q1} is implied by the following:  constants $\mathfrak  c_1, \mathfrak  c_2>0$ exist such that 
\begin{equation}\label{th_q1b}
 \mathfrak  c_1N \{1+o_{\mathbb P}(1)\} \leq \inf_{h\in\mathcal H_N} \frac{\WN(t;h) } { \min\left\{1 ,  \mathfrak m h \right\}}
  \leq  \sup_{h\in\mathcal H_N} \frac{\WN(t;h) } { \min\left\{1 ,  \mathfrak m h \right\}} \leq \mathfrak  c_1^{-1} N \{1+o_{\mathbb P}(1)\} . 
\end{equation}
and
\begin{equation}\label{th_q1a}
\mathfrak c_2\{1+o_{\mathbb P}(1)\} \leq \inf_{h\in\mathcal H_N}\frac{N\mathfrak m h}{\mathcal N_\mu (t;h) } 
	\leq \sup_{h\in\mathcal H_N}\frac{N\mathfrak m h}{\mathcal N_\mu (t;h) }
\leq \mathfrak c_2^{-1}\{1+o_{\mathbb P}(1)\}.
\end{equation}
Indeed, \eqref{th_q1b} and \eqref{th_q1a} imply  
\begin{align*}%\label{th_q1c}
 \WN(t;h)^{-1} - N^{-1} &\leq  \mathfrak  c_1^{-1} \max\left\{0,  (N \min\left\{1 ,  \mathfrak m h \right\})^{-1} 
- N^{-1}  \right\}\{1+o_{\mathbb P}(1)\}  \\ &\leq  \mathfrak  c_1^{-1}    \mathfrak  c_2     \max\left\{ h^{2 H_t} ,   \mathcal N_\mu (t;h)^{-1 }\right\} \{1+o_{\mathbb P}(1)\}, 
\end{align*} 
with the $o_{\mathbb P}(1)$ terms uniform with respect to $h\in\mathcal H_N$. 
The justification of \eqref{th_q1b} and \eqref{th_q1a} is provided in \cite{supplement}. 
Let us provide a brief insight on how \eqref{th_q1b} is obtained.  Here, $\WN(t;h)$  is a Binomial variable with  $N$ trials and the success probability is a non-decreasing function of $h$. The property \eqref{th_q1b} follows by suitably bounding $\EE[\WN(t;h)]$ and using   Chernoff's  inequality  on a grid of points in $\mathcal H_N$. The uniformity with respect to all $h\in\mathcal H_N$ is obtained using the monotonicity of $\WN(t;h)$ with respect to $h$.   
Finally, to complete the proof in the independent design case, it suffices first to notice that, from above, we can also deduce 
\begin{equation}\label{eq:not_easy_app}
  \min_{\mathcal H_N}   \left\{ h^{2 H_t} + \mathcal N^{-1}_\mu (t;h) \right\}  \asymp_\PP  \min_{\mathcal H_N}   \left\{ h^{2 H_t} + (N\mathfrak m h)^{-1}\right\}.
\end{equation}
The minimum on the RHS is attained by $h$ with the rate $(N\mathfrak m )^{-1/\{2 H_t+1\}}$,  and \eqref{th_q1}  guarantees that $h_\mu^*$ has the same rate with probability tending to 1.
Next, by \eqref{mfka} and \eqref{mfkb},  uniformly over $h\in\mathcal H_N$, we have $h^{2\widehat H_t} =  h^{2H_t} \{1+o_{\mathbb P} (1)\}$. The rate of $\widehat \mu_N^* (t) - \widetilde \mu_N (t) $ follows. For the rate of  $\widehat \mu_N^* (t) - \mu (t) $, we simply add the parametric rate of  $\widetilde \mu_N (t) - \mu (t) $.

With a common design, $\WN(t;h)$ can only take the values 0 or $N$.  The penalty introduced by $\WN(t;h)^{-1}- N^{-1} $  thus constrains the bandwidth to be greater than or equal to the lengths of the intervals $[T_m^{(i)}, T_{m+1}^{(i)}]$ including $t$. By condition \eqref{mfkb_2}, this means that the rate of convergence of  $\widehat \mu_N^* (t) - \widetilde \mu_N (t) $ could not be faster than $O_{\mathbb P} (  \Mmu^{- H_t }) $. This aspect is automatically included in the definition of $\mathcal{R}_\mu(t;h)$ because, under the constraint $\mathfrak m h \geq c_L/C_U$, using \eqref{eq:not_easy_app},
$$
  \min_{\mathcal H_N}  \left\{ h^{2H_t} + \mathcal N^{-1}_\mu (t;h) \right\}   
  \asymp_\PP 
   \min_{\mathcal H_N, \mathfrak m h \geq c_L/C_U}   \left\{ h^{2H_t}+ (N\mathfrak m h)^{-1}\right\}  
  \asymp 
 \max\left\{ \Mmu^{- \; 2H_t }, (N\mathfrak m)^{2H_t/(2H_t+1)}\right\}.  
$$
Finally, $H_t$ can be replaced by $\widehat H_t$ using again  $h^{2\widehat H_t} =  h^{2H_t} \{1+o_{\mathbb P} (1)\}$. \end{proof}

\begin{proof}[Proof of Theorem \ref{thm:optim_gamma}]
	For simplicity, we assume $\widehat \Gamma_N^* $ is built with the uniform kernel. Since the $s\neq t$ are fixed,  we can consider that $\sup \mathcal H_N< \lvert s-t \rvert/2$. We can also assume  $c_L\mathfrak m \geq 2$. Similarly to the case of the mean function estimation, we show that, 	in the independent design case, 
	\begin{equation}\label{th_q1ab_cov}
	\WN(s,t;h) \asymp_\PP \min \left\{1,  (\mathfrak m h)^2\right\} \quad \text{and} \quad 
	\frac{1}{\mathcal N_\Gamma (t|s;h) } \asymp_\PP
	\frac{\min \{1,  (\mathfrak m h)^{-1}\}}{N\min \{ 1,  (\mathfrak m h)^{2}\}} 
	=\frac{N^{-1}}{  \min \{ \mathfrak m h,  (\mathfrak m h)^{2}\} }, 
\end{equation}
	uniformly with respect to $h\in\mathcal H_N$. Here the uniform $\asymp_\PP$ is a short notation for inequalities as in \eqref{th_q1b} and \eqref{th_q1a}. As a consequence of \eqref{th_q1ab_cov}, for which 
the justification  is given in \cite{supplement}, we get
		\begin{equation}\label{th_qq1}
		\WN(s,t;h)^{-1}-  N^{-1} \leq  \min\left[\max\left\{ h^{2 H_s} ,\mathcal N^{-1}_\Gamma (s \mid t;h) \right\},  \max\left\{ h^{2 H_t} ,\mathcal N^{-1}_\Gamma (t \mid s;h) \right\} \right]  O_{\mathbb P}(1),
	\end{equation}
	with the $O_{\mathbb P}(1)$ term uniform with respect to $h\in\mathcal H_N$. From \eqref{th_q1ab_cov} and \eqref{th_qq1}, we deduce that, 
	\begin{equation}\label{eq:not_easy_b}
		\min_{\mathcal H_N}   \left\{  \mathcal R_\Gamma  (s \mid t;h)  + \mathcal R_\Gamma  (t \mid s;h) \right\}\asymp _{\mathbb P} \min_{\mathcal H_N} \left\{  h^{2 H_s}+ h^{2 H_t} + (N\min \{\mathfrak mh, (\mathfrak mh)^2\} )^{-1}\right\}.
	\end{equation}
Then, the rate of $h^*_\Gamma$ follows. Plugging $h^*_\Gamma$ into the risk bound of 	$ \widehat \gamma_N (s,t;h)$, we get 
\begin{equation}\label{eq:easy_q}
\EEMT \left[ \left\{   \widehat \gamma_N (s,t;h^*_\Gamma) - \widetilde \gamma_N (s,t)  \right\}^2 \right] =   O_{\mathbb P} \left(  (N\Mmu^2)^{- \; \frac{ 2H(s,t)}{2 \{H(s,t)+1\}} } +   (N\Mmu)^{- \; \frac{ 2H(s,t)}{2 H(s,t)+1} } \right).
\end{equation}
For the rate of $\EEMT[ \{   \widehat \gamma_N (s,t;h^*_\Gamma) -   \gamma (s,t)   \}^2 ] $, it suffices to add $O_\PP (N^{-1})$ to the right hand side in \eqref{eq:easy_q} to account for  $\EEMT[ \{   \widetilde \gamma_N (s,t) -   \gamma (s,t)   \}^2 ] $. The rate of $ \widehat \gamma_N (s,t;h^*_\Gamma) -   \gamma (s,t) $ follows by Markov inequality. 
To get the rate \eqref{eq:rate_G*}, it remains to show that  $\{ \widehat \Gamma_N (s,t;h) -  \Gamma (s,t)\} - \{   \widehat \gamma_N (s,t;h) - \gamma (s,t)   \}$ is negligible, uniformly with respect to $h\in\mathcal H_N$.  The justification is given  in \cite{supplement}.

From above and the arguments used in the proof for the mean, with common design we have
	\begin{align}
		\widehat \Gamma_N^* (s,t) - \Gamma (s,t) &= O_{\mathbb P} \left( \max\left\{ (N\Mmu^2)^{- \; \frac{ H(s,t)}{2 \{H(s,t)+1\}} } ,  (N\Mmu)^{- \; \frac{ H(s,t)}{2 H(s,t)+1} } , \Mmu^{-  H (s,t) }\right\}+ N^{-1/2}\right) \\
		& = O_{\mathbb P} \left(  \Mmu^{-  H (s,t)}\!+ N^{-1/2}\right) ,
	\end{align}
	with the last equality implied by \eqref{zdwr}. 
	The rate of $\widehat \Gamma_N^* (s,t) $ in \eqref{eq:rate_CG*} follows.
	\end{proof}

\begin{proof}[Details on \eqref{Dtilde_d}]
Let 
$
\widetilde D_t (\mathfrak  u_1 ,  \mathfrak {u}_2) = \widetilde \Gamma_N ( t- \mathfrak  u_1 , t+ \mathfrak {u}_2) -  \widetilde \Gamma_N (t,t).
$
For the bound in \eqref{Dtilde_d}, we use the assumptions: $\sup_{t\in\mathcal T} \EE (X_t^4) <\infty$ and 
 a constant $c$ exists such that 
\begin{equation}\label{eq:mom4}
\EE [\{X_s- X_t\}^4]  \leq c \EE^2 [\{X_s- X_t\}^2],\quad s, t\in\mathcal T.
\end{equation}
We then have
\begin{align*}
\EE [\widetilde D_t (\mathfrak  u_1 ,  \mathfrak {u}_2)^2] &\leq 2  
\EE \left[ \left\{\frac{1}{N} \sum_{i=1}^N \left(\{ \Xp{i}_t \}^2 - \Xp{i} _{t- \mathfrak  u_1} \Xp{i} _{t+ \mathfrak {u}_2} \right) \right\}^2\right]  \\ 
&\quad+ 2 
\EE 
\left[ 
\left(  \!
\left\{\frac{1}{N} \sum_{i=1}^N \Xp{i}_t \right\}^2 
  - \left\{ \frac{1}{N} \sum_{i=1}^N \Xp{i} _{t- \mathfrak  u_1} \right\}  \left\{ \!\frac{1}{N} \sum_{i=1}^N \Xp{i} _{t+ \mathfrak {u}_2} \!\right\} 
\right)^{\!2} 
\right]  = : 2D_1+2 D_2.
\end{align*}
By \ref{H:equivalent}, \eqref{eq:mom4}, and Jensen and Cauchy-Schwarz inequalities, a constant $C_1 $ exists, depending on $L_t$,  $S$, and $c$ appearing in \eqref{eq:mom4}, such that
$
D_1 \leq \EE [    \{ X_t^2  - X _{t- \mathfrak  u_1} X _{t + \mathfrak  u_2}    \}^2] \leq  C_1   \mathfrak d_t^{2H_t},
$
provided $  0\leq \mathfrak  u_1, \mathfrak  u_1 \leq \mathfrak d_t$. 
On the other hand, by similar arguments, 
\begin{multline*}
D_2\leq 2  \EE^{1/2} 
\left[ 
\left(\frac{1}{N} \sum_{i=1}^N \left\{  \Xp{i}_t - \Xp{i}_{t- \mathfrak  u_1} \right\}\right) ^{4}
\right] \EE^{1/2} \left[ 
\left\{\frac{1}{N} \sum_{i=1}^N  \Xp{i}_t  \right\} ^{4} \right]  \\
+2   \EE^{1/2} 
\left[ \left\{\frac{1}{N} \sum_{i=1}^N  \Xp{i}_{t- \mathfrak  u_1}\right\} ^{4}\right] \EE^{1/2} 
\left[ \left(\frac{1}{N} \sum_{i=1}^N \left\{  \Xp{i}_t -  \Xp{i}_{t+ \mathfrak  u_2} \right\}\right)^{4}
\right]\leq C_2  \mathfrak d_t^{2H_t} ,
\end{multline*}
for some constant $C_2 $. 
Gathering facts, we deduce  \eqref{Dtilde_d}. 
\end{proof}

%\medskip

\begin{proof}[Proof of Lemma \ref{lemma_mfbm}]
By construction, $\EE\left[W_A(A(t))\right]=0 $, and the  covariance function of $X$ is 
\begin{equation}\label{eq:cov_A}
\operatorname{Cov}_A(s,t) =
 D(H_{s} ,H_{t} )\left[ A(s)^{H_{s}+H_{t}} +  A(t)^{H_{s}+H_{t}} - \lvert A(t)-A(s)\rvert^{H_{s}+H_{t}}\right] ,\quad s, t\geq 0.
\end{equation}
Moreover, we show  in \cite{supplement} that, for any $t$ and $u, v \in\Ostar(t)$, we have
$\EE[ (X_u - X_v)^2 ] \approx \{A^\prime (t)\}^{2H_t} \lvert u-v \rvert^{2H_t}$.
To match  \ref{H:equivalent}, we define $A(\cdot)$ such that $\{A^\prime (t) \}^{H_{t}} = L_t $, and thus
$
A(t) = A(0) + \int_{0}^t L_s^{1/H_{s}} ds.
$
\end{proof}

\end{appendix}
%%%%%%%%%%%%%%%%%%%%%%%%%%%%%%%%%%%%%%%%%%%%%%
%% Multiple Appendixes:                     %%
%%%%%%%%%%%%%%%%%%%%%%%%%%%%%%%%%%%%%%%%%%%%%%

%\begin{appendix}
%\section{???}
%
%\section{???}
%
%\end{appendix}

%%%%%%%%%%%%%%%%%%%%%%%%%%%%%%%%%%%%%%%%%%%%%%
%% Support information, if any,             %%
%% should be provided in the                %%
%% Acknowledgements section.                %%
%%%%%%%%%%%%%%%%%%%%%%%%%%%%%%%%%%%%%%%%%%%%%%
\section*{Acknowledgements}
The authors thank the reviewers for their comments, which helped to improve the manuscript. The authors  thank Groupe Renault and the ANRT (French National Association for Research and Technology) for their financial support via the CIFRE convention No.~2017/1116.

%%%%%%%%%%%%%%%%%%%%%%%%%%%%%%%%%%%%%%%%%%%%%%
%% Funding information, if any,             %%
%% should be provided in the                %%
%% funding section.                         %%
%%%%%%%%%%%%%%%%%%%%%%%%%%%%%%%%%%%%%%%%%%%%%%
\section*{Funding}

S. Golovkine was partially supported by Science Foundation Ireland under Grant No.~19/FFP/7002 and co-funded under the European Regional Development Fund.

%%%%%%%%%%%%%%%%%%%%%%%%%%%%%%%%%%%%%%%%%%%%%%
%% Supplementary Material, including data   %%
%% sets and code, should be provided in     %%
%% {supplement} environment with title      %%
%% and short description. It cannot be      %%
%% available exclusively as external link.  %%
%% All Supplementary Material must be       %%
%% available to the reader on Project       %%
%% Euclid with the published article.       %%
%%%%%%%%%%%%%%%%%%%%%%%%%%%%%%%%%%%%%%%%%%%%%%
\section*{Supplementary Material}

In the Supplement, we provide some additional technical arguments, proofs and simulation results. In Section \ref{secA:comp_p}, we provide details for proofs presented in the main text. In Section \ref{sec:details_sm}, we provide details on some quantities and equations from the main text. A CLT for the pointwise  estimator of the mean function is presented in Section \ref{sec:TCL}.  Additional simulation results, including a case with smooth  paths,   are gathered in Section \ref{sub:simulation_design_SM}.

%%  MR numbers will be added by VTeX.                      %%
%%                                                         %%
%%  Use \cite{...} to cite references in text.             %%
%%                                                         %%
%%%%%%%%%%%%%%%%%%%%%%%%%%%%%%%%%%%%%%%%%%%%%%%%%%%%%%%%%%%%%

%% if your bibliography is in bibtex format, uncomment commands:
\bibliographystyle{apalike}
\bibliography{ref_short}

%% or include bibliography directly:
% \begin{thebibliography}{}
% \bibitem[\protect\citeauthoryear{???}{???}]{b1}
% \end{thebibliography}

\makeatletter\@input{main_supp_ref.tex}\makeatother

\end{document}

\maketitle
\begin{abstract}
In this supplementary material, we provide some additional technical arguments, proofs and simulation results. In Section \ref{secA:comp_p}, we provide details for proofs presented in the main text. In particular, details on the proofs of Theorems \ref{thm:optim_mu} and \ref{thm:optim_gamma} are given.
Section \ref{sec:details_sm} contains details on some quantities and equations from the main text, while in Section \ref{sec:TCL} a Central Limit Theorem (CLT) for the pointwise estimator of the mean function is proved. Additional simulation results are gathered in Section \ref{sub:simulation_design_SM}.
\end{abstract}

\tableofcontents

%%%%%%%%%%%%%%%%%%%%%%%%%%%%%%%%%%%%%%%%%%%%%%
%%%% Main text entry area:

%%%%%%%%%%%%%%%%%%%%%%%%%%%%%%%%%%%%%%%%%%%%%%
%% Single Appendix:                         %%
%%%%%%%%%%%%%%%%%%%%%%%%%%%%%%%%%%%%%%%%%%%%%%
% \begin{appendix}
% %\section*{???}%% if no title is needed, leave empty \section*{}.
% \end{appendix}
%%%%%%%%%%%%%%%%%%%%%%%%%%%%%%%%%%%%%%%%%%%%%%
%% Multiple Appendixes:                     %%
%%%%%%%%%%%%%%%%%%%%%%%%%%%%%%%%%%%%%%%%%%%%%%
%\begin{appendix}
%
\renewcommand{\theequation}{SM.\arabic{equation}}
% !TeX root=bj-sample-SM.tex

\section{Complements for the proofs}\label{secA:comp_p}

Below, $c,c_1, C, C_1,\mathfrak c, \mathfrak c_1, \ldots$ are constants with possibly different values at different occurrences. Moreover, the symbol $\lesssim$ (resp. $\gtrsim$) means that the left side is bounded above  (resp.  below) by a positive constant times the right side. Below, we use the notation $\EE_M[\cdot] = \EE[\cdot\mid M_1,\ldots,M_N]$.

\begin{proof}[Complements for the proof of Theorem \ref{thm:optim_mu}] %4.1]
We provide here a formal justification for the following properties:
two constants $0< \mathfrak  c_1, \mathfrak  c_2<1$ exist such that 
\begin{equation}\label{th_q1b_SM}
 \mathfrak  c_1N \{1+o_{\mathbb P}(1)\}  \leq \inf_{h\in\mathcal H_N} \frac{\WN(t;h) } { \min\left\{1 ,  \mathfrak m h \right\}}\leq  \sup_{h\in\mathcal H_N} \frac{\WN(t;h) } { \min\left\{1 ,  \mathfrak m h \right\}}
 \leq  \mathfrak  c_1^{-1} N \{1+o_{\mathbb P}(1)\} ,
\end{equation}
and
\begin{equation}\label{th_q1a_SM}
	\mathfrak c_2\{1+o_{\mathbb P}(1)\} \leq \inf_{h\in\mathcal H_N}\frac{N\mathfrak m h}{\mathcal N_\mu (t;h) }
	\leq \sup_{h\in\mathcal H_N}\frac{N\mathfrak m h}{\mathcal N_\mu (t;h) }
	\leq \mathfrak c_2^{-1}\{1+o_{\mathbb P}(1)\} 	,
\end{equation} 
and of \eqref{eq:not_easy_app} in the main text.   In reply to a Reviewer's remark, we prove \eqref{th_q1b_SM} and \eqref{th_q1a_SM} in a slightly more general framework. 
Let 
$$
\overline{M }= \frac{1}{N} \sum_{i=1}^N  M_i, \quad\text{such that } \mathfrak m = \EE(\overline{M }).
$$ 
The $M_i$ are independent, but we do not need to impose them to have the same distribution. However, for simplicity, we still assume \eqref{mfkb}. For each $1\leq i\leq N$, we denote by $g_i$, the density of the independent variables $\Tnm\in\mathcal T$, $1\leq m \leq M_i$. Moreover, the variables $\Tnm$ are drawn independently for each curve $i$.  Assume that positive constants $C_{g,L}, C_{g,U}>0$ exist such that, for all $t\in\mathcal T$ and for all $1\leq i \leq N,$
\begin{equation}\label{bgfm1}
C_{g,L} \leq g_i(t) \leq C_{g,U},
\end{equation}
and all $g_i$ are Hölder  continuous on $\mathcal T$ with the same exponent and Hölder constant. We thus allow the observation times  $\Tnm$ to be drawn independently with different distributions for different curves. 
Let 
\begin{equation}\label{def_pig}
	p_i(t;h)  =\int_{t-h}^{t+h}g_i(u)du .
\end{equation}
Under the conditions on the bandwidth  range $\mathcal H_N$ and the $g_i$, we have  $p_i(t;h)= 2hg_i(t)\{1+o(1)\}$, uniformly with respect to $h$ and $i$. 

To show the lower bound  in \eqref{th_q1a_SM}, recall that, with the NW estimator 
\begin{equation}\label{eq:Ncal1_SM}
\mathcal N_{i} (t;h) =  \frac{w_i (t;h)}{\max_{1\leq m\leq M_i} \lvert W_{m}^{(i)}(t;h) \rvert } \quad \text{and} \quad
 \mathcal N_{\mu} (t;h)^{-1} =  \frac{1}{\WN^2 (t;h)}\sum_{i=1}^N  \frac{w_i (t;h) }{\mathcal N_{i} (t;h)} .
\end{equation}
(Recall that the rule $0/0=0$ applies for $w_i/\mathcal N_i$). 
We simplify the notation in the following: $\mathcal N_{\mu} (t;h)$, $\mathcal N_{i} (t;h)$  and $w_i (t;h) $ become $\mathcal N_{\mu} $, $\mathcal N_{i}$   and $w_i$, respectively. Moreover, for simplicity, we assume that the NW kernel is built with the uniform density. The general case can be handled similarly using a positive lower bound for the kernel $K$ on a non-degenerate sub-interval of $[-1,1]$. With a uniform kernel, we have 
$$
\mathcal N_{i} = \sum_{m=1}^{M_i} \mathbf{1}\left\{ \lvert \Tnm-t \rvert \leq h\right\}.
$$ 
By Cauchy-Schwarz inequality, 
\begin{equation}\label{zebn}
\frac{1}{ \mathcal N_{\mu}} =  \frac{1}{\WN^2(t;h)}\sum_{i=1}^N  \frac{w_i  }{\mathcal N_{i}}\geq  \frac{1}{S_N(t;h)} \quad \text{with} \quad S_N(t; h) = \sum_{i=1}^N \mathcal N_i .
\end{equation}
Note that $S_N(t;h)$ is a sum of $M_1+\ldots+M_N$ independent Bernoulli variables with parameters
$$
\underbrace{p_1(t;h),\ \ldots \ ,p_1(t;h)}_{M_1 \text{ times}},\ \ldots \ , \underbrace{p_N(t;h),\ \ldots \ ,p_N(t;h)}_{M_N \text{ times}}.
$$
We have
\begin{equation}\label{zebn2}
C_{g,L} N\mathfrak m h \times   \{1+o(1)\} \leq \EE[S_N(t;h)] = \sum_{i=1}^N p_i(t;h) \EE(M_i ) \leq  C_{g,U} N\mathfrak m  h \times \{1+o(1)\}.
\end{equation}
Recall that we impose $N\mathfrak m \times \min \mathcal H_N \rightarrow \infty$. 
By Chernoff's inequality, for any $0 \leq \delta < 1$, 
$$
\mathbb P \left[ \left\lvert \frac{S_N(t;h) }{\EE[S_N(t;h)]} - 1\right\rvert > \delta  \right] \leq 2 \exp(-\delta^2 C_{g,L} N\mathfrak m \min \mathcal H_N  /3).
$$
We can choose  $\delta$ such that 
\begin{equation}\label{eq:choice_delta}
\delta^2 = C_\delta \frac{\log(N\mathfrak m) }{N\mathfrak m \min \mathcal H_N } ,
\end{equation}
with $C_\delta$ some large constant. 
If $h_1,\ldots, h_J$ is an equidistant grid on $\mathcal H_N$ of $J$ points, with $N\mathfrak m \leq J<  N\mathfrak m+1$, we deduce 
\begin{equation}\label{zebn3}
\mathbb P \left[ \sup_{1\leq j\leq J} \left\lvert\frac{S_N(t;h_j) }{\EE[S_N(t;h_j)]} - 1\right\rvert > \delta  \right] \leq 2 \exp\left[\log(N\mathfrak m) -\delta^2 C_{g,L} N\mathfrak m \min \mathcal H_N  /3\right],
\end{equation}
and the exponential bound tends to zero when $C_\delta$ is sufficiently large.  
Next, the supremum over the grid can be extended  over $\mathcal H_N$ using the Lipschitz continuity of the map $h\mapsto \EE[S_N(t;h)]$ and the monotonicity of the maps $h\mapsto S_N(t;h)$ and  $h\mapsto \EE[S_N(t;h)]$. Finally, by \eqref{zebn}, we write 
$$
\frac{N\mathfrak m h}{\mathcal N_\mu (t;h) }\geq \frac{N\mathfrak m h}{\EE[S_N(t;h)] }\times \frac{\EE[S_N(t;h)]}{S_N(t;h)}\times \frac{S_N(t;h)}{\mathcal N_\mu (t;h) } \geq  \frac{N\mathfrak m h}{\EE[S_N(t;h)] }\times \frac{\EE[S_N(t;h)]}{S_N(t;h)},
$$
and we deduce the first inequality in \eqref{th_q1a_SM} from \eqref{zebn2} and \eqref{zebn3}. 

For the last inequality  in \eqref{th_q1a_SM},  with a uniform kernel we have by definition
\begin{equation}\label{simpl_cas}
	\frac{\WN(t;h)}{\mathcal N_\mu (t;h)} = \frac{1}{\WN(t;h)} \sum_{i=1}^N \frac{w_i(t;h)}{\mathcal N_i(t;h)}  \leq 1. 
\end{equation} 
From this, \eqref{mfkb} and the first inequality in \eqref{th_q1b_SM}, that we will justify below, we deduce that 
\begin{equation} 
	\mathcal N_\mu (t;h)^{-1 } \leq \frac{\mathfrak c_1 ^{-1} c_L^{-1}\{1+o_{\mathbb P}(1)\}}{N  \mathfrak m h} ,
\end{equation}
with some $o_{\mathbb P}(1)$ term which does not depend on $h$, provided $\min\{1,\mathfrak m h\} =\mathfrak m h$.
For the case $\mathfrak m h >1$, let us note that, given $M_i$, $\mathcal N_i(t;h)$ is a Binomial random variable with parameters $M_i$ and $p_i(t;h)$. It can be show, see also \eqref{df_chaw}, that if $U$ follows a binomial distribution $B(n,p)$, then 
$$
\EE \left[\frac{\mathbf 1\{U>0\}}{U}\right] \leq \frac{2}{(n+1)p} - \frac{2n}{n+1}q^n \leq\frac{2}{(n+1)p},\quad q=1-p.
$$
Applying this for each $\mathcal N_i$, we derive a bound for the expectation of $\sum_{i=1}^N   w_i   \mathcal N_{i}^{-1}$. By the arguments used to derive \eqref{zebn3}, the $\sum_{i=1}^N   w_i   \mathcal N_{i}^{-1}$ concentrates around its expectation. Bounding this expectation and using 
 again \eqref{mfkb} and the first inequality in \eqref{th_q1b_SM}, 
we deduce 
\begin{equation}\label{new_star*}
	\mathcal N_\mu (t;h)^{-1 } \leq  c_L ^{-1} C_{g,L}^{-1}   \{(\mathfrak m+1) h\}^{-1}   \times \WN(t;h)^{-1} \{1+o_{\mathbb P}(1)\}
	\leq  
	(\mathfrak c_1  c_L  C_{g,L} )^{-1} \frac1{N   \mathfrak m h} \{1+o_{\mathbb P}(1)\},
\end{equation}
with the $o_{\mathbb P}(1)$ rate uniform with respect to $h\in\mathcal H_N$. The justification of the last inequality in \eqref{th_q1a_SM} is now complete.

Next, to show \eqref{th_q1b_SM}, note that, given $M_i$,  the indicator $w_i$ is a Bernoulli variable with parameter
\begin{equation}\label{up_bd1a}
\pi_i(t;h) = 1 -  \{1- p_i(t;h)\}^{M_i}.
\end{equation}
Let us notice that,  for any $M>0$ and for any $u\in(0,1)$,
$$
- M \frac{u}{1-u} \leq \log (1-u)^{M} < -uM.
$$
Assuming, without loss of generality, that $p_i(t;h)\leq 1/2$, for all $ h\in\mathcal H_N$ and for all $i$, we  deduce
\begin{equation}\label{bounds_pii}
1- \exp(-M_i  p_i(t;h))  \leq \pi_i(t;h) \leq 1- \exp(-2M_i  p_i(t;h)).
\end{equation}
By \eqref{bgfm1}, we have, for all $h\in\mathcal H_N$ and $1\leq i\leq N,$
$$
2 C_{g,L}h \leq p_i(t;h) \leq 2C_{g,U}h.
$$
From this and \eqref{mfkb}, we have, for all $h\in\mathcal H_N$ and $1\leq i\leq N,$
\begin{equation}\label{bounds_pii2}
1- \exp(-2 C_{g,L} c_L \mathfrak m h   )  \leq \pi_i(t;h) \leq 1- \exp(-4 C_{g,U} C_U \mathfrak m h),
\end{equation}
from which we deduce, for all $h\in\mathcal H_N$,
\begin{align}
1- \exp(-2 C_{g,L} c_L \mathfrak m \min \mathcal H_N    )  &\leq
1- \exp(-2 C_{g,L} c_L \mathfrak m h   ) \nonumber\\
&\leq\frac{ \EE_M [\mathcal W_N(t;h)] }{N} = \frac{1}{N} \sum_{i=1}^N \pi_i(t;h) \nonumber\\
&\leq 1- \exp(-4 C_{g,U} C_U \mathfrak m h) \nonumber\\
&\leq 1- \exp(-4 C_{g,U} C_U \mathfrak m \max \mathcal H_N ). \label{bounds_pii3}
\end{align}
Condition \eqref{eq:cdg_H} imposes $N\mathfrak m \min \mathcal H_N \rightarrow \infty$. Let us now consider the case 
$\mathfrak m \min \mathcal H_N \rightarrow 0$, the arguments for the case $\liminf\{ \mathfrak m \min \mathcal H_N\} > 0$ being straightforward. Since $1- \exp(-x) = x\{1+o(1)\}$ when $x$ decreases to zero, we deduce \eqref{th_q1b_SM} with 
$\EE[\mathcal W_N(t;h)]$ instead of $\mathcal W_N(t;h)$. Next, similarly to the justification of \eqref{th_q1a_SM}, we use the Chernoff's exponential bound and a grid on $\mathcal H_N$ to replace $\EE[\mathcal W_N(t;h)]$ by $\mathcal W_N(t;h)$. The property \eqref{th_q1b_SM} follows and we thus complete the proof of Theorem~\ref{thm:optim_mu}.
\end{proof}

%%%%%%%%%%%%%%%%%%%%%%%%%%%%%
%%%%%%%%%%%%%%%%%%%%%%%%%%%%%

% !TeX root=bj-sample-SM.tex

\begin{proof}[Complements for the proof of Theorem \ref{thm:optim_gamma}]
We here provide the arguments for \eqref{th_q1ab_cov} in the main text.  For simplicity, we assume $\widehat \Gamma_N^* $ is built with the uniform kernel. Recall that $s\neq t$ are fixed and without loss of generality, we  consider $\sup \mathcal H_N< \lvert s-t \rvert/2$. We can also assume  $c_L\mathfrak m \geq 2$.

For the justification of the first part of \eqref{th_q1ab_cov}, we start by proving the following result.

\begin{lemSM}\label{WN_order_SMb} 
	Assume that \eqref{mfkb} and \eqref{bgfm1} hold true. Then constants $\underline c_W,\overline c_W\in(0,1]$  exist  such that 
	\begin{equation}\label{sdyu_bis}
		\underline c_W N \leq \frac{ \mathbb E_{M} \left[\mathcal W_N(s,t;h) \right]}{\min\{1,(\mathfrak m h)^2\}} \leq \overline c_W N,\quad  \forall h\in\mathcal H_N, \lvert s-t \rvert> 2h.
	\end{equation}
\end{lemSM}

\begin{proof}[Proof of Lemma SM.\ref{WN_order_SMb}]Without loss of generality, we can assume 
	$\sup_{h\in \mathcal H_N}  p_i(t;h) < 1/8$, for all $1\leq i \leq N$, where $p_i(t;h)$ are defined in  \eqref{def_pig}.
	Since $w_i(s;h)$ and $ w_i(t;h) $ are indicator functions,  we have
$$
w_i(s;h) w_i(t;h) \geq 1 - \{1 -w_i(s;h)\} - \{1 -w_i(t;h)\},
$$
and thus 
\begin{equation}
	\pi_i(s,t;h) \geq  1 - \{1-p_i(s;h)\}^{M_i} - \{1-p_i(t;h)\}^{M_i}, 
\end{equation}
where $	\pi_i(s,t;h) \coloneq\EE[w_i(s;h) w_i(t;h) ]$. We first consider the case $\mathfrak mh\gtrsim 1$. We can check that  $ \{1-p_i(s;h)\}^{M_i} + \{1-p_i(t;h)\}^{M_i}\leq 1/2$, provided $\mathfrak mh>C$ for some large constant $C$. We can then focus on investigating the case $\mathfrak mh \asymp 1$. By definition of the multinomial distribution, we have 
\begin{equation}
	\pi_i(s,t;h) =  \sum_{l+l^\prime=0 }^{M_i-2 }\frac{M_i! \, \times p_i(s;h)^{l+1} p_i(t;h)^{l^\prime+1} \{1-p_i(s;h) - p_i(t;h)\}^{M_i-2- (l+l^\prime)}}{(l+1)!\,(l^\prime+1)! \,(M_i - 2  - (l+l^\prime))!}.
\end{equation}
Let 
$$
\underline \pi_i(s,t;h) \coloneq \frac{M_i!}{  (M_i - 2 )!}  p_i(s;h)  p_i(t;h)  \{1-p_i(s;h) - p_i(t;h)\}^{M_i-2 } ,
$$
such that, for all $s, t, h$,
\begin{equation}
	\underline \pi_i(s,t;h) \leq 	\pi_i(s,t;h). 
\end{equation}
By construction and a continuity argument, we can check that  for any $ c_L, c_U>0$, a constant $\underline c$ exists (depending on  $ c_L, c_U$) such that 
$$
\inf_{s,t}\underline \pi_i(s,t;h) \geq \underline c >0, \quad \forall h \in\mathcal H_N \; \text{ such that } \mathfrak m h \in[ c_L, c_U].
$$
Gathering all these facts, it is clear that when $\mathfrak m h \gtrsim  1$, $\pi_i(s,t;h)  $ is uniformly bounded from below and above by positive constants, the arguments are complete for the case where $\min\{1,\mathfrak m h \} \asymp  1.$

For the case $\mathfrak m h \lesssim 1$,  whenever $\lvert s-t \rvert >2h$, we  will bound $\pi_i(s,t;h)$ from below and from above using  $\underline \pi_i(s,t;h)$. Next, modulo suitable constants, we bound $\underline \pi_i(s,t;h)$ from below and from above by $(\mathfrak m h)^2$ multiplied by suitable constants. Let
$$
4C_{g,L} h \leq \underline g(s,t;h) = p_i(s;h) + p_i(t;h)\leq 4C_{g,U} h,
$$
and note that 
\begin{align}
	\{1-p_i(s;h) - p_i(t;h)\}^{M_i-2 }  &= \exp\left((M_i-2)  \log\{1-\underline g(s,t;h) \} \right) \\
	&\geq  \exp\left( - (M_i-2)  \frac{\underline g(s,t;h)}{1-\underline g(s,t;h)} \right) \\
	&\geq \exp\left( - (C_{U} \mathfrak m -2)  \frac{4C_{g,U} h}{1-4C_{g,L} h} \right) \\
	&\geq 
	\exp\left( - C_1 \mathfrak m h \right),
\end{align}
for some constant $C_1>0$. Uniformly bounding $p_i(s;h)$ and $p_i(t;h)$ from below, we deduce that a constant $C_2>0$ exits such that, for all $s, t, h$,
\begin{equation}
	C_2 \exp( - C_1) \times (\mathfrak m h)^2  \leq C_2 (\mathfrak m h)^2 \exp\left( - C_1 \mathfrak m h \right) \leq 	\underline \pi_i(s,t;h) \leq 	\pi_i(s,t;h). 
\end{equation}
We now show that when $\mathfrak m h \lesssim  1$, modulo a constant, $\underline \pi_i(s,t;h) $ is also an upper bound for $ \pi_i(s,t;h) $.
We rewrite
\begin{align}
	\pi_i(s,t;h) &= \underline \pi_i(s,t;h)  
	\times\sum_{l+l^\prime=0 }^{M_i-2 }C(l,l^\prime ) \times \frac{(M_i -2) !}{l! \, l^\prime! \, (M_i - 2  - (l+l^\prime))!} \times q_1 ^l \times q_2^{l^\prime}\times q_3^{M_i - 2  - (l+l^\prime)}\\
	&\leq \underline \pi_i(s,t;h)  \times  \sup_{s,t,h}\max_{l,l^\prime }C(l,l^\prime ) ,
\end{align}
where
$$
q_1 = \frac{p_i(s;h) }{1- \{p_i(s;h)+p_i(t;h)\}},  \qquad  q_2=\frac{ p_i(t;h)}{1- \{p_i(s;h)+p_i(t;h)\}} \qquad   	q_3 = \frac{1- 2\{p_i(s;h)+p_i(t;h)\} }{1- \{p_i(s;h)+p_i(t;h)\}} ,
$$
with $q_1,q_2,q_3 \in (0,1)$, $q_1 + q_2 + q_3 = 1$, and
$$
C(l,l^\prime ) =\frac{1}{(l+1)(l^\prime+1)}
\times \frac{\{1-\{p_i(s;h) + p_i(t;h)\}\}^{M_i-2- (l+l^\prime)} }{\{1-2\{p_i(s;h) + p_i(t;h)\}\}^{M_i-2- (l+l^\prime)} }\leq   \frac{\{1-\{p_i(s;h) + p_i(t;h)\}\}^{M_i} }{\{1-2\{p_i(s;h) + p_i(t;h)\}\}^{M_i} } .
$$
It remains to  bound the RHS in the last display and show that $C(l,l^\prime )$ is bounded by a constant, provided $\mathfrak m h \lesssim 1$. Let 
$$
\overline g(s,t;h) = \frac{p_i(s;h) + p_i(t;h)}{1-2\{p_i(s;h) + p_i(t;h)\} } \leq 2\{p_i(s;h) + p_i(t;h)\} ,
$$
such that, for all $l,l^\prime \leq M_i-2$ and for all $s,t,h$,
\begin{align}
	C(l,l^\prime ) &\leq  \{1+\overline g(s,t;h) \}^{M_i} \\
	&= \exp\left(M_i  \log\{1+\overline g(s,t;h) \} \right) \\
	&\leq  \exp\left(C_{U}\mathfrak m  \log\{1+\overline g(s,t;h) \} \right) \\
	&\leq \exp\left(C_{U}\mathfrak m  \overline g(s,t;h) \right) \\ 
	&\leq \exp\left(8 C_{g,U}C_{U}\mathfrak m  h \right).
\end{align}
Thus, $C(l,l^\prime )$ is uniformly  bounded and a  constant $c_5>0$ exists such that, for all $s, t, h$,
\begin{equation}
	\underline \pi_i(s,t;h) \leq 	\pi_i(s,t;h) \leq c_5 	\underline \pi_i(s,t;h) \leq c_5 C_{U}^2	(\mathfrak m h)^2. 
\end{equation}
The proof of Lemma SM.\ref{WN_order_SMb} is now complete. 
\end{proof}

The next step is to show that $\mathcal W_N(s,t;h)$ concentrates around its conditional mean given the  integers $M_1,\ldots,M_N$, uniformly with respect to $h\in\mathcal H_N$. We then deduce that $\mathcal W_N(s,t;h)/\EE[\mathcal W_N(s,t;h)]=1+o_\PP(1)$ uniformly with respect to $h\in\mathcal H_N$. Let us sketch the arguments for getting this uniform rate. Recalling that $\mathcal W_N(s,t;h)$ is a sum of independent Bernoulli variables, by Bernstein's inequality applied with a given $h$, we get an exponential bound for the conditional probability of the event 
$$
\left\{\left\lvert\frac{\mathcal W_N(s,t;h)}{\EE_M [\mathcal W_N(s,t;h)]} - 1 \right\rvert >   \frac{ y }{\sqrt{\EE_M[\mathcal W_N(s,t;h)]}}\right\},  \quad y>0, 
$$
given $M_1,\ldots,M_N$. We next note that, for any $s,t$ and any $h< h^\prime$, by definition 
$$
w_i(s;h) w_i(t;h) \leq w_i(s;h^\prime ) w_i(t;h^\prime ).
$$
This property allows to build a set of brackets, for the family of indicators $w_i(s;h) w_i(t;h)$ indexed by $h$, by taking a grid of  $(N\mathfrak m) ^2$ values in $\mathcal H_N$. By applying Bernstein's inequality for each $h$ in the grid, and using the Boole's inequality, we deduce 
\begin{equation}\label{Boole11}
\sup_{h\in\mathcal H_N}\left\lvert \frac{\mathcal W_N(s,t;h)}{\EE_M [\mathcal W_N(s,t;h)]} - 1 \right\rvert	=  o_\PP(1).
\end{equation}
This, combined with \eqref{sdyu_bis} and condition \eqref{mfkb}, yields the first part of \eqref{th_q1ab_cov}.

For the second part of \eqref{th_q1ab_cov}, we start by showing  that there exists a constant $\mathfrak c_1>0$ such that 
\begin{equation}\label{th_qq1a}
\inf_{h\in\mathcal H_N}\frac{N \min \{\mathfrak mh, (\mathfrak mh)^2\} }{\mathcal N_\Gamma (t \mid s;h) }\geq  \mathfrak c_1 \{1+o_{\mathbb P}(1)\}.
\end{equation}
Using the fact that the harmonic mean is less than or equal to the mean, we obtain
\begin{equation}\label{th_qq1a1}
\frac{1}{\mathcal N_\Gamma (t \mid s;h) }\geq \frac{c_i (t;h)}{\sum_{i=1}^N  w_i (s;h) w_i(t;h)  \mathcal N_{i} (t;h)  },
\end{equation}
with $w_i(t;h)$, $c_i (t;h)$ and  $ \mathcal N_{i} (t;h)$  defined in \eqref{def_wi}, \eqref{eq:ci} and \eqref{eq:Ncal1}, respectively. For the case we consider, for all $i$, we have  $c_i\equiv w_i$.  To justify \eqref{th_qq1a}, it suffices to prove that a positive constant $c_{\mathcal N} $ exists such that 
\begin{equation}\label{th_q1a3}
\frac{\sum_{i=1}^N  w_i (s;h) w_i(t;h)\mathcal N_{i} (t;h) }{N \min \{\mathfrak mh, (\mathfrak mh)^2\}  }\leq c_{\mathcal N}\{1+ o_{\mathbb P}(1) \} ,
\end{equation}
with the $o_{\mathbb P}(1)$ uniform with respect to $h\in\mathcal H_N$. 
Let us notice that, in the case of a NW estimator with a uniform kernel,  
\begin{equation}
	\sum_{i=1}^N w_i (s;h) w_i(t;h)\mathcal N_{i} (t;h) = \sum_{i=1}^Nw_i (s;h) \sum_{m=1}^{M_i} \mathbf{1}\left\{\lvert \Tnm-t \rvert \leq h\right\}\\ 
	 = \sum_{i=1}^N S^{(i)},
\end{equation} 
with 
$
S^{(i)} = S^{(i)}(h)  = w_i (s;h) \sum_{1\leq m \leq M_i }  \mathbf{1}\{\lvert\Tnm-t \rvert\leq h\}.
$
We thus need to suitably bound the sum of $ S^{(i)}(h)$ from above. 
Let $\PP_M, \EE_M$, and $\operatorname{Var}_M$,
denote the conditional probability,  expectation and variance, respectively, given $M_1,\ldots,M_N$.
We then have
\begin{align*}
\EE_M[S^{(i)}] &=   \sum_{m=1}^{M_i} 
\EE_M \left[ w_i (s;h)  \mathbf{1}\left\{\lvert \Tnm-t \rvert\leq h\right\}  \right] \\ 
&= \sum_{m=1}^{M_i} \EE_M \left[ \mathbf{1}\left\{\sum_{1\leq m^\prime \neq m \leq M_i } \mathbf{1}\left\{\lvert T^{(i)}_{m^\prime} -s \rvert \leq h\right\} \geq 1\right\} \mathbf{1}\left\{\lvert \Tnm-t \rvert\leq h\right\}   \right] \\
&= \sum_{m=1}^{M_i} \EE_M \left[ \mathbf{1}\left\{\sum_{1\leq m^\prime \neq m \leq M_i }  \mathbf{1}\left\{\lvert T^{(i)}_{m^\prime} -s \rvert\leq h \right\} \geq 1\right\}  \right]  \times \EE_M \left[ \mathbf{1}\left\{ \lvert \Tnm-t \rvert \leq h \right\} \right] \\
&= \left[  1 -  \{1- p_i(t;h)\}^{M_i-1}\right] \times  M_i p_i(t;h) \\
&=\{1+o(1)\} \times \pi_i(s;h)  \times  M_i p_i(t;h),
\end{align*}
where  $p_i(t;h) =\int_{t-h}^{t+h}g_i(u)du $ and $\pi_i(s;h)$ is defined as in \eqref{up_bd1a}. The $o(1)$ term is uniform with respect to $h$. 
Moreover, 
\begin{align*}
\{S^{(i)}\}^2 - S^{(i)} &=  w_i (s;h) \hspace{-1em}\sum_{1\leq m^\prime \neq m \leq M_i } 
 \hspace{-1em}\mathbf{1}\left\{\lvert \Tnm-t \rvert\leq h\right\}\mathbf{1}\left\{\lvert T^{(i)}_{m^\prime} -t \rvert \leq h\right\}  \\
&=  \mathbf{1}\left\{\sum_{\substack{1\leq m^{\prime\prime}  \leq M_i \\ m^{\prime\prime} \not\in\{ m, m ^\prime\} }} \hspace{-1em}\mathbf{1}\left\{\lvert T^{(i)}_{m^{\prime\prime}} -s \rvert \leq h\right\} \geq 1\right\}  \sum_{1\leq m^\prime \neq m \leq M_i } \hspace{-1em} \mathbf{1}\left\{\lvert \Tnm-t \rvert \leq h\right\}\mathbf{1}\left\{\lvert T^{(i)}_{m^\prime} -t \rvert \leq h\right\},
\end{align*}
and thus, 
\begin{align}
\EE_M[\{S^{(i)}\}^2] &= \EE_M[S^{(i)}  ]+ \left[  1 -  \{1- p_i(t;h)\}^{M_i-2}\right]\times M_i (M_i-1)  p_i^2(t;h) \nonumber\\
&= \{1+o(1)\}\times  \pi_i(s;h)  \times  M_i p_i(t;h)\times\{1+M_i p_i(t;h)\}, \label{thg1}
\end{align}
with the $o(1)$ term uniform with respect to $h$. We deduce that
\begin{align}
\operatorname{Var}_M[S^{(i)} ] &= \EE_M[\{S^{(i)}\}^2 ] -  \EE^2_M[S^{(i)}  ] 
\\ &= \{1+o(1)\}\times  \pi_i(s;h)    M_i p_i(t;h) \times\left[  1+M_i p_i(t;h) - \pi_i(s;h)  M_i p_i(t;h)\right] \\
&=  \{1+o(1)\}\times  \pi_i(s;h)    M_i p_i(t;h) +   \{1+o(1)\}\times  \pi_i(s;h)   \{1- \pi_i(s;h)  \} \{ M_i p_i(t;h) \}^2.
\end{align}
Let us recall  the following notation:
given $\varphi_1 $, $\varphi_2 $, positive functions of $M_i$ and $h$, 
$$
\varphi_1  \lesssim \varphi_2  \quad\Leftrightarrow  \quad \exists C>0   \quad \text{ a constant such that } \quad \varphi_1  \leq C \varphi_2  ,
$$
and 
$$
\varphi_1  \asymp \varphi_2  \quad\Leftrightarrow  \quad \varphi_1  \lesssim \varphi_2  \quad\text{and} \quad\varphi_2  \lesssim \varphi_1.
$$
With this notation, 
$\EE_M[S^{(i)}  ] \asymp \pi_i(s;h)  \times  \mathfrak m h$,
and
\begin{equation}\label{ger1c}
\EE_M\left[\sum_{i=1}^N S^{(i)} \right] \asymp  \EE_M\left[\mathcal W_N(s;h) \right]  \times  \mathfrak m h,
\end{equation}
and thus, by \eqref{bounds_pii3}, for all $h\in\mathcal H_N$,
\begin{equation}\label{bounds_pii3_b}
 N\mathfrak m h \left\{ 1- \exp(-2 C_{g,L} c_L \mathfrak m h    )  \right\} 
\lesssim \EE_M\left[\sum_{i=1}^N S^{(i)} \right]
\lesssim N\mathfrak m h \left\{ 1- \exp(-4 C_{g,U} C_U \mathfrak m h) \right\}.
\end{equation}
On the other hand, 
$\operatorname{Var}_M[S^{(i)} ] \asymp \pi_i(s;h)  \times  \mathfrak m h + \pi_i(s;h) \{1- \pi_i(s;h)  \} \times  (\mathfrak m h)^2$.
By \eqref{bounds_pii2}, we deduce 
$$
 \left\{1- \exp(-2 C_{g,L} c_L \mathfrak m h   )  \right\}  \mathfrak m h \{ 1 +\exp(-4 C_{g,U} C_U \mathfrak m h) \mathfrak m h  \} 
\lesssim \operatorname{Var}_M[S^{(i)} ]
$$
and
$$
\operatorname{Var}_M[S^{(i)} ] \lesssim \left\{1- \exp(-4 C_{g,U} C_U \mathfrak m h   )  \right\}  \mathfrak m h \{ 1 +\exp(-2 C_{g,L} c_L \mathfrak m h) \mathfrak m h  \} .
$$
Since for any $c>0$, the map $x\mapsto x\exp(-cx)$, $x\geq 0$ is bounded, we deduce 
\begin{equation}\label{ger1bb}
 \left\{1- \exp(-2 C_{g,L} c_L \mathfrak m h   )  \right\}  \mathfrak m h \lesssim \operatorname{Var}_M[S^{(i)}  ] \lesssim \left\{1- \exp(-4 C_{g,U} C_U \mathfrak m h   )  \right\}  \mathfrak m h.
\end{equation}
Let us note that
\begin{equation}\label{reeq1n}
\EE_M[S^{(i)}] \asymp\operatorname{Var}_M[S^{(i)}  ] \asymp \mathfrak mh \times \min \{1, \mathfrak mh\} .
\end{equation}
It remains to show that the sum of $S^{(i)}(h)$ concentrates around a quantity which allows us to deduce \eqref{th_q1a3}. 
Let  $A = A(h) >0$ to be determined below, and let
$$
\mathcal A =\mathcal A(h) =  \left\{\max_{1\leq i \leq N} S^{(i)} \leq A\right\} \quad \text{and} \quad S^{(i)}_A = S^{(i)}_A (h) = S^{(i)}\mathbf{1}_{\mathcal{A}}.
$$
Let 
$$
E_{M,A} = E_{M,A}(h) \coloneq \EE_M\left[\sum_{ i=1} ^NS_A^{(i)} \right] \leq  \EE_M \left[\sum_{ i=1} ^NS^{(i)} \right] \quad \text{and} \quad V_{M,A} = V_{M,A}(h) \coloneq\operatorname{Var}_M \left[\sum_{ i=1} ^NS_A^{(i)} \right].
$$
By definition, 
$$
V_{M,A}   \lesssim \sum_{1=1}^N \operatorname{Var}_M[S^{(i)}  ] \asymp  N \mathfrak m h \times \min \{1, \mathfrak mh\}  \eqqcolon  \Omega_N(h) \longrightarrow \infty.
$$
Indeed, we have
\begin{align}
\operatorname{Var}_M[S_A^{(i)}] &= \EE_M [ \{S^{(i)} \}^2 ] -  \EE_M [ \{ S^{(i)}  \}^2 \mathbf{1}_{\overline{\mathcal{A}}}]
-\left\{\EE_M [ S^{(i)}  ]  - \EE_M [ S^{(i)}  \mathbf{1}_{\overline{\mathcal{A}}}] \right\}^2 \\
&\leq \operatorname{Var}_M[S^{(i)}  ] + 2 \EE_M [ S^{(i)}  ]   \EE_M [ S^{(i)}  \mathbf{1}_{\overline{\mathcal{A}}}] .
\end{align}
Herein, for any set $B$, $\overline B$ denotes its complement. By \eqref{reeq1n}, we deduce 
$$
\operatorname{Var}_M[S_A^{(i)}  ]  \lesssim  \operatorname{Var}_M[S^{(i)}  ] ,
$$
provided a constant exists such that $\EE_M [ S^{(i)}  \mathbf{1}_{\overline{\mathcal{A}}}]\leq C $ for all $\mathfrak m$ and $h$. By the Cauchy-Schwarz inequality and \eqref{thg1},
$$
 \EE_M [ S^{(i)} (h) \mathbf{1}_{\overline{\mathcal{A}}(h)}] \leq  \EE_M^{1/2} [ \{ S^{(i)}(h)  \}^2]\times \PP (\overline{\mathcal{A}}(h)) \lesssim \mathfrak m h \times \PP (\overline{\mathcal{A}}) \leq  \mathfrak m  \times \PP (\overline{\mathcal{A}}(\min \mathcal H_N)) \longrightarrow 0.
$$
The convergence to zero follows from \eqref{rezt6} below. Next, by the Bernstein inequality applied to the $ S_A^{(i)}$, for each $h\in\mathcal H_N$,
\begin{equation}
\PP_M \left[ \sum_{ i=1} ^N S_A^{(i)}(h) >  E_{M,A} (h) + \Omega_N(h)  \right]  
\leq \exp\left( - \frac{\Omega_N(h) ^2/2}{V_{M,A}(h)  + A(h) \Omega_N(h) /3} \right).
\end{equation}
To derive bounds for the concentration probability of  the sum of $ S^{(i)}(h)  $,  it suffices to take 
$A$ such that
$$
\sqrt{V_{M,A}(h) } \ll\Omega_N(h) \quad  \text{and} \quad A(h)\ll \Omega_N(h).
$$
Let 
$$
A(h)=  \frac{\Omega_N(h) }{c_A \log(\mathfrak m)},
$$
with $c_A$ some large constant. 
Consider $\mathcal G_N$ a uniform  a grid in $\mathcal H_N$ with mesh of rate $1/N\mathfrak m$. By \eqref{mfkb} and the taking $c_A$ sufficiently large, we deduce
that a constant $0< C < c_A  $ exists such that 
\begin{align}
\PP_M \left[ \sup_{h\in\mathcal H_N}\frac{1} {\Omega_N(h) } \sum_{ i=1} ^N S_A^{(i)}(h)>  C \right]  &\leq \PP_M \left[ \sup_{h\in\mathcal G_N}\frac{1} {\Omega_N(h) } \sum_{ i=1} ^N S_A^{(i)}(h)>  C/2 \right]  \\
 &\leq \exp\left( \log (\lvert\mathcal G_N\rvert) - c_A \log (\mathfrak m) \right) \\
 &\leq \exp\left(   -  (c_A-C) \log (\mathfrak m) \right) \longrightarrow 0.
\end{align}
Here, 
 $\lvert \mathcal G_N \rvert$ denotes the cardinal of $\mathcal G_N$. Finally, we have
\begin{align}
\PP_M \left[ \sup_{h\in\mathcal H_N}\frac{1} {\Omega_N(h) } \sum_{ i=1} ^N S^{(i)}(h)>  C \right] &= \PP_M \left[ \sup_{h\in\mathcal H_N}\frac{1} {\Omega_N(h) } \sum_{ i=1} ^N S^{(i)}(h) \{\mathbf 1 _{\mathcal A(h)}+\mathbf 1 _{\overline{\mathcal A}(h)} \} >  C \right] \\
 &\leq \PP_M \left[ \sup_{h\in\mathcal G_N}\frac{1} {\Omega_N(h) } \sum_{ i=1} ^N S_A^{(i)}(h)>  C /4 \right]  \\
  &\qquad + \PP_M \left[ \sup_{h\in\mathcal H_N} \mathbf{1}_{\overline{\mathcal{A}} (h)} >  0\right]  \\
 &\leq \exp\left( - (c_A - C) \log (\mathfrak m) \right) +  \PP_M \left[ \sup_{h\in\mathcal H_N} \mathbf{1}_{\overline{\mathcal{A}} (h)} >  0\right] .
\end{align}
Next, let $h_{j(h)} $, with $1\leq j(h) \le J$, be the point in the grid $\mathcal G_N$ such that $h_{j(h)-1}  \leq h < h_{j(h)} $. Using the monotonicity of the $ S^{(i)} (h)$ and $ \Omega_N(h)$, with respect to $h$, we have
\begin{equation}
\left\{\max_{1\leq i \leq N}
\frac{ S^{(i)} \left(h_{j(h)-1} \right) }{\Omega_N \left(h_{j(h)}\right) }
\geq \frac{1}{c_A\log(\mathfrak m )}\right\}
\subset \overline{\mathcal{A}} (h)  \subset
\left\{\max_{1\leq i \leq N} 
 \frac{ S^{(i)} \left(h_{j(h)} \right) }{\Omega_N \left(h_{j(h)-1} \right) }
\geq \frac{1}{c_A\log(\mathfrak m )}\right\}.
\end{equation}
This implies 
\begin{align}
\PP_M \left[ \sup_{h\in\mathcal H_N} \mathbf{1}_{\overline{\mathcal{A}} (h)} >  0\right] &\leq 
  \sum_{j=2}^J \PP_M \left[ \max_{1\leq i \leq N} 
  \frac{ S^{(i)} (h_j) }{\Omega_N \left( h_{j-1} \right) }
  \geq \frac{1}{c_A\log(\mathfrak m )}\right] \nonumber\\
 & \leq\sum_{j=2}^J \sum_{i=1}^N  \PP_M \left[ S^{(i)} (h_j) 
  \geq \frac{\Omega_N(h_{j-1} )}{c_A\log(\mathfrak m )}\right] \nonumber\\
 & \leq \sum_{j=2}^J \sum_{i=1}^N 
  \PP_M \left[ 
  \sum_{m=1}^{M_i} \mathbf{1}\left\{\lvert \Tnm-t \rvert \leq h_j\right\}
  \geq M_i \EE \left[\mathbf{1}\left\{\lvert \Tnm-t \rvert \leq h_j\right\} \right] \times (1+ \delta_{ij}) \right] \nonumber\\
 & \leq J\times N \times \exp\left(-  c_LC_{g,L}\min_{1\leq i \leq N}\min_{2\leq j \leq J}\left[\frac{\delta_{ij}}{2+\delta_{ij}} \times\delta_{ij} \mathfrak m h_j\right]\right), \label{rezt6}
\end{align}
where for the last inequality, we use Chernoff's inequality, and $c_L$ and $C_{g,L}$ are the constants in \eqref{mfkb} and \eqref{bgfm1}, respectively. Here, 
\begin{equation}
\delta_{ij} = \frac{\Omega_N(h_{j-1} )/ \{c_A\log(\mathfrak m )\}}{M_i \EE \left[\mathbf{1}\left\{\lvert\Tnm-t\rvert\leq h_j\right\} \right]}
\geq C  \frac{ N   \min \{1, \mathfrak m h_j \} }{c_A\log(\mathfrak m )}
\geq C \frac{ N  \min \{1, \mathfrak m \min \mathcal H_N \} }{c_A\log(\mathfrak m )} \longrightarrow \infty,
\end{equation}
for some constant $C>0$. Moreover, by the condition $N\{\mathfrak m \min \mathcal H_N \}^2/\log^2 (N\mathfrak m) \rightarrow \infty$, we have 
$$
\delta_{ij} \mathfrak m h_j \geq C \frac{ N  \mathfrak m \min \mathcal H_N \min \{1, \mathfrak m \min \mathcal H_N \} }{c_A\log(\mathfrak m )} \gg \log(JN).
$$
This implies that the exponential bound in \eqref{rezt6} tends to zero. 
Gathering facts, we deduce \eqref{th_q1a3} and then \eqref{th_qq1a}.

%%%%%%%%%%%%%%%%%%%%%%%%%%%
%%%%%%%%%%%%%%%%%%%%%%%%%%%
%%%%%%%%%%%%%%%%%%%%%%%%%%%

Next, we show   that it exists a constant $\mathfrak C_1$ such that 
\begin{equation}\label{th_qq1a8}
	\sup_{h\in\mathcal H_N}\frac{N \min \{\mathfrak mh, (\mathfrak mh)^2\} }{\mathcal N_\Gamma (t \mid s; h) }\leq  \mathfrak C_1 \{1+o_{\mathbb P}(1)\}.
\end{equation}
For this purpose, we show that a constant $\overline c_{\mathcal N}>0$  exists  such that, for all $h < \lvert s-t \rvert /2$,
\begin{equation}\label{sdyu_g-SM}
\frac{1}{   \mathcal N_\Gamma (s \mid t;h) } 
	\leq \frac{\overline c_{\mathcal N} \{1+o_{\mathbb P}(1)\}}{\max\{1,\mathfrak m h\}\times \EE_M [ \mathcal{W}_N(s,t;h)]},
\end{equation}
with the $o_{\mathbb P}(1)$ uniform with respect to $h\in\mathcal H_N$. To  get \eqref{th_qq1a8}, we note that by definition
\begin{equation}\label{eq_891}
	\frac{\mathcal{W}_N(s,t;h)} {\mathcal N_{\Gamma} (t \mid s;h)} = \frac{1}{\mathcal{W}_N(s,t;h)} \sum_{i=1}^N \frac{w_i(s;h)w_i(t;h)} {\mathcal N_{i} (t;h)} \leq \frac{1}{\mathcal{W}_N(s,t;h)} \sum_{i=1}^N w_i(s;h)w_i(t;h) \leq 1.
\end{equation}
In the case $\max\{1, \mathfrak m h\}\asymp 1$, it suffices to apply Lemmas SM.\ref{WN_order_SMb} and \eqref{Boole11}. We now investigate the case $\mathfrak m h \gtrsim1$. Let us recall a result of 
\cite{chao72}: if $X$ is a non-degenerate Binomial random variable $B(n,p)$, then 
\begin{equation}\label{df_chaw}
\EE[(1+U)^{-1}] = \frac{1 - q^{n+1}}{(n+1)p},\quad \text{where}\quad q=1-p.
\end{equation}
On the other hand, we can write
$$
\EE\left[\frac{1}{1+U}\right] = \mathbb P (U=0)+ \EE\left[\frac{\mathbf 1\{U\geq 1\}}{U(1+1/U)}\right] \geq  \mathbb P (U=0)+\frac{1}{2}  \EE\left[\frac{\mathbf 1\{U>0 \}}{U}\right] .
$$
We deduce 
\begin{equation}\label{inv_bin_mom}
	\EE\left[\frac{\mathbf 1\{U>0 \}}{U}\right] \leq \frac{2}{(n+1)p} - \frac{2(1+np)}{(n+1)p}q^n\leq \frac{2}{(n+1)p} - \frac{2n}{n+1}q^n.	
\end{equation}
Let us now compute the expectation of $w_i(s;h)w_i(t;h) / \mathcal N_{i} (t;h)$ given $M_i$.
Let us notice that, with an uniform kernel $K$, the distribution of $\mathcal N_{i} (t;h)$ given 
$$
n^{(i)}(s;h) \coloneq \sum_{1\leq m \leq M_i }  \mathbf{1}\left\{\lvert \Tnm-s \rvert \leq h\right\} = k, \quad k=1,2,\ldots ,M_i-1,
$$
is a Binomial distribution with parameters $B(M_i-k, p_i(t;h)/\{1-p_i(s;h)\})$, while $n^{(i)}(s;h)$ is a Binomial distribution with parameters $B(M_i, p_i(s;h))$. We can write 
\begin{align}
	\EE_M \left[ \frac{w_i(s;h)w_i(t;h)}{\mathcal N_{i} (t;h)} \right] &= \EE_M \left\{ \EE_M \left[ \frac{w_i(t;h)}{\mathcal N_{i} (t;h)} \Bigm|  n^{(i)}(s;h)\right] \right\} \\
	&= \sum_{k=1}^{M_i-1}  \EE_M \left[ \frac{w_i(t;h)}{\mathcal N_{i} (t;h)} \Bigm|  n^{(i)}(s;h)=k\right]  \times \mathbb P_M \left( n^{(i)}(s;h)=k\right).
\end{align}
Letting $\mathcal N_{i} (t;h)$ to be the Binomial distribution $U$ in \eqref{inv_bin_mom}, we have, for all $1\leq k \leq M_i-1$,
\begin{align}
	\EE_M \left[ \frac{w_i(t;h)}{\mathcal N_{i} (t;h)} \Bigm|  n^{(i)}(s;h)=k\right]  
	&\leq 2 \frac{1-p_i(s;h)}{(M_i-k+1)p_i(t;h)}
	-2 \frac{M_i-k}{M_i-k+1} \left\{1 - \frac{p_i(t;h)}{1-p_i(s;h) }\right\}^{M_i-k} \\
	&\leq \frac{2}{(M_i-k+1)p_i(t;h)}.
\end{align}
Let us note that
$$
\sum_{k=0}^{M_i}  \frac{1}{M_i-k+1 }\times \mathbb P_M \left( n^{(i)}(s;h)=k\right) = \EE\left[\frac{1}{1+Y}\right],
$$
where $Y = M_i - n^{(i)}(s;h)$  is a Binomial distribution with parameters $B(M_i, 1-p(s;h))$.
Gathering facts,   we deduce that it exists a constant $\overline {\mathcal C}_{\mathcal N}$ such that for all  $h\in\mathcal H_N$ with $\mathfrak m h \gtrsim 1$,
\begin{equation}\label{harmo_m}
	\EE_M \left[ \frac{w_i(s;h)w_i(t;h)}{\mathcal N_{i} (t;h)} \right] \leq \frac{2}{p_i(t;h)} \times \frac{1 -\{1-p_i(s;h)\}^{M_i+1} }{(M_i+1) \{1-p_i(s;h)\}} \leq \frac{\overline {\mathcal C}_{\mathcal N}}{\mathfrak m h}.
\end{equation}
Next, since $\mathcal N_{i} (t;h)$ is a non-negative integer, the conditional variance of $w_i(s;h)w_i(t;h) / \mathcal N_{i} (t;h)$ given $M_i$ is smaller that its conditional expectation given $M_i$. Finally, by Bernstein's inequality applied for a grid of values and the Boole's inequality, the sum of the $w_i(s;h)w_i(t;h) / \mathcal N_{i} (t;h)$ deviates from its mean with exponentially small probability given the $M_i$, uniformly with respect to  $h$. 

%%%%%%%%%%%%%%%%%%%%%%%%%%%
%%%%%%%%%%%%%%%%%%%%%%%%%%%
%%%%%%%%%%%%%%%%%%%%%%%%%%%

From \eqref{th_q1ab_cov}, we deduce that
\begin{equation}\label{th_qq1_SM}
	\frac{1}{\WN(s,t;h)} - \frac{1}{N} \leq  \min\left[\max\left\{ h^{2 H_s} ,\;\mathcal N^{-1}_\Gamma (s \mid t;h) \right\},  \max\left\{ h^{2 H_t} ,\;\mathcal N^{-1}_\Gamma (t \mid s;h) \right\} \right]  O_{\mathbb P}(1),
\end{equation}
uniformly with respect to $h\in\mathcal H_N$. The arguments are similar to show \eqref{th_q1} in the proof of Theorem \ref{thm:optim_mu} and are thus omitted.

To complete the proof of Theorem \ref{thm:optim_gamma}, it remains to show that, in the independent design case,  $   \widehat \gamma_N (s,t;h) - \widetilde \gamma_N (s,t)   $  is the leading term of $    \widehat \Gamma_N (s,t;h) - \widetilde \Gamma_N (s,t)  $ with $\widetilde \Gamma_N (s,t) = \widetilde \gamma_N (s,t) - \widetilde \mu_N (s)\widetilde \mu_N (t)$. Without loss of generality, we can consider $\mu(t)=0$, for all $t\in \mathcal T$. If this is not the case, we replace $X^{(i)}_t$ by $X^{(i)}_t - \mu(t)$ for the theory. In the adaptive procedure, the mean is not supposed known, and is estimated nonparametrically. We can write $\widehat \mu_N (t;h) = \widetilde \mu_N (t) + \left\{\widehat \mu_N (t;h)- \widetilde \mu_N (t) \right\}$ and using Theorem \ref{thm:optim_mu}, we deduce that it exists a constant $C_1>0$, that may depend on $t$, such that 
$$
\widehat \mu_N (t;h) - \widetilde \mu_N (t)  =\{ h^{ H_t} + (N\mathfrak m h)^{-1/2}\} \{C_1+o_\PP(1)\},
$$ 
where the convergence in probability $o_\PP(1)$ is uniform with respect to $h\in\mathcal H_N$.
Next, since $\lvert \widetilde \mu_N (s)\rvert$ and $\lvert\widetilde \mu_N (t) \rvert = O_\PP( N ^{-1/2})$, we deduce that two constants $C_2, C_3 > 0$ exist, that may depend on $s$ and $t$, such that
\begin{multline}
	\widehat \mu_N (s;h)\widehat \mu_N (t;h) - \widetilde \mu_N (s) \widetilde \mu_N (t) \\= O_\PP   \left( N ^{-1/2}\right) \times
	\left\{ h^{ H_s} + h^{ H_t} + (N\mathfrak m h)^{-1/2} \right\} \times   \{C_2+o_\PP(1)\} \\ +   \left\{ h^{ H_s} + (N\mathfrak m h)^{-1/2} \right\} \left\{ h^{ H_t} + (N\mathfrak m h)^{-1/2} \right\}  \times   \{C_3+o_\PP(1)\},
\end{multline}
uniformly with respect to $h\in\mathcal H_N$. Gathering facts, we deduce 
\begin{equation}
 \widehat \Gamma_N (s,t;h) - \widetilde \Gamma_N (s,t) = \{\widehat \gamma_N (s,t;h) - \widetilde \gamma_N (s,t) \}+
o_{\mathbb P} \left(  (N\Mmu^2)^{- \; \frac{ H(s,t)}{2 \{H(s,t)+1\}} } +   (N\Mmu)^{- \; \frac{ H(s,t)}{2 H(s,t)+1} } \right),
\end{equation}
with $H(s,t) = \min\{H_s,H_t\}$, and thus
\begin{equation}
\widehat \Gamma_N (s,t;h) -  \Gamma (s,t)= \{\widehat \gamma_N (s,t;h) - \widetilde \gamma_N (s,t) \} +
	o_{\mathbb P} \left(  (N\Mmu^2)^{- \; \frac{ H(s,t)}{2 \{H(s,t)+1\}} } +   (N\Mmu)^{- \; \frac{ H(s,t)}{2 H(s,t)+1} } +  N^{-1/2}\right),
\end{equation}
uniformly with respect to $h\in\mathcal H_N$.\end{proof}

%%%%%%%%%%%%%%%%%%%%%%%%%%%%%
%%%%%%%%%%%%%%%%%%%%%%%%%%%%%

\begin{proof}[Complements for the proof of Lemma \ref{lemma_mfbm}] %5.1]
Here, we provide a formal justification for the following property: for any $t$ and $u, v \in\Ostar(t)$, we have
\begin{equation}\label{loc_reg_m_SM}
\EE\left[ (X_u - X_v)^2 \right] \approx \{A^\prime (t)\}^{2H_t} |u-v|^{2H_t}.
\end{equation}
The precise meaning of this approximation of the second order moments of the increments is described in \ref{H:equivalent}. 
First, it exists a constant $C$, such that, for all $ u,v \in \Ostar(t)$,
\begin{equation}\label{loc_reg_m_SM_1}
0\leq  \frac1{2}- D(H_{u}, H_{v} )\leq C | H^\prime_t|^2 |u-v|^{2}.
\end{equation}
To prove this double inequality, note that the map $(x,y)\mapsto D(x,y)$ admits partial derivatives of any order on $(0,1)\times (0,1)$. Let
$$
g(x) = \log(\Gamma (2x+1)) + \log( \sin(\pi x)) \eqqcolon  g_1(x) - g_2(x).
$$
We notice that $g^{\prime\prime}(x) < 0$, for any $x\in (0,1)$. Indeed, using the expression of the derivative of the digamma function \citep[p.~260]{AS1964}, we have
$$
g^{\prime\prime}(x) = 4\sum_{k\geq 0} \frac1{(2x+1+k)^2} - \frac{\pi^2}{\sin^2 (\pi x)} 
= g_1^{\prime\prime}(x) - g_2^{\prime\prime}(x).
$$
We deduce that $g^{\prime\prime}$ is decreasing on $[1/2, 1)$ and, since $g^{\prime\prime}(0+) = -\infty$, the function $g^{\prime\prime}_1$ is decreasing on $(0,1/2]$ with
$$
g^{\prime\prime}_1(0) =2\pi^2/3, \;\; g^{\prime\prime}_1(1/4) =4 \{\pi^2/2 -1\}, \;\; g^{\prime\prime}_1 (1/2) = 4 \{\pi^2/6 - 1\} ,
$$
and the function $g^{\prime\prime}_2$ is decreasing on $(0,1/2]$ with
$$
g_2^{\prime\prime}(0+) =\infty, \;\; g_2^{\prime\prime}(1/4) =2\pi^2, \;\; g^{\prime\prime}_2(1/2) = \pi^2 ,
$$
we conclude that $g^{\prime\prime}<0$ on $(0,1]$. In other words, $x\mapsto g(x)$ is log-concave, and for all $0<x\neq y \leq 1$,
$$
2D(x,y)=\frac{\sqrt{\exp(g(x))\times \exp(g(y))}} {\exp(g((x+y)/2))} <1.
$$
The left-hand side of \eqref{loc_reg_m_SM_1} follows. Next, since $2D(x,x)\equiv 1$, we deduce that,  for any $x\in(0,1)$, the first order derivative of  $y\mapsto D(x,y)$ is equal to zero at $y=x$. By Taylor expansion, given a small value $r>0$, a constant $C_{x,r}$, depending on $x$ and $r$, then exists such that, for all $0\leq \lvert x-y \rvert \leq r$,
$$
 \frac1{2}- D(x, y) = D(x,x) -  D(x,y) \leq C_{x,r} \lvert x-y\rvert^2.
$$
Finally, using the fact that $ \lvert H_u-H_v \rvert \approx \lvert H^\prime _u \rvert \lvert u-v \rvert$ when $u-v$ is close to zero, we deduce the right-hand side of  \eqref{loc_reg_m_SM_1}. 

For any $t$ and $u, v \in\Ostar(t)$,  let us write
\begin{align*}
\EE[ (X_u - X_v)^2 ] &= \EE[ X_u^2]  +\EE[ X_v^2 ] -2\EE[ X_u X_v ]   \\ 
 &\hspace{-4em} = A(u)^{2H_u}+A(v)^{2H_v}- 2 D(H_u,H_v )\left[ A(u)^{H_u+H_v} +  A(v)^{H_u+H_v} - \lvert A(v)-A(u) \rvert^{H_u+H_v}\right] \\
 &\hspace{-4em}= \left\{  A(u)^{2H_u} - 2 D(H_u,H_v ) A(u)^{H_u+H_v}  \right\} + \left\{  A(v)^{2H_v} - 2 D(H_u,H_v ) A(v)^{H_u+H_v}  \right\} \\ 
 &+ 2 D(H_u,H_v )  \lvert A(v)-A(u) \rvert^{H_u+H_v} \\
 &\hspace{-4em}\eqqcolon  D_1(u \mid v) + D_1(v \mid u)+2D_2(u,v).
\end{align*}
Next, let $\mathcal T \subset [0,\infty)$ be a compact interval, and for any real-valued function $B$ defined on $\mathcal T$, let $\|B\|_{\mathcal T,\infty} = \sup_{t\in\mathcal T} B(t)$. In the case $t>0$, where $A$ stays away from zero on $\Ostar(t)$, we can write
$$
D_1(u \mid v)  = A(u)^{2H_u} -  A(u)^{H_u+H_v} + R_1(u \mid v),
$$
with 
\begin{equation}
\left\lvert R_1(u \mid v)  \right\rvert \leq  \left\{ 1- 2D(H_u,H_v)\right\}
\|A^H\|^2_{\mathcal T,\infty}\leq 
C\|A^H\|^2_{\mathcal T,\infty} \|H^\prime \|^2_{\mathcal T,\infty}\lvert u-v \rvert^2=
O(\lvert u-v\rvert^2) ,
\end{equation}
and
\begin{align}
A(u)^{2H_u} -  A(u)^{H_u+H_v} 
&= A(u)^{2H_u} \left[ 1- \exp\{(H_v-H_u)\log(A(u))\}  \right] \\
&= A(u)^{2H_u}  \left[ -  (H_v-H_u)\log(A(u)) +O(\lvert u-v \rvert^2)  \right] \\ 
&= A(u)^{2H_u}  \left[- H_u^\prime (v-u)\log(A(u)) +O(\lvert u-v \rvert^2)  \right] .
\end{align}
The term  $D_1(v \mid u) $ decomposed similarly, and we thus deduce
\begin{align}
 D_1(u \mid v)+D_1(v \mid u) &= (v-u) \left[A(v)^{2H_v}H_v^\prime\log(A(v))  - A(u)^{2H_u}H_u^\prime \log(A(u))  \right] +O(\lvert u-v\rvert^2)  \\
 &= O(\lvert u-v\rvert^2) .
\end{align}
The last equality is due to the fact that, by assumptions,  the map $v\mapsto A(v)^{2H_v}H_v^\prime\log(A(v))$ is continuously differentiable over $\Ostar(t)$. On the other hand, by \eqref{loc_reg_m_SM_1},
\begin{align}
D_2(u,v) &= \lvert A(v)-A(u)\rvert^{H_u+H_v} + \sup_{t\in\mathcal T} \lvert H^\prime _t\rvert \times o(\lvert u-v \rvert^2)\\ 
&= \left\lvert A^\prime (t)(v-u) + O(\lvert u-v\rvert^2) \right\rvert^{2H_t+2H^\prime_t (v-u) +  O(\lvert u-v \rvert^2) } +o(\lvert u-v \rvert^2)\\
&=  \lvert A^\prime (t)\rvert^{2H_t}\lvert v-u \rvert ^{2H_t}  \times\left\{ 1+ O\left( \lvert u-v\rvert^{\min_{t\in \mathcal T}H_t} \right)\right\} +o(\lvert u-v \rvert^2).
\end{align}
Assumption \ref{H:equivalent} is thus satisfied with any  $0 < \beta \leq \min\left\{\min_{t\in\mathcal T} H_t, 2(1-\max_{t\in\mathcal T} H_t) \right\}$.

In the case $t=0$, we only have to investigate the case $A(0)=0$. Whenever  $A(0)>0$, the previous arguments apply without any change. If $A(0)=0$, we only have to revisit the arguments for bounding 
$D_1(u \mid v)+D_1(v \mid u) $. The map  
$$
v\mapsto \zeta(v)=A(v)^{2H_v}H_v^\prime\log(A(v)),
$$ 
is no longer differentiable over $\Ostar(0)$, if $H_v \leq 1/2$. However, this map is Hölder continuous and we still have  
\begin{equation}\label{eq:Hold_l}
 D_1(u \mid v)+D_1(v \mid u) =  O\left(\lvert u-v \rvert^{1+\gamma}\right), \quad \text{for any}\quad 0< \gamma < \min(1,  2\min_{v\in\Ostar(0)}H_v),
\end{equation}
 and, given the regularity conditions imposed on the map $u\mapsto H_u$, this suffices to complete the proof of Lemma \ref{lemma_mfbm}. With the rule $0\log(0)=0$,  the Hölder continuity we use is  
\begin{equation}\label{reztp}
\sup_{0\leq u<v\leq \Delta_*/2} \frac{\left\lvert \zeta(u)- \zeta(v)\right\rvert }{\lvert u-v\rvert^{\gamma} }\leq C<\infty ,
\end{equation}
for some $C$ depending on $\gamma$, $H$ and $\Delta_*$ and decreasing to zero with $\Delta_*$. 
Under our assumptions, we have 
$$
\zeta (v)  = [ A^\prime(0)v ]^{2H_0}  \times \zeta_1 (v) \times H_0^\prime \times  \{\log(v) + \log(A^\prime(0)\}   \times   \{1+o(1)\}, 
$$
with
$$
 \zeta_1 (v)  = v^{2H_0^\prime v \{1+o(1)\}}.
$$
All the  $o(1)$ terms in the last display can be uniformly bounded, with respect to $v\in\Ostar(0)$, by a constant times $\Delta_*$, this constant depending only on the bounds of  $A^\prime$, $\lvert H^\prime \rvert$ and $\lvert H^{\prime\prime} \rvert$ near the origin. 
Let us first notice that for any $0<\gamma < \min(1, 2H_0)$, it exists a constant $c$  such that 
\begin{equation}\label{reztp2}
\sup_{0\leq u<v\leq \Delta_*/2} \frac{\left\lvert u^{2H_0} \log(u)- v^{2H_0} \log(v)\right\rvert }{\lvert u-v \rvert^{\gamma} }\leq c<\infty.
\end{equation}
Moreover, we have $\zeta_1(0+) = 1$,  $\zeta_1$ is bounded on $\Ostar(0)$, and for  $0<u<v\leq \Delta_*/2$, we  have 
\begin{align*}
\lvert \zeta_1 (v) -  \zeta_1 (u) \rvert &= \lvert \exp( 2H_0^\prime v \log(v) \{1+o(1)\})  - \exp( 2H_0^\prime u \log(u) \{1+o(1)\}) \rvert\\
 &\leq 2\lvert H_0^\prime\rvert  \lvert v \log(v) -  u \log(u) \rvert \leq 2c_1\lvert H_0^\prime\rvert \lvert u-v\rvert^{\gamma},
\end{align*}
for some constant $c_1$. Gathering facts, we deduce \eqref{reztp}. The justification of Lemma \ref{lemma_mfbm} is now complete. 
\end{proof}

% !TeX root=bj-sample-SM.tex

\section{Details on some equations from the main text}\label{sec:details_sm}

\subsection{Discussion of the choices of $t_1$, $t_2$ and $t_3$ in the definition (\ref{eq:tilde-Hd})}

The following discussion is inspired by a comment from a Reviewer. A first step for the construction of the local regularity estimator is the definition of a proxy value:
\begin{equation}\label{eq:tilde-Hd_SM}
	\widetilde H_{t} =  \frac{\log(\theta(t_1,t_3)) - \log(\theta(t_1,t_2))}{2\log(2)} \quad \text{if }  \Delta_* \text{ is small}.
\end{equation}
See \eqref{eq:tilde-Hd} in the main text. This definition is based on the simple choice of $t_1$, $t_2$ and $t_3$  such that 
$$
\lvert t_3-t_1 \rvert = 2 \lvert t_2-t_1 \rvert.
$$
We can more generally proceed as follows: let $t_1$ and $t_3$ be such that $[t_1, t_3] \in \mathcal{O}_{*}(t)$ and define the proxy
\begin{equation}
	\widetilde{H}_t = \frac{\log(\theta(t_1, t_3)) - \log(\theta(t_1, t_2))}{2\left\{\log(\lvert t_3 - t_1 \rvert) - \log(\lvert t_2 - t_1 \rvert)\right\}} \approx H_t,
\end{equation}
and the corresponding estimator of the exponent
\begin{equation}
	\widehat{H}_t = \frac{\log(\widehat{\theta}(t_1, t_3)) - \log(\widehat{\theta}(t_1, t_2))}{2\left\{\log(\lvert t_3 - t_1 \rvert) - \log(\lvert t_2 - t_1 \rvert)\right\}}.
\end{equation}
In practice, the choice of $t_1$, $t_2$ and $t_3$, which, here, has to be the same for all curves, can be guided by the density of the design points. The practical investigation of these aspects is left for future work. 

Finally, one can also consider a nearest neighbors idea for the choice of $t_1$, $t_2$ and $t_3$. With an independent design,  these values then become random. This idea was investigated by \cite{golo2020} and leads to alternative estimates of the local regularity which do not require preliminary smoothing. Even if it offers an elegant alternative, which avoids the choice of a smoothing parameter such as the bandwidth, the idea based on nearest neighbors idea could require more points $\Tnm$ on each curve. See discussion in the paragraph following \citet[Th.~1]{golo2020}. Moreover, the extension of the nearest neighbors idea to the regularity estimation in the case of differentiable sample paths is much more challenging.

\subsection{Details on the approximation (\ref{eq:ci2b})} %(17)}

Recall that \begin{equation}\label{eq:ci_aa}
c_i (t;h,\alpha )  = \sum_{m=1}^{M_i} \left\lvert (\Tnm - t)/h\right\rvert^{\alpha} \left\lvert W_{m}^{(i)}(t;h)\right\rvert,
\end{equation}
and, using the Nadaraya-Watson (NW) estimator,
\begin{equation}\label{eq:ci2_aa} 
 \overline{C} (t;h,\alpha)= \!
\frac{1}{\WN(t;h)  } \sum_{i=1}^{N}w_i (t;h) c_i (t;h,\alpha ).
\end{equation}
When using the NW estimator, for each $1\leq i \leq N$,
\begin{equation}\label{eq:ci2_aa2}
c_i (t;h,\alpha )  = \frac{1}{\widehat g^{(i)} (t)} \frac{1}{M_ih} \sum_{m=1}^{M_i} \left\lvert (\Tnm - t)/h\right\rvert^{\alpha} K\left((\Tnm - t)/h\right),
\end{equation}
with 
$$
\widehat g^{(i)} (t)= \frac{1}{M_ih} \sum_{m=1}^{M_i}   K\left((\Tnm - t)/h\right)\approx g(t).
$$
Here, $g$ denotes the density of the $\Tnm$. By a standard change of variables, 
$$
\EE[c_i (t;h,\alpha )\widehat g^{(i)} (t)] \approx g(t)\int|u|^{ \alpha}K(u)du.
$$
and this explains our proposal
\begin{equation}\label{eq:ci2b_aa}
\overline{C} (t;h,\alpha)\approx \int|u|^{ \alpha}K(u)du,
\end{equation}
for the NW estimator. The same arguments apply for $ \overline{\mathfrak{C}} (t\mid s;h,\alpha)$ used for estimating the covariance function. 
In the case of a local linear estimator, it suffices to use the equivalent kernels for local polynomial smoothing. Approximation \eqref{eq:ci2b_aa} could remain the same in the local linear case, but has to be changed for higher-order polynomials.  See \citet[Sec.~3.2.2]{fan_local_1996}.

\subsection{Details on the definition (\ref{eq:opt_h_cov0c})} %(29)}

Recall that  $\widetilde \gamma_N (s,t) = N^{-1} \sum_{i=1}^N \Xp{i}_s \Xp{i}_t $. Here, $\WN$ and $w_i$ are short notations for  $\WN(s,t;h) $ and $w_i (s,t;h)$, respectively. Moreover, 
$
\hXni_t - \Xp{i}_t =  
B_t^{(i)} +
 V_t^{(i)} , 
$
where  $B_t^{(i)}\coloneq\EEi[\hXni_t ]  -  \Xp{i}_t  $ and $ V_t^{(i)} : = \hXni_t - \EEi [\hXni_t ] .$ Let us define
$$
 \widetilde \gamma_W (s,t;h)   =  \frac{1}{\WN }\sum_{i=1}^N w_i \Xp{i}_s \Xp{i}_t .
$$
To explain our empirical risk bound  $\mathcal{R}_\Gamma (s \mid t;h) $ defined in \eqref{eq:opt_h_cov0c}, let us write
\begin{align*}%\label{detail_ga}
   \widehat \gamma_N (s,t;h) - \widetilde \gamma_W (s,t;h)  &= 
 \frac{1}{\WN }\sum_{i=1}^N w_i \{\hXni_s-\Xp{i}_s\} \Xp{i}_t + \frac{1}{\WN }\sum_{i=1}^N w_i \Xp{i}_s\{ \hXni_t - \Xp{i}_t\} \\
&\qquad +  \frac{1}{\WN }\sum_{i=1}^N w_i \{ \hXni_s - \Xp{i}_s\}\{ \hXni_t - \Xp{i}_t\}\\
&=  \frac{1}{\WN }\sum_{i=1}^N w_i \{ B_s^{(i)} \Xp{i}_t +\Xp{i}_s B_t^{(i)}\} \\
&\qquad + \frac{1}{\WN }\sum_{i=1}^N w_i \{ V_s^{(i)} \Xp{i}_t +\Xp{i}_s V_t^{(i)}\} \\
&\qquad + \frac{1}{\WN }\sum_{i=1}^N  w_i \{ B_s^{(i)}  B_t^{(i)}+V_s^{(i)}  V_t^{(i)} \} \\
&\qquad + \frac{1}{\WN }\sum_{i=1}^N w_i \{ B_s^{(i)}  V_t^{(i)}+V_s^{(i)}  B_t^{(i)} \}.
\end{align*}
By construction, for all $1\leq i,j \leq N$,
$$
\EEMT[V_s^{(i)}  B_t^{(j)}] =\EEMT[B_s^{(i)}  V_t^{(j)}] =0.
$$
Moreover, whenever $2h<\lvert s-t \rvert$, we have, for all $1\leq i,j\leq N$,
$$
  \EEMT[V_s^{(i)}  V_t^{(j)}] = 0.
$$
Using  these properties, the inequality $(a+b)^2 \leq 2(a^2+b^2)$, and repeated application of Cauchy-Schwarz inequality to check the negligible terms, we deduce
\begin{align*}%\label{detail_ga2}
 \EEMT\left[ \left\{ \widehat \gamma_N (s,t;h) - \widetilde \gamma_W (s,t;h)  \right\}^2 \right] 
&= 
  \EEMT\left[ \left\{  \frac{1}{\WN }\sum_{i=1}^N w_i \left( B_s^{(i)} \Xp{i}_t +\Xp{i}_s B_t^{(i)}\right)\right\}^2 \right]
  \\  
  &\hspace{-10em}+
 \EEMT\left[ \left\{  \frac{1}{\WN }\sum_{i=1}^N w_i \left\{ V_s^{(i)} \Xp{i}_t +\Xp{i}_s V_t^{(i)}\right\} \right\}^2 \right] + \text{ negligible terms}
 \\ 
 &\hspace{-12em}\leq
2 \EEMT\left[ \left\{  \frac{1}{\WN }\sum_{i=1}^N w_i  B_s^{(i)} \Xp{i}_t \right\}^2 \right]
 + 2\EEMT\left[ \left\{  \frac{1}{\WN }\sum_{i=1}^N w_i \Xp{i}_s B_t^{(i)}\right\}^2 \right]
 \\ 
 &\hspace{-10em}+ \frac{1}{\WN^2  } \sum_{i=1}^N w_i   \EEMT \left[ \left\{  V_s^{(i)} \Xp{i}_t \right\}^2 +\left\{ \Xp{i}_s V_t^{(i)} \right\}^2  \right]+ \text{ negligible terms}\\ 
&= \{G_{1}(s \mid t) +G_{1}(t \mid s) + G_2\}\{1+o_{\mathbb P}(1) \}.
\end{align*}
We can now write
\begin{align*}%\label{detail_ga2c}
\EEMT \left[ \left\{ \Xp{i}_s V_t^{(i)} \right\}^2\right] 
&= \EEMT \left[ \left\{ \Xp{i}_s \right\}^2 \left\{ \sum_{m=1}^{M_i} \enm  W_{m}^{(i)}(t;h) \right\}^2 \right] \\
&=  \EEMT \left[ \left\{ \Xp{i}_s \right\}^2 
  \sum_{m=1}^{M_i}\EEi \left[ \left\lvert \enm \right\rvert^2 \right]\left\lvert W_{m}^{(i)}(t;h) \right\rvert^2  \right] 
\\
&\leq  \sigma_{\max}^2   m_2(s)  \left\{ \max_{m} \left\lvert W_{m}^{(i)}(t)\right\rvert\times \sum_{m=1}^{M_i} \left\lvert W_{m}^{(i)}(t;h)\right\rvert   \right\},
\end{align*}
where $m_2(s) = \EE\left[  X_s^2  \right]$ and $\EEi [\cdot] = \EE [\cdot\mid M_i,  \mathcal T_{obs} ^{(i)}  , \Xp{i}] $.  
Let us recall that
\begin{equation}\label{eq:Ncal2def1_b}
    \mathcal N_{i} (t;h) =  \frac{w_i (s,t;h)}{\max_{1\leq m\leq M_i} \lvert W_{m}^{(i)}(t;h)\rvert}.
    % \quad \text{ and } \quad \mathcal N_{i} (s|t;h) =  \frac{w_i (s,t;h)}{\max_{1\leq m\leq M_i} |W_{m}^{(i)}(s;h)|}.
\end{equation}
We deduce 
$$
G_2 \leq \frac{\sigma_{\max}^2}{\WN^2}\sum_{i=1}^N  w_i \left[m_2 (t) \frac{c_i (s;h)}{ \mathcal N_{i} (s;h) } +m_2 (s)\frac{c_i (t;h)}{  \mathcal N_{i} (t;h) }\right],
$$
where the $c_i(t;h)$ are defined by equation~\eqref{eq:ci}
and the   $\mathcal N_{i} (s;h)$  and $\mathcal N_{i} (t;h)$  are defined using~\eqref{eq:Ncal2def1_b}.

To bound the terms related to the bias of $\hXni_t$, using the moment assumptions, 
by the law of large numbers, dominated convergence theorem,  we can write
\begin{align*}%\label{detail_ga2b1}
G_{1}(s \mid t)+ G_{1}(t \mid s) 
&\leq 2 \EEMT \left[\frac{1}{\WN  } \sum_{i=1}^N w_i   \left\lvert B_t^{(i)} \right\rvert^2  \times  \left\{m_2(s) +  \frac{1}{\WN  } \sum_{i=1}^N w_i  \left(  \left\lvert \Xp{i}_s \right\rvert^2 -m_2(s) \right) \right\} \right]  \\
&\hspace{-2em}+ 2 \EEMT \left[\frac{1}{\WN  } \sum_{i=1}^N w_i   \left\lvert B_s^{(i)} \right\rvert^2  \times \left\{m_2(t) +  \frac{1}{\WN  } \sum_{i=1}^N w_i  \left(  \left\lvert \Xp{i}_t \right\rvert^2 -m_2(t) \right) \right\} \right]\\
&\hspace{-5em}=  2\left\{\frac{1}{\WN  } \sum_{i=1}^N w_i  \EEMT  \left[ \left\lvert B_s^{(i)} \right\rvert^2 \right]\right\}m_2(t)  + 2\left\{\frac{1}{\WN  } \sum_{i=1}^N w_i  \EEMT  \left[ \left\lvert B_t^{(i)} \right\rvert^2 \right]\right\}m_2(s) \\
&\hspace{20em}+ \text{ negligible terms}\\
&\hspace{-5em}\leq h^{2H_t} \times 2\left\{   m_2(t)\overline{\mathfrak{C}} (s \mid t;h,2  H_s)    L_{s}^2 +  m_2(s)  \overline{\mathfrak{C}} (t \mid s;h,2  H_t)     L_{t}^2 \right\} \{1+ o_{\mathbb P}(1)\},
\end{align*}
where  $\overline{\mathfrak{C}} (t \mid s;h,  \cdot )$ is defined according to~\eqref{eq:cbar_g}.

Gathering facts, we deduce that 
\begin{align}\label{detail_ga2f}
 \EEMT\left[ \left\{ \widehat \gamma_N (s,t;h) - \widetilde \gamma_W (s,t;h)  \right\}^2 \right]  &\leq  h^{2H_t} \times  2  \left\{m_2(t)\overline{\mathfrak{C}} (s \mid t;h,2  H_s)    L_{s}^2
+  m_2(s)\overline{\mathfrak{C}} (t \mid s;h,2  H_t)  L_{t}^2\right\} \\
&\hspace{-5em}+ \frac{\sigma_{\max}^2}{\WN^2}\sum_{i=1}^N  w_i \left\{m_2 (t) \frac{c_i (s;h)}{ \mathcal N_{i} (s;h) } +m_2 (s)\frac{c_i (t;h)}{  \mathcal N_{i} (t;h)  }\right\} + \text{ negligible terms}.
\end{align}
On the other hand, we have 
\begin{align}\label{eq:Ncal2b}
 \EEMT\left[ \left\{ \widetilde \gamma_N (s,t)   - \widetilde \gamma_W (s,t;h)  \right\}^2 \right]
 &=\frac{ \operatorname{Var}[X_sX_t]}{\WN^2} \sum_{i=1}^N  \left\{w_i   - \frac{\WN }{N}  \right\}^2 \\ &= \operatorname{Var}[X_s X_t]  \left\{  \frac{1}{\WN } - \frac{1}{N } \right\}.
\end{align}
It remains to note that 
\begin{align}
 \EEMT\left[ \left\{ \widehat \gamma_N (s,t;h) - \widetilde \gamma_N (s,t)  \right\}^2 \right] 
&\leq 2 \EEMT\left[ \left\{ \widehat \gamma_N (s,t;h) - \widetilde \gamma_W (s,t;h)  \right\}^2 \right] \\
&\quad\quad + 2 \EEMT\left[ \left\{ \widetilde \gamma_W (s,t;h)  - \widetilde \gamma_N (s,t)  \right\}^2 \right] .
\end{align}

% !TeX root=bj-sample-SM.tex
\section{A Central Limit Theorem for the mean function estimator}\label{sec:TCL}  

Let us now study the  asymptotic distribution of our adaptive pointwise estimation of the mean function. In standard nonparametric curve estimation, if the regression function is Hölder continuous with some exponent $H\in(0,1)$, the optimal bandwidth  that balances the squared bias and the variance for pointwise estimation has a rate equal to  $n^{-1/(2H + 1)}$, where $n$ is the sample size.  In that context, under appropriate conditions, it is possible to derive the convergence in distribution of the kernel regression estimator to a normal limit, with the rate $n^{-H/(2H + 1)}$. In general, the expectation of the normal limit distribution is not zero due to the bias term. Note that considering a  bandwidth rate slower than $n^{-1/(2H + 1)}$ will cause this convergence in distribution to fail, as the bias term multiplied by the square root of the sample size times the bandwidth will tend to infinity. Knowing a lower bound on the regularity of the curve to be estimated allows for asymptotic distribution results which are robust to failures due to diverging bias terms.

In the context of functional data, the regularity of the mean function is necessarily equal to or larger than that of the sample paths. One consequence is that the minimax optimal rate for the mean function estimation is given by the sample path regularity, see \cite{cai2011}. This implies that, from the minimax optimality perspective, the rate of convergence in distribution for the mean function estimation must take into account the regularity of the sample paths. In the following, we follow this idea and provide the asymptotic normality of our estimator, but with a bandwidth decrease rate slightly faster than $(N\mathfrak{m})^{-1/(2H_t + 1)}$ in order to make the bias term negligible. For simplicity, we only consider the independent design case.

\begin{theorem}\label{thm:mean-clt}
	Let $t\in \mathcal T$. Assume that the conditions of Theorem \ref{thm:optim_mu}
	hold true, that the conditional variance $\sigma^2(t,x)$ of the measurement errors does not depend on $x$, denoted $\sigma^2(t)$, and that the map $t\mapsto \sigma(t)>0$ is Hölder continuous. Assume $\widehat \mu (t;h)$ is defined as in \eqref{mu1_d} with $K$ a  bounded kernel which is bounded away from zero on the  support $[-1,1]$.  Let  $h_{N}\in\mathcal H_N$, $N\geq 1$, such that 
	\begin{equation}\label{h_as_norma}
		h_N   (N\mathfrak m)^{1/(2H_t+1)} \rightarrow 0  . 
	\end{equation}
	Moreover, 
	\begin{equation}\label{eq:sigma_clt}
	V_N (t;h_N)\coloneq	\frac{\sigma^2(t)}{\mathcal W_N(t;h_N)}\sum_{i=1}^{N} w_i(t;h_N)  \sum_{m=1}^{M_i}\left[W_{m}^{(i)}(t;h_N)\right]^2 
		\overset{\PP}{\longrightarrow} \Sigma(t) \in [0,\infty).
	\end{equation}
	Then 
	\begin{equation*}
		\sqrt{\mathcal W_N(t;h_N)}\left\{\widehat \mu_N(t;h_{N}) - \mu(t)\right\}\overset{d}{\longrightarrow} \mathcal{N}\left(0, \operatorname{Var} (t) +\Sigma(t) \right).
	\end{equation*}
\end{theorem}

Condition \eqref{h_as_norma} renders the bias term negligible and thus avoids the usual mean correction used in kernel regression. In the dense case ($N\ll \mathfrak m ^{2H_t}$), if in addition $\mathfrak m  h_N \rightarrow \infty$, then $\Sigma(t)=0$ and the $\widehat \mu_N(t;h_{N}) $ has the same asymptotic distribution as the infeasible empirical mean function estimator $\widetilde \mu (t)$ obtained with the $X^{(i)}_t$, $1\leq i\leq N$.  As expected, in the sparse case ($N\gg \mathfrak m  ^{2H_t}$), the rate of convergence in distribution given by $\EE[\mathcal W_N(t;h_N)]^{-1/2}$, is slower than $N^{-1/2}$. Moreover,  $\Sigma(t)>0$.

\begin{proof}[Proof of Theorem~\ref{thm:mean-clt}]
	Recall the definition 
	\begin{equation}
		\widetilde \mu_W (t;h) = \frac{1}{\mathcal W_N(t;h_N)}\sum_{i=1}^N  w_i (t;h) X_t^{(i)} .
	\end{equation}
	We can then write 
	\begin{equation}
		\widehat \mu_N(t;h) - \mu(t) = \left\{\widehat \mu_N(t;h) - \widetilde  \mu_W(t;h)\right\}+ \left\{\widetilde \mu_W(t;h) -   \mu(t)\right\}  \eqqcolon T_{N1}(t;h) + T_{N2}(t;h). 
	\end{equation}
 Using the Hölder continuity of $t\mapsto \sigma(t)$ and  the arguments used to justify Theorem~\ref{thm:optim_mu}, since $h_N$ decreases to zero faster than the bandwidth which realizes the balance between the squared bias and the variance, and $K(\cdot)$ is bounded and bounded away from zero, we have 
	\begin{equation}
T_{N1}(t;h_N) = \{1+O (h)\}\frac{\sigma (t)}{\mathcal W_N(t;h_N)}  \sum_{n=1}^N   w_i(t;h_N) \sum_{m=1}^{M_i}W_{m}^{(i)}(t;h_N) \unm \eqqcolon   \frac{1+ O(h)}{\sqrt{\mathcal W_N(t;h_N)}}  \times \mathcal U_N (t;h).
	\end{equation}
A bounded $K$, bounded away from zero, is a convenient assumption to guarantee that for all $h, i, m$,
$$\max_m W^{(i)}_m(t;h) /  \min_m W^{(i)}_m(t;h)\leq C,$$
for some constant $C$,  and thus the rates on both sides of \eqref{eq:var_8} are the same. 
By Lyapunov CLT for independent variables, conditionally given the $M_i$ and $\{\Tnm,  1 \leq m \leq M_i\}$, $1\leq i \leq N$, we have 
	$$
	\frac{\mathcal U_N (t;h_N)}{\sqrt{V_N (t;h_N)}}\overset{d}{\longrightarrow} N(0,1),
	$$
	with $V_N (t;h_N) $  defined in \eqref{eq:sigma_clt}.
	This implies that, for any event for which the sequence $V_N (t;h_N)$ converges to $\Sigma(t)$, we get, for all $u\in\RR$,
	\begin{equation}\label{zelm}
		\EE_{M,T}\left[\exp\left\{-iu\sqrt{\mathcal W_N(t;h_N) }\;T_{N1} (t;h_N)\right\}\right] 
		\longrightarrow \exp(-u^2\Sigma(t) / 2).
	\end{equation}
	Note that $\EEMT[\cdots]$ on the left hand side is a bounded sequence of random variables. Since $V _N (t;h_N) - \Sigma(t) = o_\PP(1)$, and the convergence in probability is characterized by the fact that 
	every sub-sequence  has a further sub-sequence which converges almost surely, we deduce that the convergence in \eqref{zelm} holds in probability. By the Dominated Convergence Theorem for a sequence of bounded random variables that converges in probability, we get, for all $u\in\RR$,
	\begin{equation}\label{zelm2}
		\EE \left[\exp\left\{-iu\sqrt{\mathcal W_N(t;h_N) }\;T_{N1} (t;h_N)\right\}\right] 
		\longrightarrow \exp(-u^2\Sigma(t) / 2),
	\end{equation}
	which means 
	\begin{equation}\label{tcl_2}
	\sqrt{\mathcal W_N(t;h_N) }\;T_{N1} (t;h_N) \overset{d}{\longrightarrow} \mathcal{N}\left(0, \Sigma(t) \right).
	\end{equation}
Applying again standard CLT, we also have 
	\begin{equation}\label{tcl_1}
		\sqrt{\mathcal W_N(t;h_N)}	T_{N2}(t;h_N) \overset{d}{\longrightarrow} \mathcal{N}\left(0,\operatorname{Var}(t) \right).
	\end{equation}
Let us now point out that,  by the independence between the noise and the sample paths of $X$, the sequences  $T_{N1} (t;h_N)$ and $T_{N2} (t;h_N)$  are independent. From this, \eqref{tcl_1} and \eqref{tcl_2}, we get  
	$$
	\sqrt{\mathcal W_N(t;h_N)}\; \{T_{N1} (t;h_N)+T_{N2} (t;h_N) \} \overset{d}{\longrightarrow} \mathcal{N}\left(0, \Sigma(t)+\operatorname{Var}(t) \right).
	$$
	Let us note that, by our results above, $\Sigma(t)=0$ if $\mathfrak mh_N \rightarrow \infty$.   \end{proof}

% !TEX root=../main_supp.tex

\section{Additional simulation results} % (fold)
\label{sub:simulation_design_SM}

Let us recall that we simulate data sets using the data generating process defined in Section \ref{lemma_NicoH} in the main text, with a Hurst index function $H_t$ and a time deformation function $A_t$ estimated on the Power Consumption data set, to which we add a mean curve also fitted to the real data set.
The estimates $\widehat{H}_t$ and the estimates of the mean and covariance functions are obtained using the same data, i.e. we did not use a \emph{learning sample} for $\widehat{H}_t$.

\begin{figure}[h!]
    \centering
 \begin{subfigure}[p!]{0.49\textwidth}
     \centering
     \includegraphics[scale=0.5]{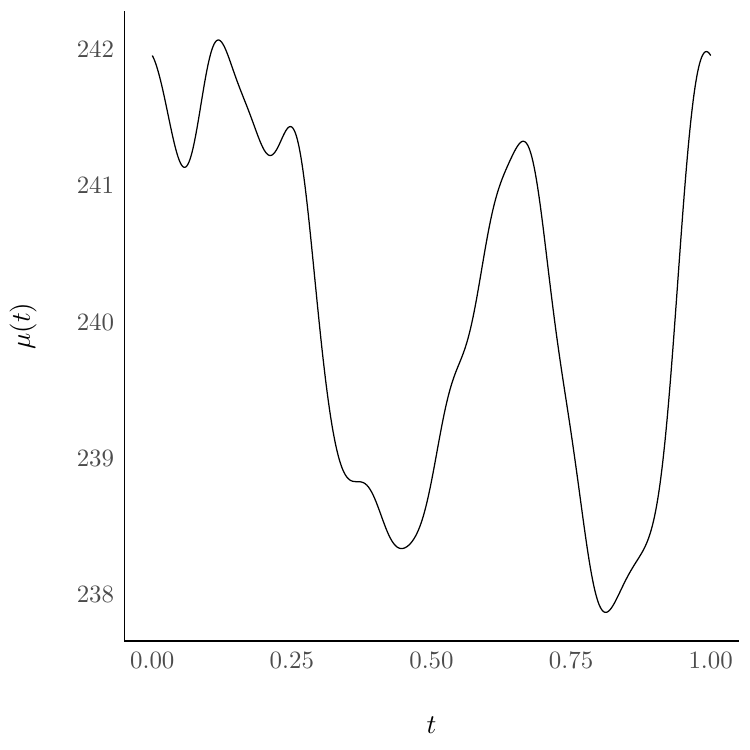}
     \caption{Mean curve $\mu(\cdot)$ for \emph{Experiment 4}.}
     \label{fig:mean_exp4}
\end{subfigure}
\hfill
\begin{subfigure}[p!]{0.49\textwidth}
     \centering
     \includegraphics[scale=0.5]{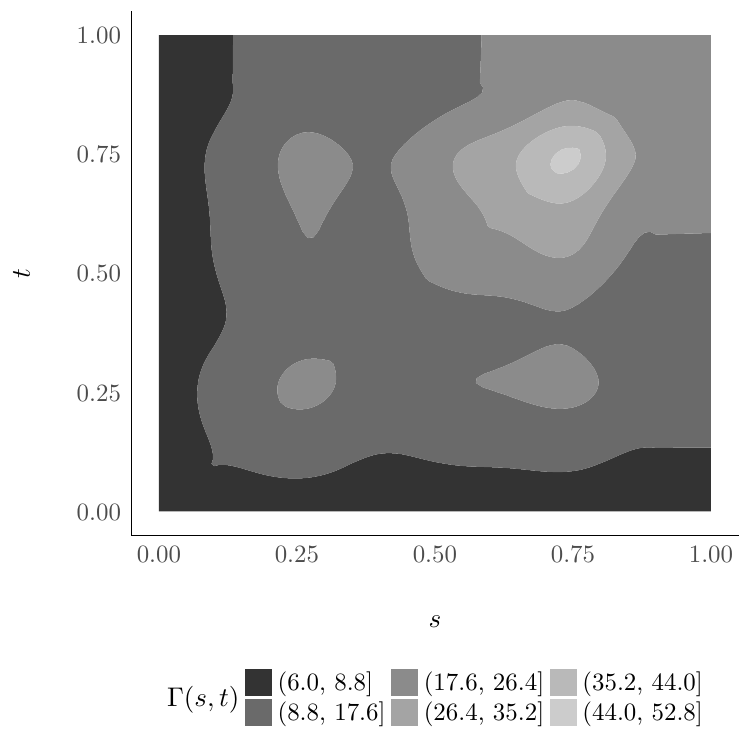}
     \caption{Covariance surface $\Gamma(\cdot, \cdot)$ for \emph{Experiment 5}.}
     \label{fig:cov_exp5}
\end{subfigure}
\caption{Description of the modifications for \emph{Experiment 4} and \emph{5}.}
\end{figure}

We consider eight experiments, each of them replicated 500 times. For each experiment, except specifically specified, we consider $N \in \{50, 100, 200\}$, $\mathfrak{m} \in \{20, 30, 40, 50\}$ and that the number of points per curve $M_i$ has a Poisson distribution with mean $\mathfrak{m}$. In \emph{Experiment 1}, we assume that the sampling points are uniformly distributed in $\mathcal{T}$, the standard deviation of the noise is $\sigma = 0.5$, the regularity of the mean function is $s = \exp(-6)$, the number of Fourier basis functions for the estimation of $H_t$ and $L_t$ is $9$, and $\varpi = 2.5$. All the other experiments are designed starting  from \emph{Experiment 1} and modifying one parameter at a time.  In \emph{Experiment 2} and \emph{Experiment 3}, we consider $\sigma = 0.25$ and $\sigma = 1$, respectively. We set $s = \exp(-3)$ for \emph{Experiment 4} resulting in a smoother mean function $\mu$ (see Figure \ref{fig:mean_exp4}). We used only $7$ functions in the Fourier basis in \emph{Experiment 5}, that is a smoother estimation of $H_t$ and $L_t$ and resulting in a smoother covariance surface $\Gamma$ (see Figure \ref{fig:cov_exp5}). For \emph{Experiment 6}, the distribution of the sampling points is a mixture of beta distributions $0.5 \mathcal{B}(1, 2) + 0.5 \mathcal{B}(2, 1)$. For \emph{Experiment 7}, we set $\varpi = 1$. Finally, in \emph{Experiment 8}, we apply our approach to the case of differentiable trajectories that we obtain by integrating the sample paths generated as in \emph{Experiment 1}.

The results from \emph{Experiment 1}, with the $\text{ISE}_{0}$ criterion, are presented in the main text. We present the results from  \emph{Experiment 1}, with the $\text{ISE}_{0.05}$ criterion, and the results of the other seven experiments below. The results for the mean function are in Section \ref{sub:mean_estimation}, while the results for the covariance function can be found in Section \ref{sub:covariance_estimation}. 

% subsection simulation_design (end)

\subsection{Mean estimation} % (fold)
\label{sub:mean_estimation}

\captionsetup[subfigure]{labelformat=empty}

\begin{figure}[h!]
     \centering
     \includegraphics[scale=0.6]{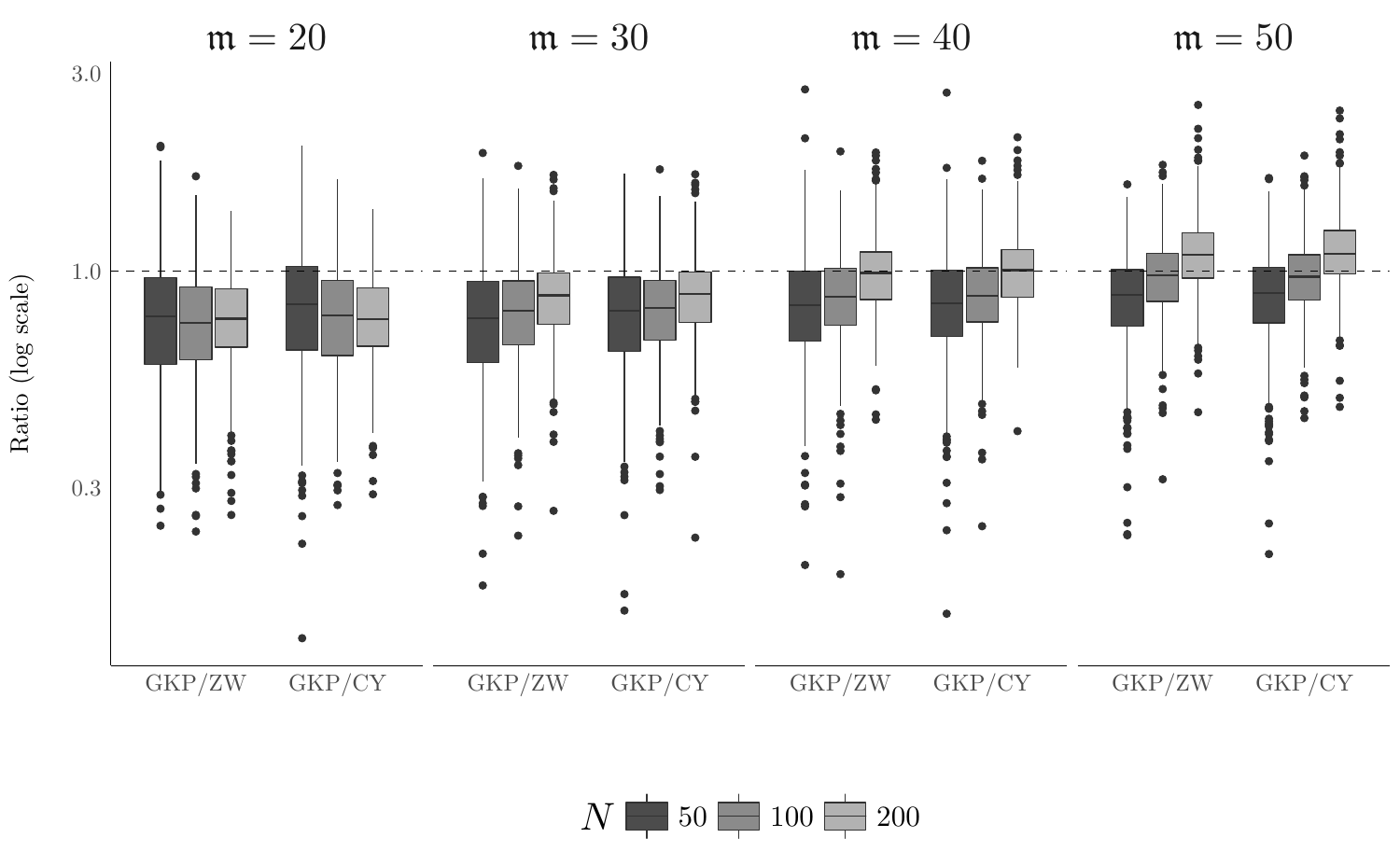}
     \caption{Results for the estimation of $\mu$ for \emph{Experiment 1}. The ratios are computed using $\text{ISE}_{0.05}$.}
     \label{fig:mean_true_exp1_ise005}
\end{figure}

\begin{figure}[h!]
	\centering
     \includegraphics[scale=0.6]{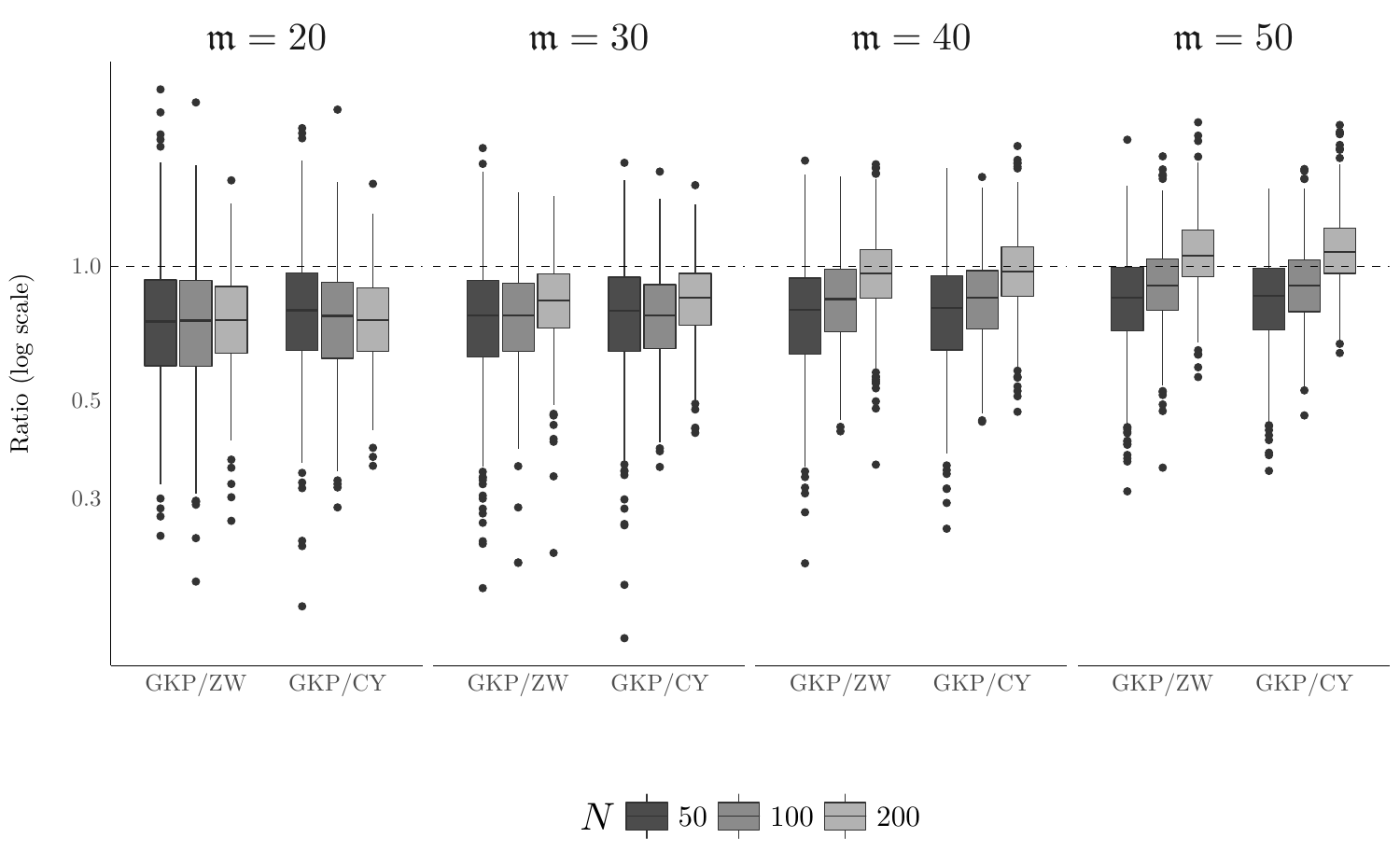}
     \caption{Results for the estimation of $\mu$ for \emph{Experiment 2} (noise std  $\sigma = 0.25$). The ratios are computed using $\text{ISE}_{0}$.}
     \label{fig:mean_true_exp2}
\end{figure}

\begin{figure}[h!]
	\centering
     \includegraphics[scale=0.6]{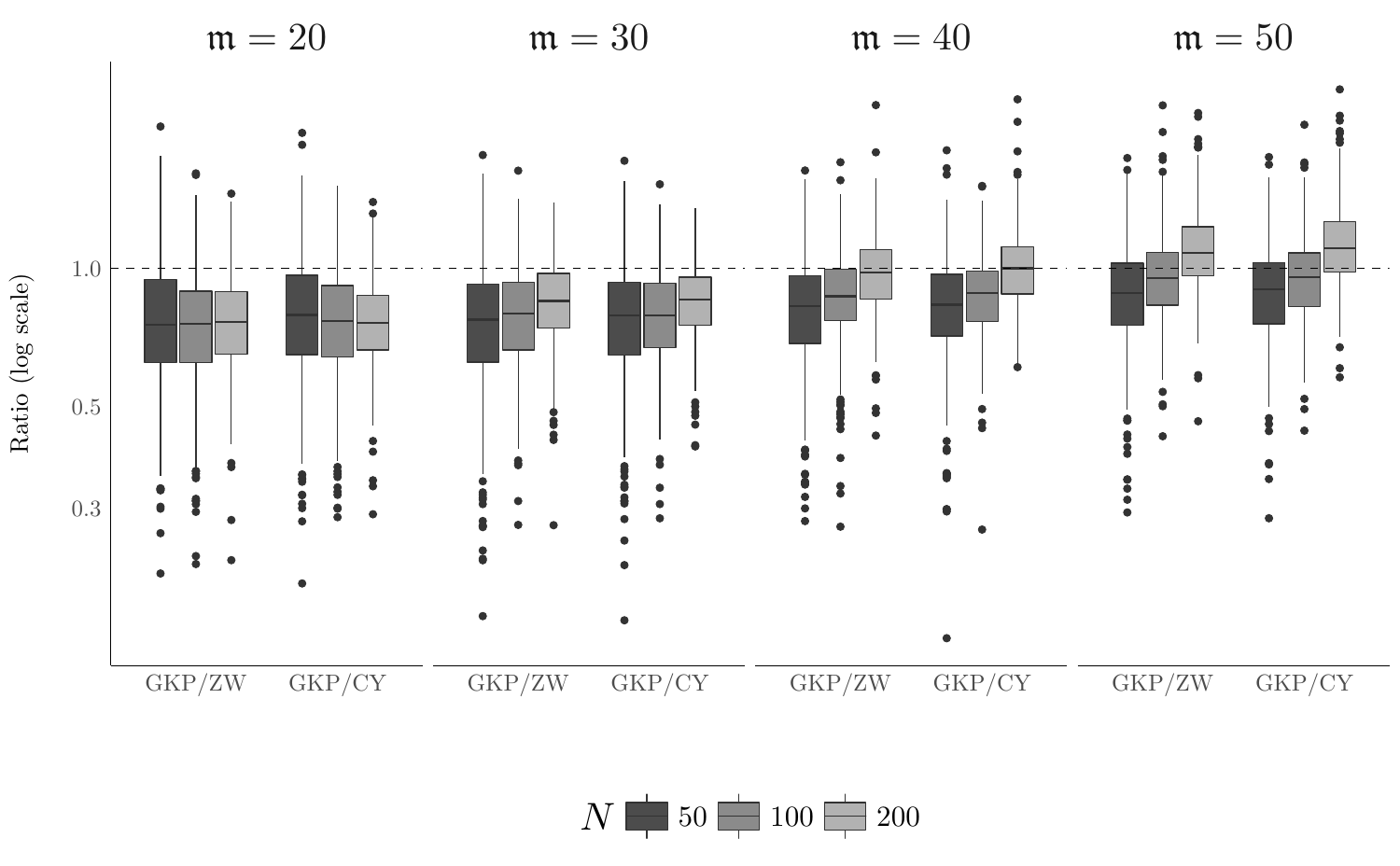}
     \caption{Results for the estimation of $\mu$ for \emph{Experiment 3} (noise std  $\sigma = 1$). The ratios are computed using $\text{ISE}_{0}$.}
     \label{fig:mean_true_exp3}
\end{figure}

\begin{figure}[h!]
	\centering
     \includegraphics[scale=0.6]{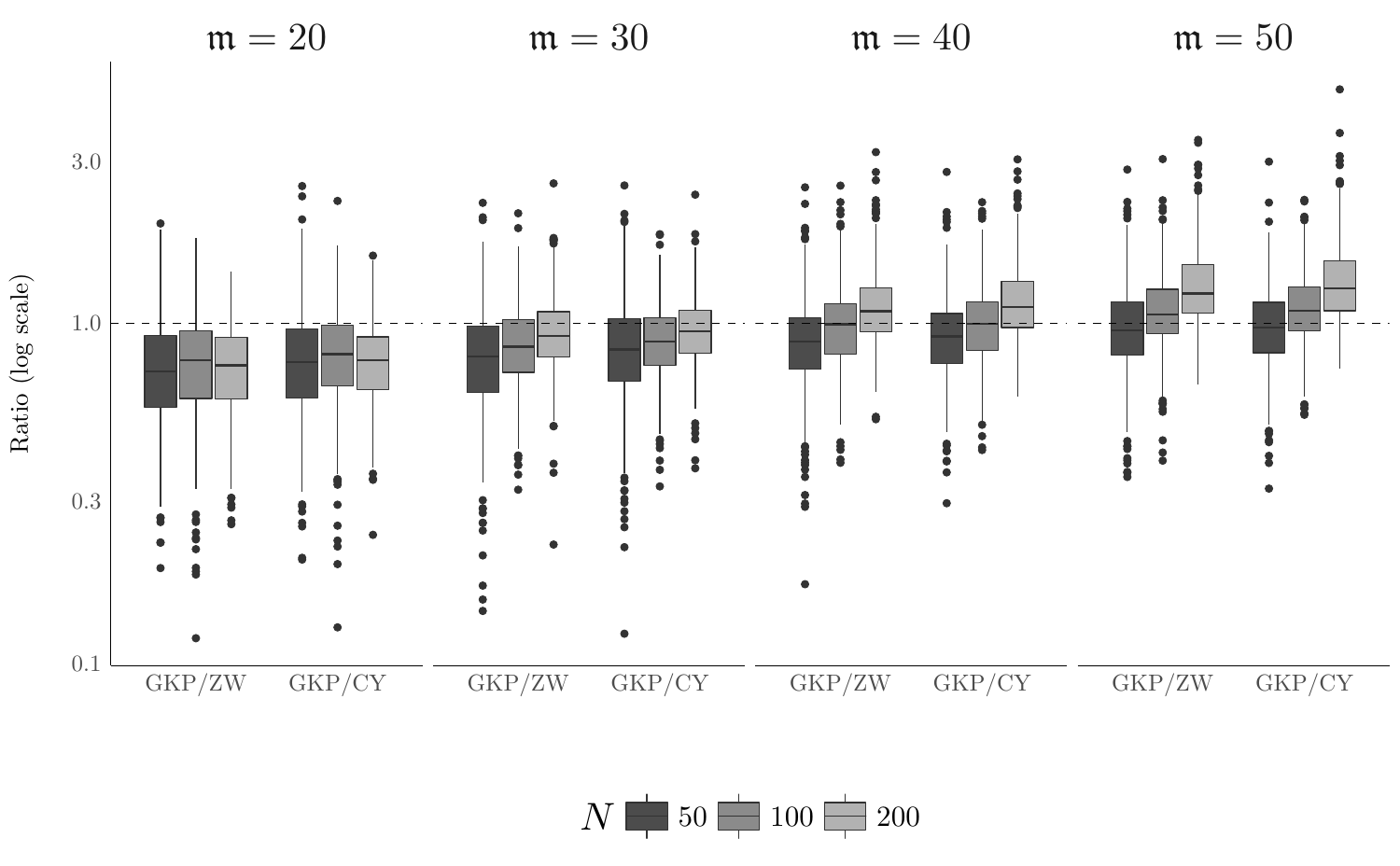}
     \caption{Results for the estimation of $\mu$ for \emph{Experiment 4} (smoother true  mean curve $\mu$). The ratios are computed using $\text{ISE}_{0}$.}
     \label{fig:mean_true_exp4}
\end{figure}

\begin{figure}[h!]
	\centering
     \includegraphics[scale=0.6]{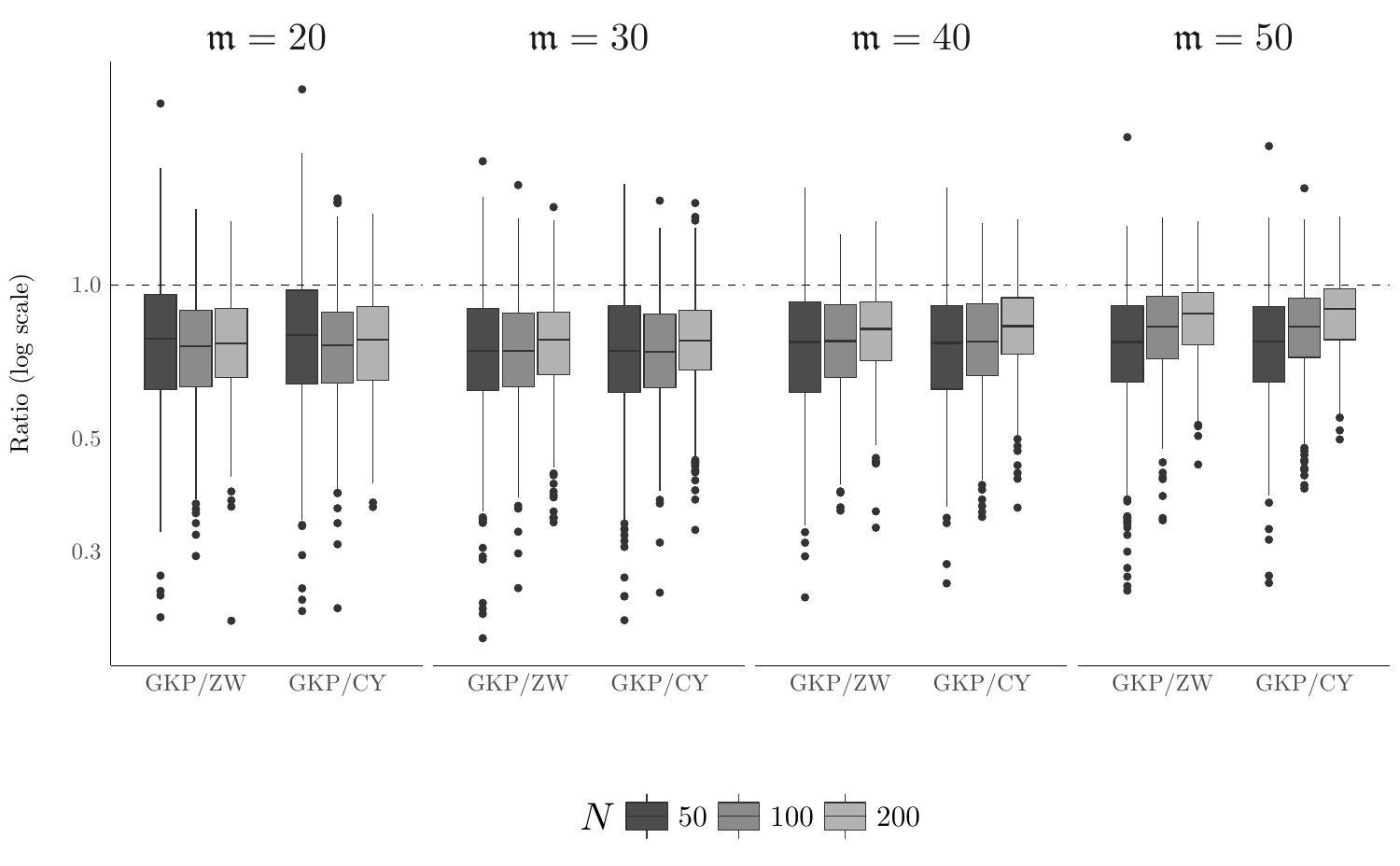}
     \caption{Results for the estimation of $\mu$ for \emph{Experiment 5} (smoother maps $H$ and $L$). The ratios are computed using $\text{ISE}_{0}$.}
     \label{fig:mean_true_exp5}
\end{figure}

\begin{figure}[h!]
	\centering
     \includegraphics[scale=0.6]{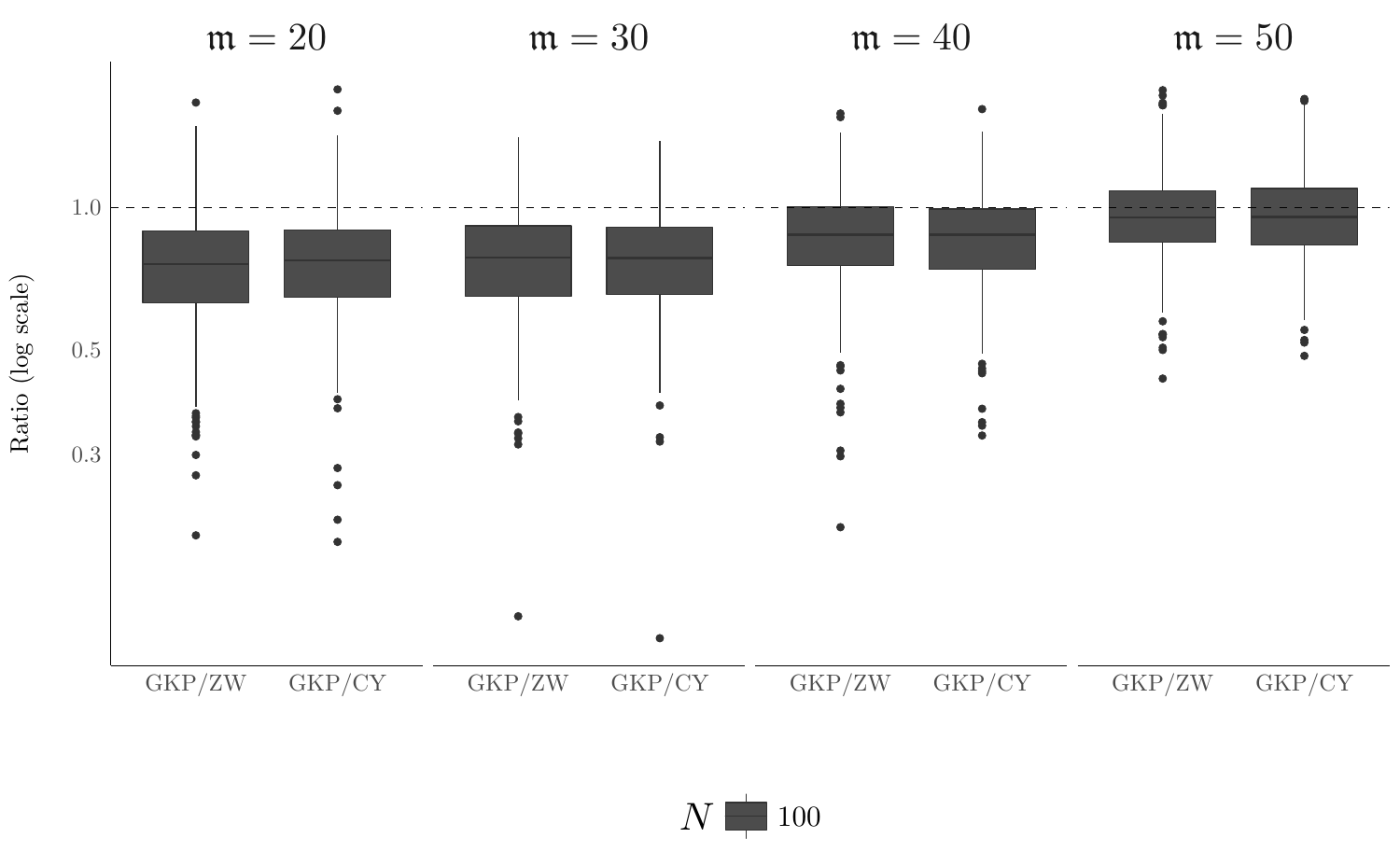}
     \caption{Results for the estimation of $\mu$ for \emph{Experiment 6} (the density of $\Tnm$ is a beta mixture). The ratios are computed using $\text{ISE}_{0}$.}
     \label{fig:mean_true_exp6} 
\end{figure}

\begin{figure}[h!]
	\centering
     \includegraphics[scale=0.6]{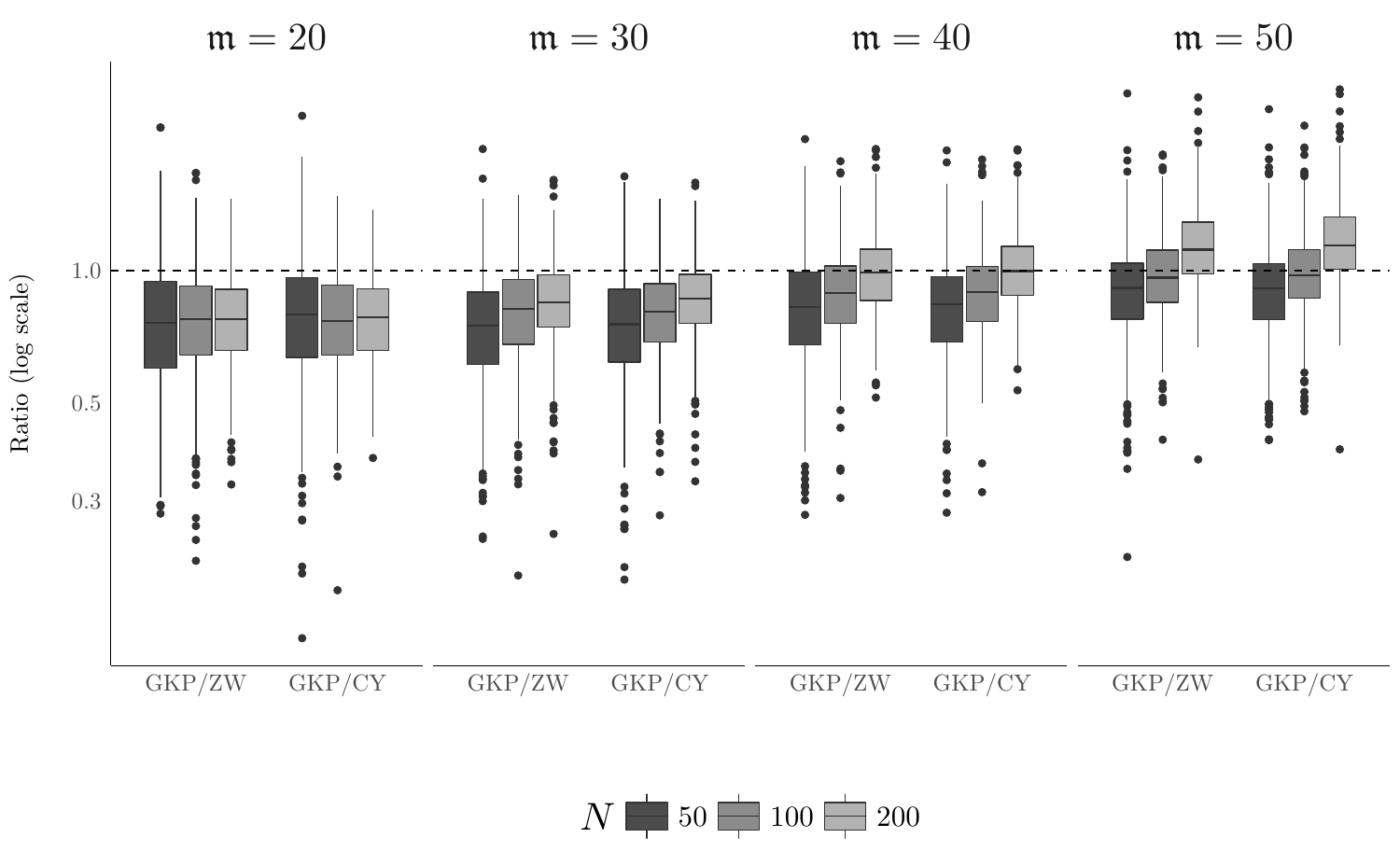}
     \caption{Results for the estimation of $\mu$ for \emph{Experiment 7} (std of $X(0)$ is $\varpi=1$). The ratios are computed using $\text{ISE}_{0}$.}
     \label{fig:mean_true_exp7}
\end{figure}

% !TEX root=../main_supp.tex

\subsection{Covariance estimation} % (fold)
\label{sub:covariance_estimation}

\begin{figure}[h!]
	\centering
	
	\includegraphics[scale=0.6]{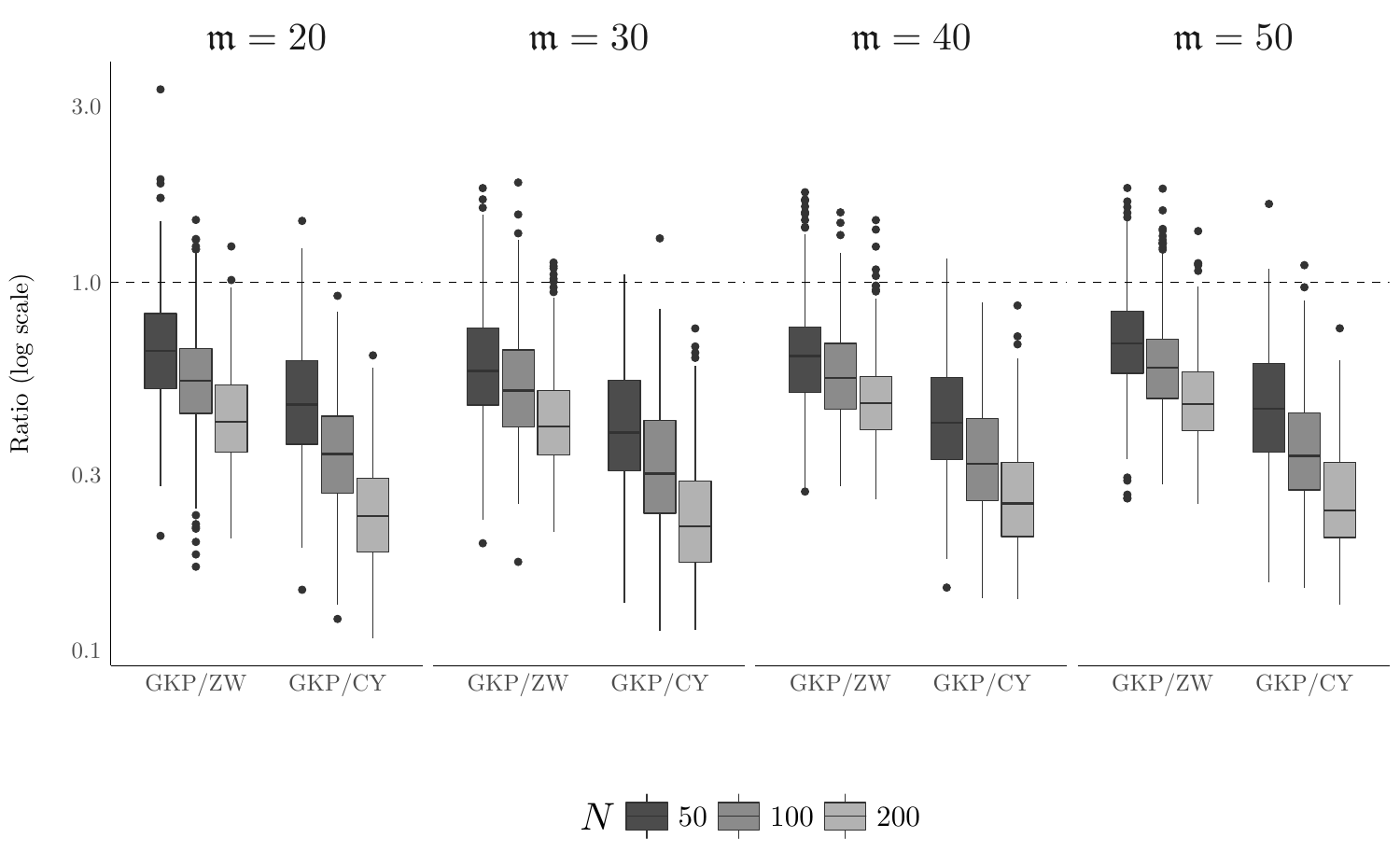}
	\caption{Results for the estimation of $\Gamma$ for \emph{Experiment 1}. The ratios are computed using $\text{ISE}_{0.05}$.}
	\label{fig:cov_true_exp1_supp}
\end{figure}

\begin{figure}[h!]
	\centering
	\includegraphics[scale=0.6]{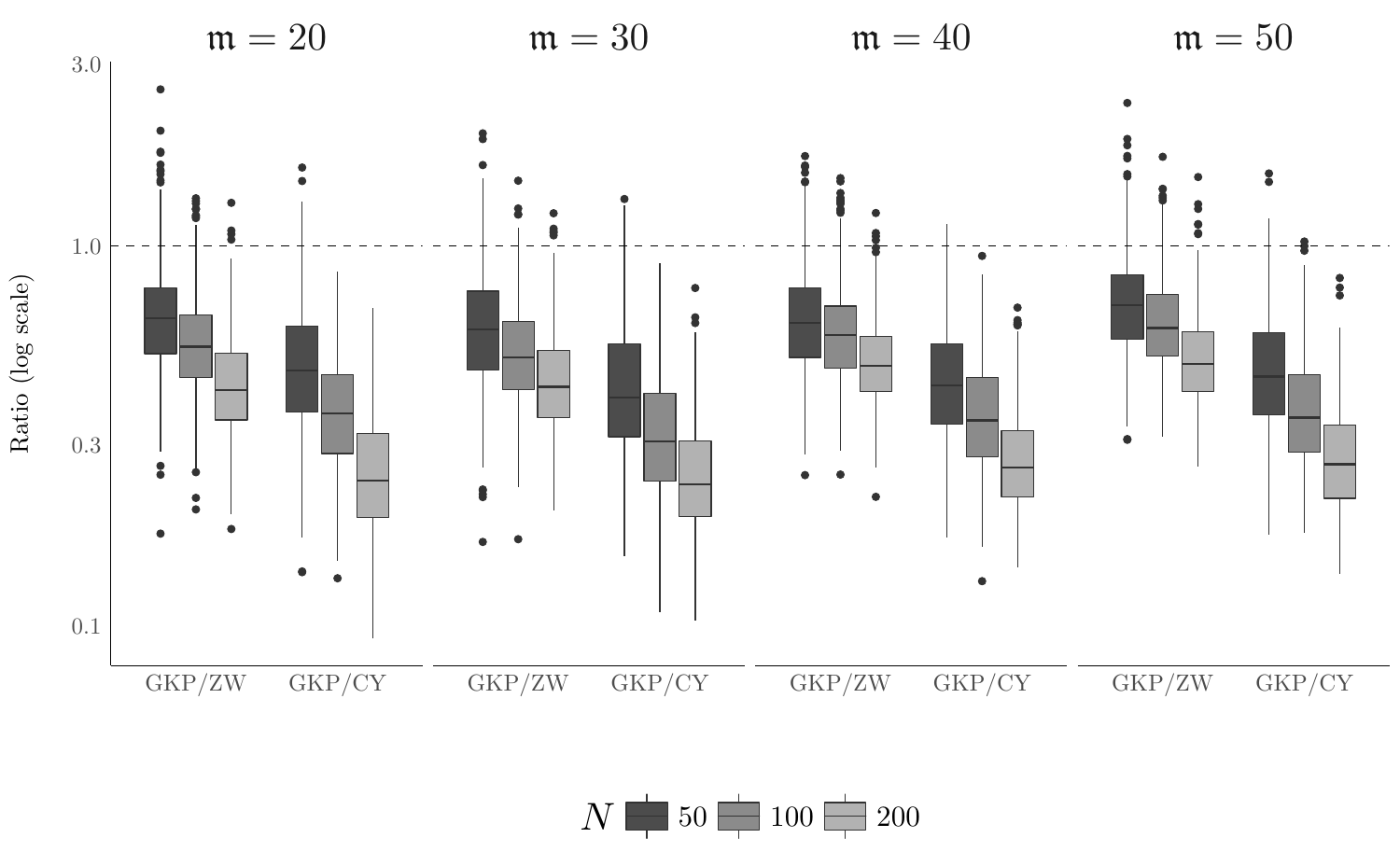}
	\caption{Results for the estimation of $\Gamma$ for \emph{Experiment 2} (noise std  $\sigma = 0.25$). The ratios are computed using $\text{ISE}_{0}$.}
	\label{fig:cov_true_exp2}
\end{figure}

\begin{figure}[h!]
	\centering
	\includegraphics[scale=0.6]{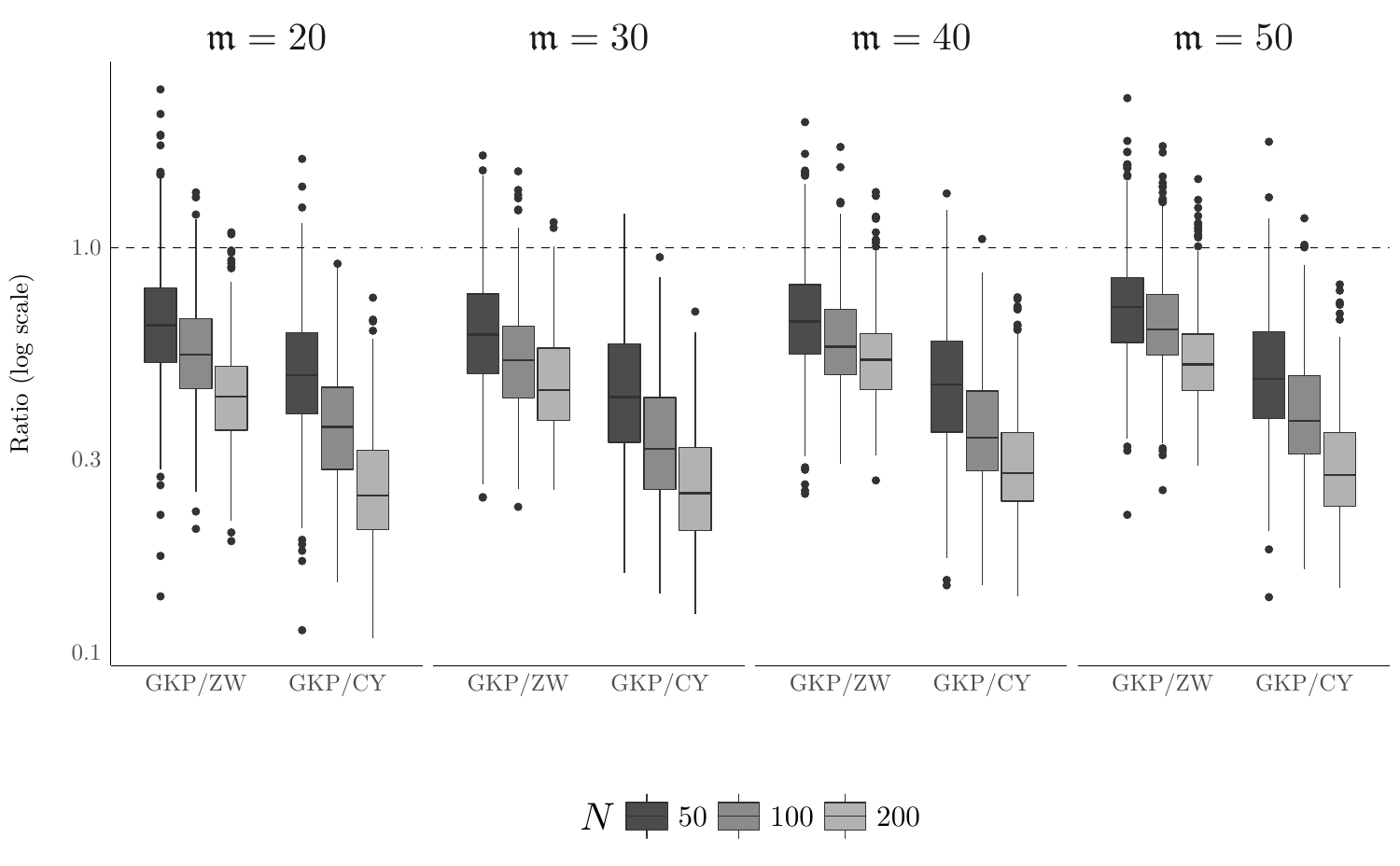}
	\caption{Results for the estimation of $\Gamma$ for \emph{Experiment 3} (noise std  $\sigma = 1$). The ratios are computed using $\text{ISE}_{0}$.}
	\label{fig:cov_true_exp3}
\end{figure}

%%%%%%%%%%%%%%%%%%

\begin{figure}[h!]
	\centering
	\includegraphics[scale=0.6]{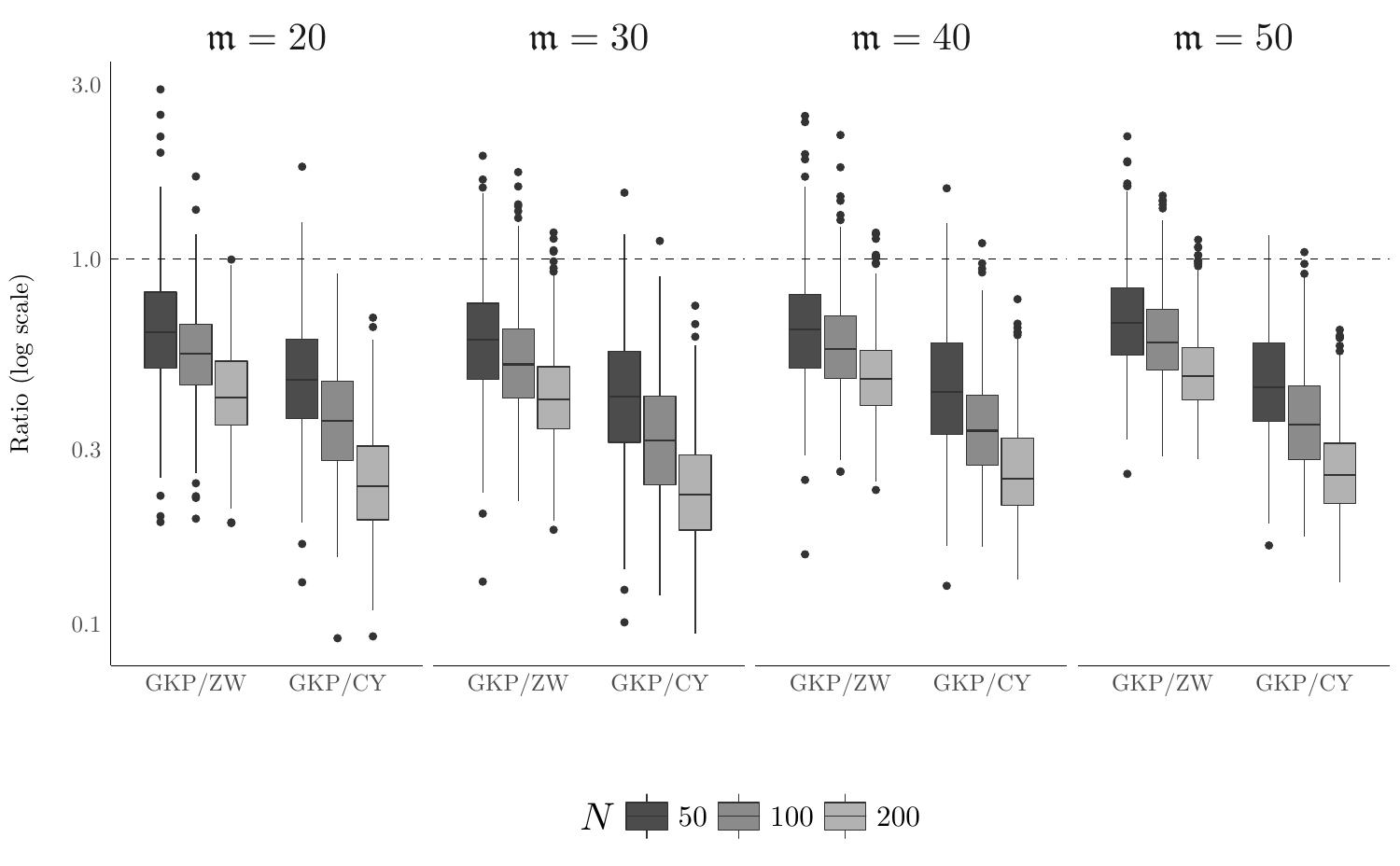}
	\caption{Results for the estimation of $\Gamma$ for \emph{Experiment 4} (smoother true mean $\mu$). The ratios are computed using $\text{ISE}_{0}$.}
	\label{fig:cov_true_exp4}
\end{figure}

\begin{figure}[h!]
	\centering
	\includegraphics[scale=0.6]{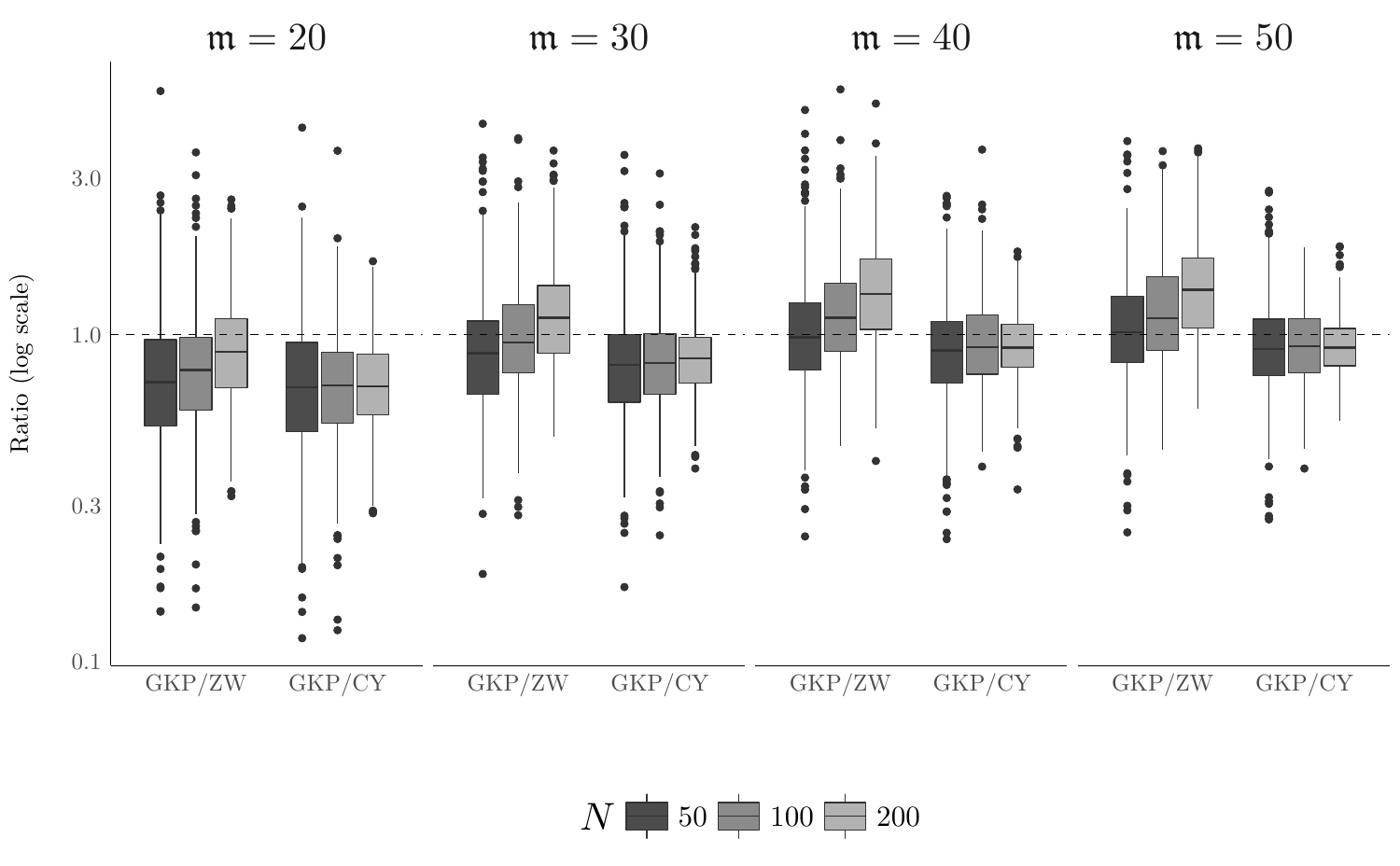}
	\caption{Results for the estimation of $\Gamma$ for \emph{Experiment 5} (smoother maps $H$ and $L$). The ratios are computed using $\text{ISE}_{0}$.}
	\label{fig:cov_true_exp5}
\end{figure}

%%%%%%%%%%%%%%%%%
\begin{figure}[h!]
	\centering
	\includegraphics[scale=0.6]{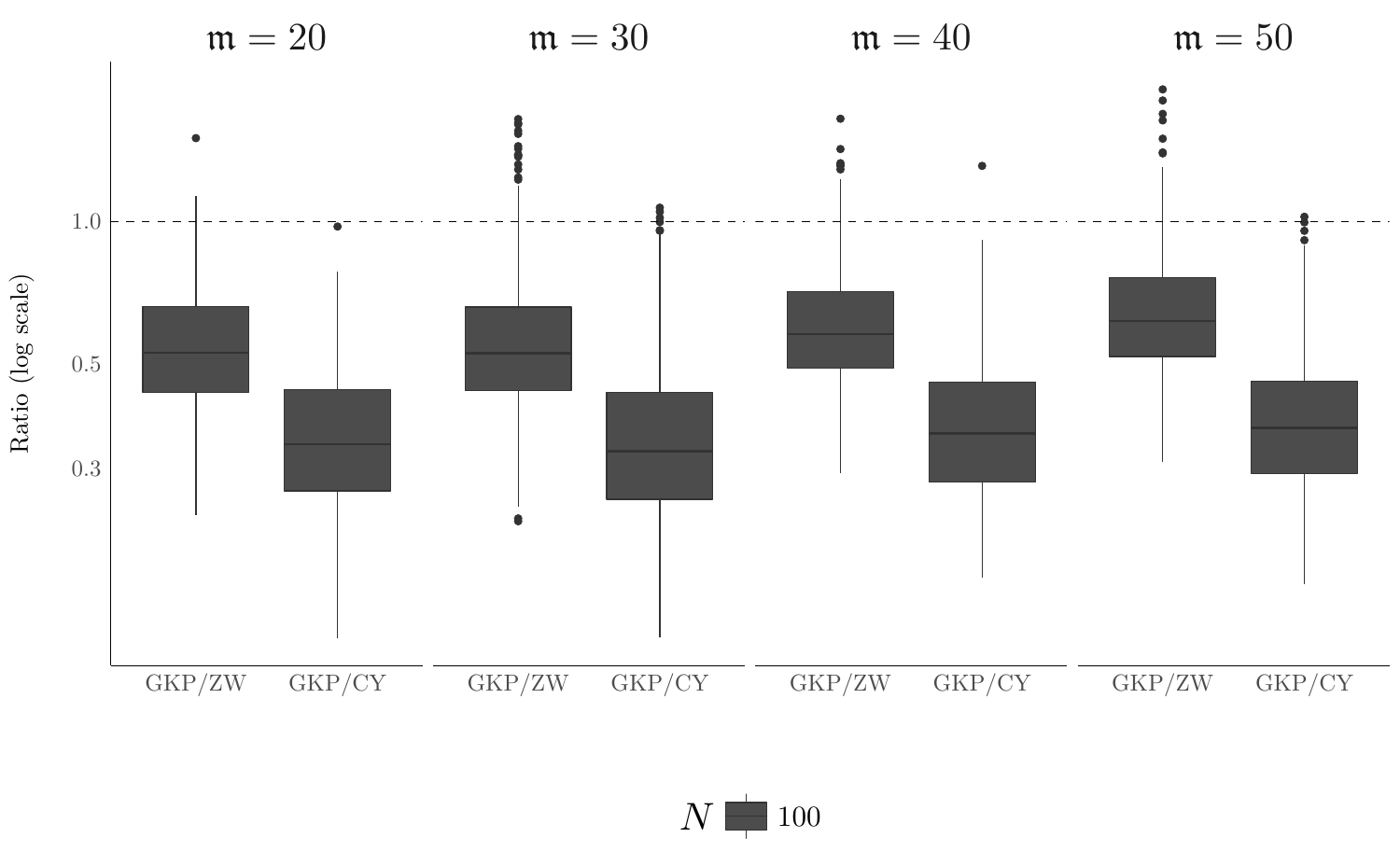}
	\caption{Results for the estimation of $\Gamma$ for \emph{Experiment 6} (the density of the $\Tnm$ is a beta mixture). The ratios are computed using $\text{ISE}_{0}$.}
	\label{fig:cov_true_exp6}
\end{figure}

\begin{figure}[h!]
	\centering
	\includegraphics[scale=0.6]{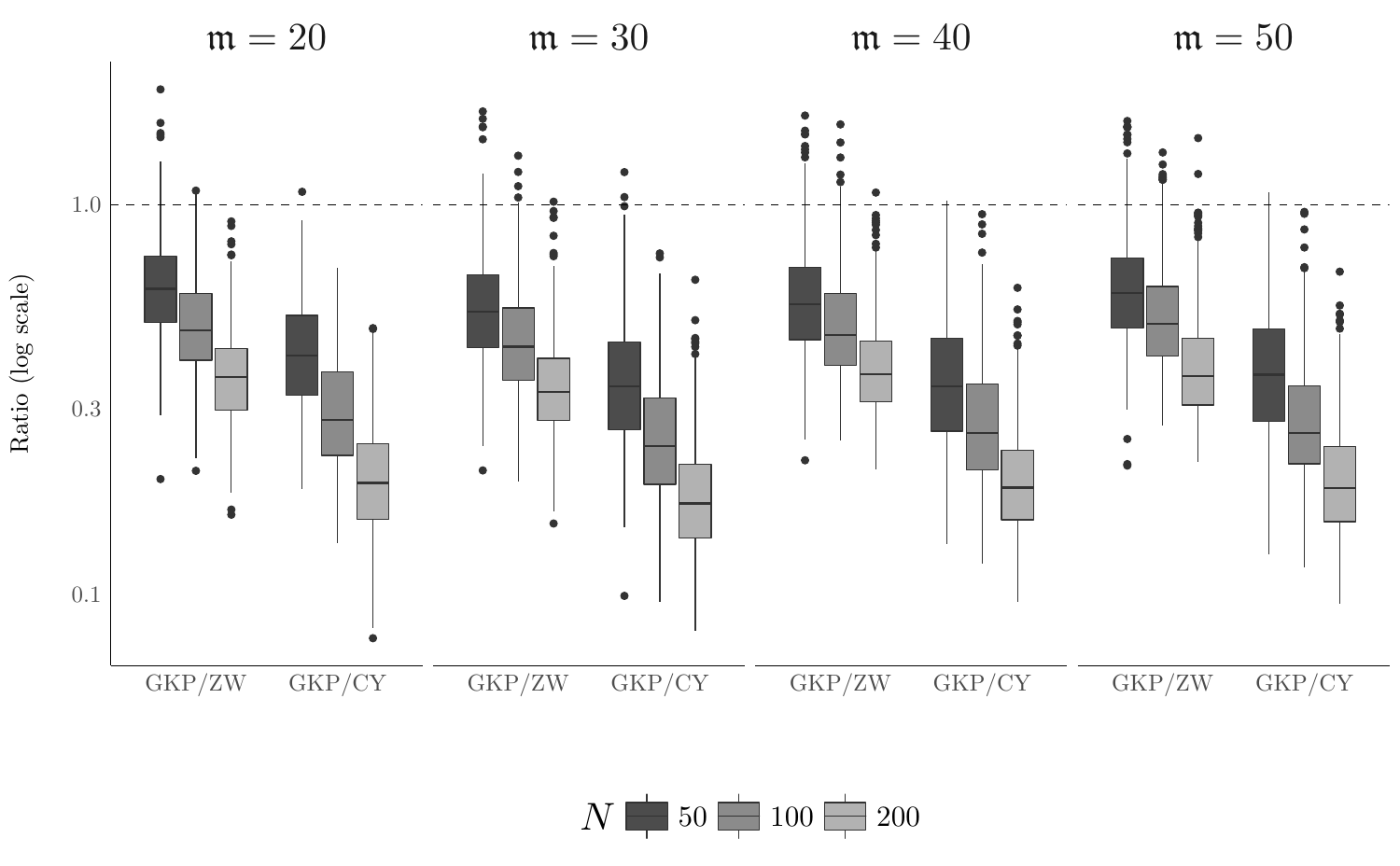}
	\caption{Results for the estimation of $\Gamma$ for \emph{Experiment 7} (std of $X(0)$ is $\varpi =1$). The ratios are computed using $\text{ISE}_{0}$.}
	\label{fig:cov_true_exp7}
\end{figure}

% !TEX root=../main_supp.tex

\subsection{Case of differentiable curves} % (fold)
\label{sub:case_of_differentiable_curves}

Let us note that, for any $d\geq 1$,  we can use  $X$ as  in \eqref{eq:X_md} to define a process which, almost surely, has $d-$times differentiable sample paths and the derivatives of order $d$ satisfy \ref{H:equivalent}. Indeed, it suffices to define
$$
X(t) = \int_0^{t}\int_0^{s_1}\cdots \int_0^{s_{d-1}} X(s_{d}) \mathrm{d}s_d \cdots \mathrm{d}s_2 \mathrm{d}s_1, \quad t\geq 0.
$$

\begin{figure}[h!]
	\centering
	\begin{subfigure}[b]{0.48\textwidth}
		\centering
		\includegraphics[width=\textwidth]{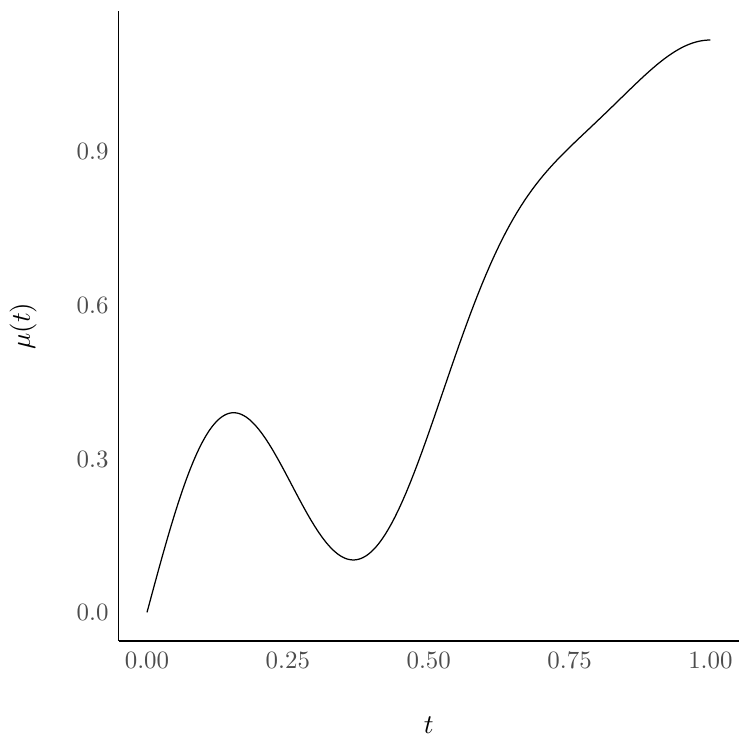}
		\caption{Mean curve $\mu(\cdot)$.}
		\label{fig:mean_curve_deriv}
	\end{subfigure}
	\hfill
	\begin{subfigure}[b]{0.48\textwidth}
		\centering
		\includegraphics[width=\textwidth]{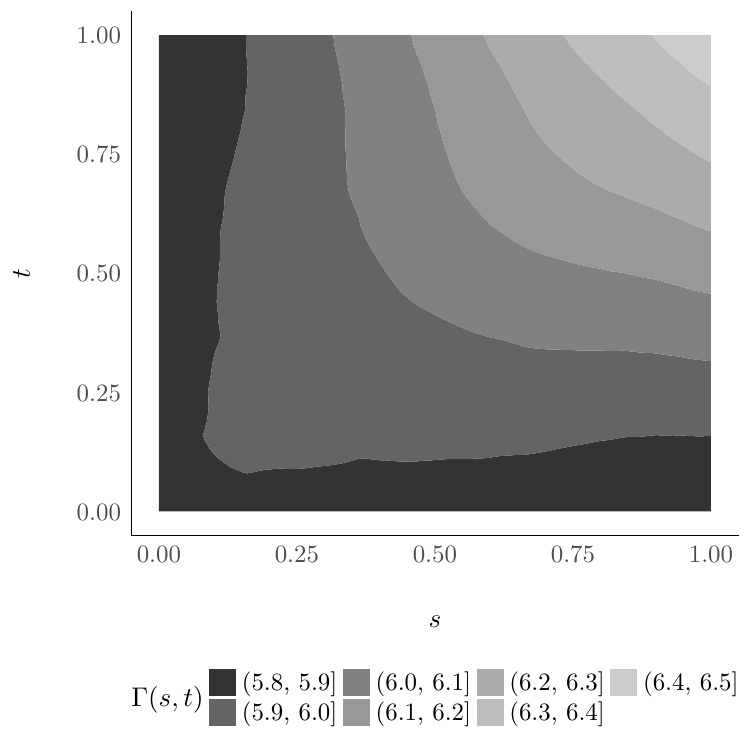}
		\caption{Covariance surface $\Gamma(\cdot, \cdot)$.}
		\label{fig:cov_curve_deriv}
	\end{subfigure}
	\\
	\begin{subfigure}[b]{0.48\textwidth}
		\centering
		\includegraphics[width=\textwidth]{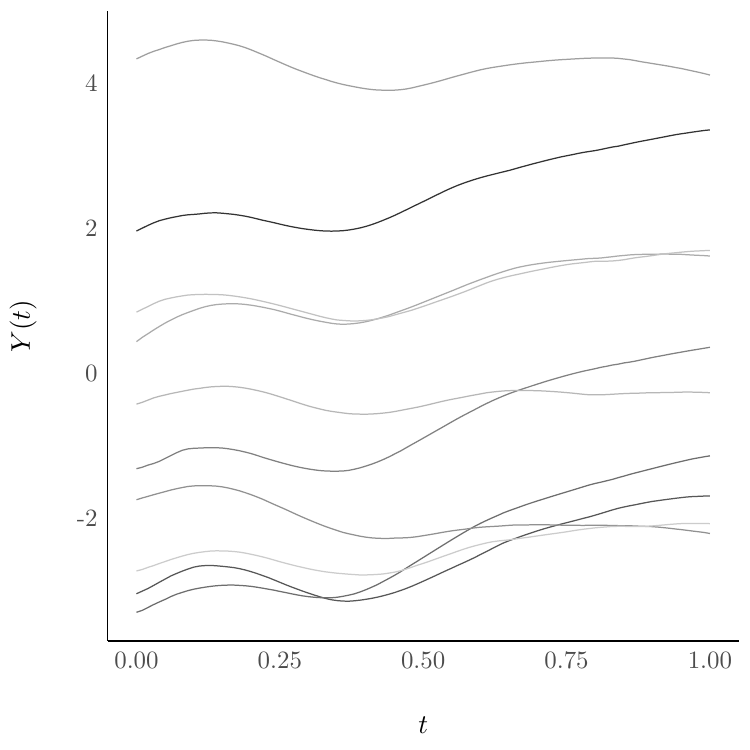}
		\caption{Curves $X^{(i)}$.}
		\label{fig:sample_curve_deriv}
	\end{subfigure}
	\hfill
	\begin{subfigure}[b]{0.48\textwidth}
		\centering
		\includegraphics[width=\textwidth]{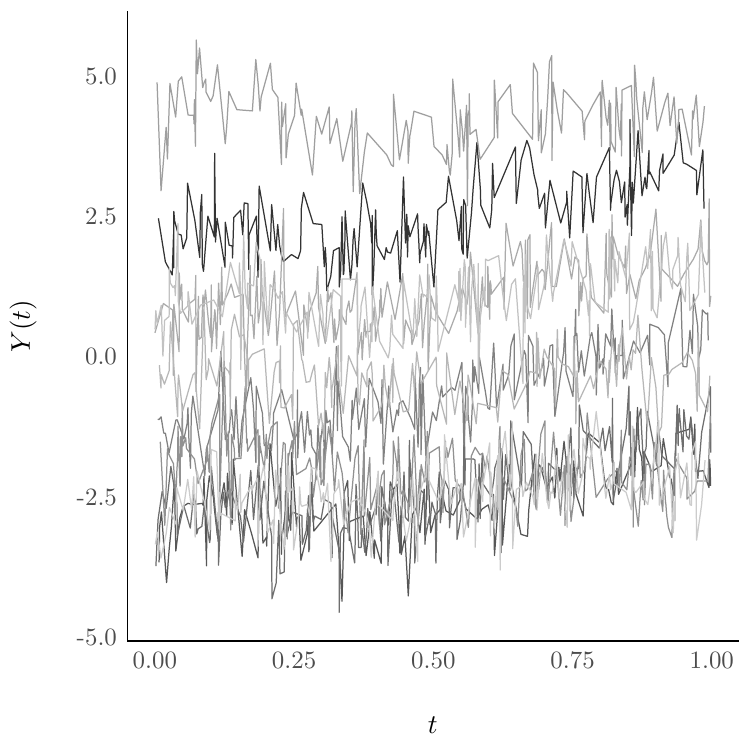}
		\caption{Noisy curves $Y^{(i)}$.}
		\label{fig:sample_noisy_curve_deriv}
	\end{subfigure}
	\caption{Description of the simulated data set with differentiable curves.}
	\label{fig:dataset_description_deriv}
\end{figure}

We consider the case of the estimation of the mean function for differentiable curves ($d=1$), referred to as \emph{Experiment 8}. More precisely, we generate curves as in \emph{Experiment 1} and perform numerical integration such that the regularity of the curves is larger than one, and the Hurst index function $H_t$ is defined on the sample path of the first derivative, for all $t \in [0, 1]$. See also \cite{golo2020} for the formal definition of the local regularity for the case of differentiable sample paths. In this experiment, the mean curve is not learned from the Power Consumption data set but generated as follows:
\begin{equation}
	\mu(t) = \sqrt{2}\sum_{k = 1}^{5} z_k\frac{\sin((k - 1/2)\pi t)}{(k - 1/2)\pi}, \quad(z_1, \dots, z_5) = (1.37, -0.56,  0.36,  0.63,  0.40).
\end{equation}
The values $z_k$ are obtained as random draws of a standard Gaussian distribution.  

We plot the mean curve $\mu(\cdot)$ in Figure \ref{fig:mean_curve_deriv} and the covariance matrix $\Gamma(\cdot, \cdot)$ in Figure \ref{fig:cov_curve_deriv}. A random sample of curves generated according to our simulation setup is plotted in Figure \ref{fig:sample_curve_deriv} without noise and in Figure \ref{fig:sample_noisy_curve_deriv} with noise.

\begin{figure}[h!]
	\centering
	\includegraphics[scale=0.6]{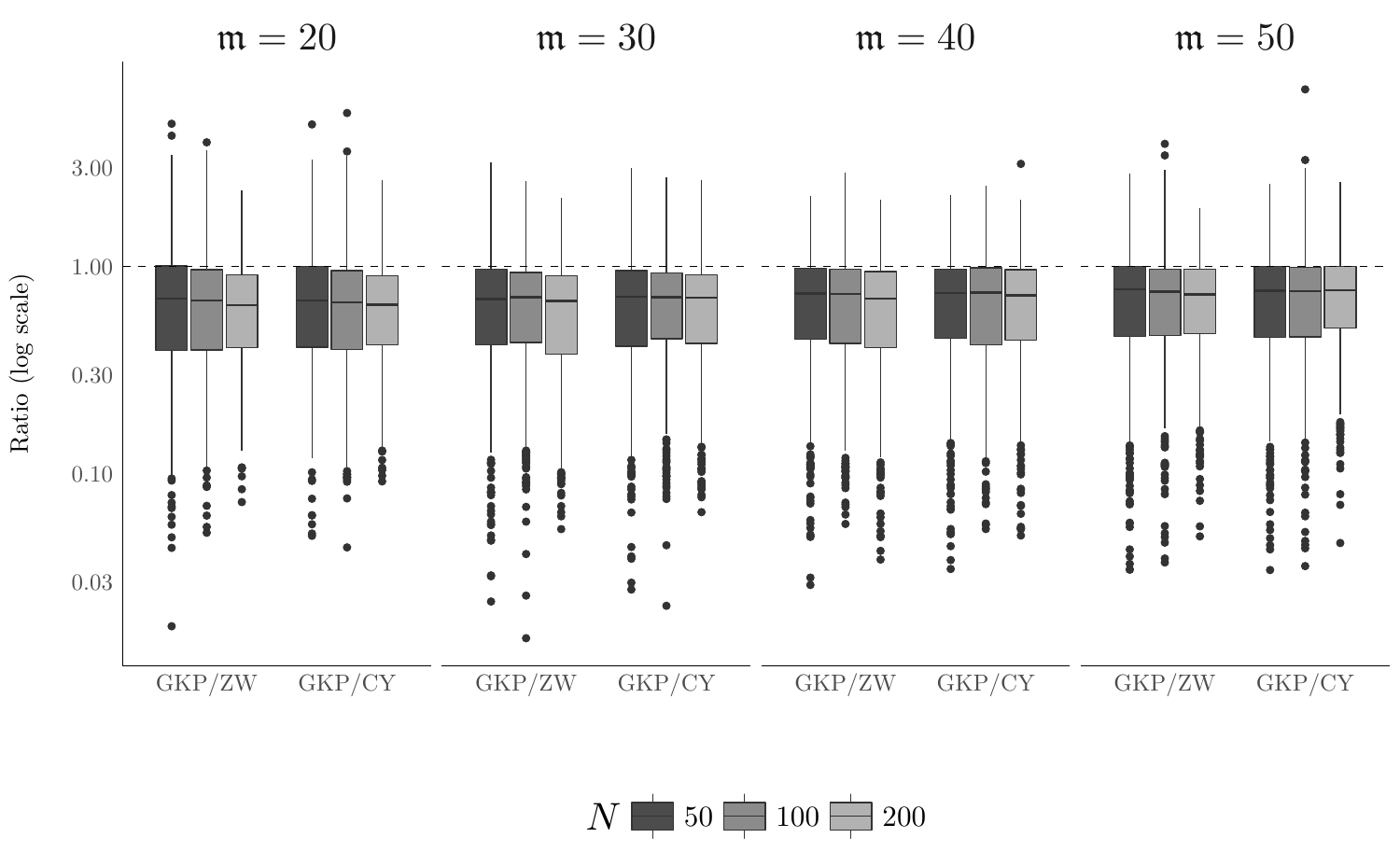}
	\caption{Results for the estimation of $\mu$ for \emph{Experiment 8}. The ratios are computed using $\text{ISE}_{0}$.}
	\label{fig:mean_true_exp8}
\end{figure}

As we assumed that the curves are differentiable, we first estimate their derivatives using local polynomials of degree $2$ with bandwidth $3 / \mathfrak{m}$. The estimation of the Hurst index function $\widehat{H}_t$ is then performed on the set of estimated derivative curves. Finally, our bandwidth selection methodology is run with  $q^2_1h^{2(1 + \widehat{H}_t)}$ as the first term in the definition of $\mathcal{R}_\mu (t;h)$ in \eqref{eq:opt_h_mean0b}. The results are plotted in Figure~\ref{fig:mean_true_exp8}, on a logarithmic scale. The ratios are obtained using $\text{ISE}_0$. Our estimator outperforms the competitors for every pair $(N, \mathfrak{m})$.

\newpage
%% if your bibliography is in bibtex format, uncomment commands:
\bibliographystyle{apalike}
\bibliography{ref_short}

%% or include bibliography directly:
% \begin{thebibliography}{}
% \bibitem[\protect\citeauthoryear{???}{???}]{b1}
% \end{thebibliography}

\makeatletter\@input{main_ref.tex}\makeatother